\documentclass[leqno,11pt,a4paper]{amsart}
\pdfoutput=1
\usepackage[full]{textcomp}
\usepackage{newpxtext}
\usepackage{comment}
\usepackage{cabin} 
\usepackage[varqu,varl]{inconsolata} 
\usepackage[bigdelims,vvarbb]{newpxmath} 
\usepackage[cal=cm,bb=esstix,bbscaled=1.05]{mathalfa} 
\usepackage[margin=2.4cm]{geometry}
\usepackage{savesym}
\let\savedegree\bigtimes
\let\bigtimes\relax
\usepackage{mathabx}
\usepackage{mathrsfs}
\usepackage{amsmath}
\let\bigtimes\savedegree
\usepackage[usenames,dvipsnames]{color}
\usepackage{hyperref}
\hypersetup{
 colorlinks,
 linkcolor={blue!90!black},
 citecolor={red!80!black},
 urlcolor={blue!50!black}
}
\usepackage{tikz}
\usepackage{tikz-cd}
\usepackage{cite}
\usepackage{enumitem}
\usepackage{amsfonts}

\setlist[enumerate]{labelsep=*, leftmargin=1.5pc}
\setlist[enumerate]{label=\normalfont(\roman*), ref=\roman*}
\newtheorem{thm}{Theorem}[section]
\newtheorem{lemma}[thm]{Lemma}
\newtheorem{cor}[thm]{Corollary}
\newtheorem{prop}[thm]{Proposition}

\newtheorem{innercustomthm}{Theorem}
\newenvironment{customthm}[1]
  {\renewcommand\theinnercustomthm{#1}\innercustomthm}
  {\endinnercustomthm}

\newtheorem{innercustomcor}{Corollary}
\newenvironment{customcor}[1]
  {\renewcommand\theinnercustomcor{#1}\innercustomcor}
  {\endinnercustomcor}

\theoremstyle{definition}
\newtheorem{example}[thm]{Example}
\newtheorem{notation}[thm]{Notation}
\newtheorem{remark}[thm]{Remark}
\newtheorem{definition}[thm]{Definition}
\newtheorem*{definition*}{Definition}
\newtheorem*{notation*}{Notation}
\newtheorem*{organization*}{Organization}
\newtheorem*{theorem*}{Theorem}
\numberwithin{equation}{section}

\newcommand{\kk}{\mathbf{k}}
\newcommand{\aaa}{\mathbf{a}}
\newcommand{\bb}{\mathbf{b}}
\newcommand{\Gr}{\mathbf{Gr}}
\newcommand{\HH}{\mathcal{H}}
\newcommand{\XX}{\mathcal{X}}
\newcommand{\EE}{\mathcal{E}}
\newcommand{\PP}{\mathbf{P}}
\newcommand{\inn}{\mathrm{in}}

\newcommand{\TT}{\mathbf{T}}

\newcommand{\bg}{\boldsymbol{\gamma}}
\newcommand{\tZ}{\tilde{Z}}
\newcommand{\hZ}{\hat{Z}}

\newcommand{\sss}{\mathfrak{S}}
\newcommand{\QQ}{\mathbf{Q}}

\newcommand{\ZZ}{\mathbf{Z}}

\DeclareMathOperator{\GL}{GL}
\DeclareMathOperator{\Cl}{Cl}
\DeclareMathOperator{\spn}{span}
\DeclareMathOperator{\Hom}{Hom}
\DeclareMathOperator{\Hilb}{Hilb}
\DeclareMathOperator{\Proj}{Proj}
\DeclareMathOperator{\Spec}{Spec}
\DeclareMathOperator{\sat}{sat}
\DeclareMathOperator{\Bl}{Bl}
\DeclareMathOperator{\Sym}{Sym}

\DeclareMathOperator{\Pic}{Pic}
\DeclareMathOperator{\Nef}{Nef}
\DeclareMathOperator{\Eff}{Eff}

\makeatletter
\@namedef{subjclassname@2020}{%
  \textup{2020} Mathematics Subject Classification}
\makeatother
\begin{document}
\author[R.\,Ramkumar]{Ritvik~Ramkumar}
\address{Department of Mathematics\\University of California at Berkeley\\Berkeley, CA\\94720\\USA}
\email{ritvik@math.berkeley.edu}
\title{The Hilbert scheme of a pair of linear spaces}
\maketitle
\begin{abstract} Let $\HH_{a,b}^n$ denote the component of the Hilbert scheme whose general point parameterizes an $a$-plane union a $b$-plane meeting transversely in $\PP^n$. We show that $\HH_{a,b}^n$ is smooth and isomorphic to successive blow ups of $\Gr(a,n) \times \Gr(b,n)$ or $\text{Sym}^2 \Gr(a,n)$ along certain incidence correspondences. We classify the subschemes parameterized by $\HH_{a,b}^n$ and show that this component has a unique Borel fixed point. We also study the birational geometry of this component. In particular, we describe the effective and nef cones of $\HH_{a,b}^n$ and determine when the component is Fano. Moreover, we show that $\HH_{a,b}^n$ is a Mori dream space for all values of $a,b,n$. 
\end{abstract}

\setcounter{section}{-1}
\section{Introduction}

The Hilbert scheme $\Hilb^{P(t)} \mathbf{P}^n$, which parameterizes closed subschemes of $\mathbf{P}^n$ with a fixed Hilbert polynomial $P(t)$, introduced by Grothendieck \cite{g}, has attracted a lot of interest. Although their singularities are known to be arbitrarily complicated \cite{v}, the cases when they are smooth or have smooth components have been extensively studied. Early on these smooth components were used to solve numerous enumerative problems \cite{ess} and with major advances in the minimal model program \cite{bchm},  they are also a source of examples with rich birational structure. Fogarty \cite{f} proved that $\Hilb^{m} \mathbf{P}^2$ is smooth and Arcara, Bertram, Coskun and Huizenga \cite{abch} proved that its a Mori dream space and described the stable base decomposition of its effective cone in numerous cases. Piene and Schlessinger \cite{ps} showed that $\Hilb^{3t+1} \mathbf{P}^3$ has two smooth components that meet transversely and described the points of the component  corresponding to twisted cubics explicitly. Chen \cite{chen} proved that the component corresponding to the twisted cubics is the flip of $\widebar{\mathcal{M}}_{0,0}(\PP^3,3)$ over the Chow variety. Avritzer and Vainsencher \cite{av} proved that the component corresponding to elliptic quartics in $\Hilb^{4t} \mathbf{P}^3$ is smooth and isomorphic to a double blow up of $\mathbf{Gr}(1,9)$; Gallardo, Huerta and Schmidt \cite{ghs} computed its effective cone. Chen, Coskun and Nollet \cite{ccn} showed that the component corresponding to a pair of codimension two linear spaces meeting transversely is smooth and isomorphic to a blow of $\text{Sym}^2 \mathbf{Gr}(n-2,n)$. They also completely worked out its Mori theory. It is thus very interesting to find components of Hilbert schemes that are smooth and describe their birational geometry.

Let $\kk$ be an algebraically closed field with $\mathrm{char} \, \kk \ne 2$ and let $d \geq c \geq 2$. Let $X$ be the union of an $(n-c)$-dimensional plane and an $(n-d)$-dimensional plane meeting  transversely in $\mathbf{P}^n$. The Hilbert polynomial of $X$ is
\begin{equation*}
P_{n-c,n-d}^n(t) = \binom{n-c+t}{t}+\binom{n-d+t}{t}-\binom{n-c-d+t}{t}.
\end{equation*}
There is an integral component of $\Hilb^{P_{n-c,n-d}^n(t)} \, \mathbf{P}^n$, denoted $\mathcal{H}_{n-c,n-d}^n$ or $\mathcal{H}_{n-c,n-d}(\mathbf{P}^n)$, whose general point parameterizes $X$, see Proposition \ref{component}.

\smallskip

We begin with the natural rational map
\begin{equation} \label{XI}
\Xi: \mathbf{Gr}(n-c,n) \times \mathbf{Gr}(n-d,n) \dashrightarrow \HH_{n-c,n-d}^n, \quad (\Lambda,\Lambda') \mapsto [I_{\Lambda}I_{\Lambda'}].
\end{equation}
If $c=d$, the rational map is $\mathfrak{S}_2$-equivariant where $\mathfrak{S}_2$ is the group of order $2$. It acts on $\Gr(n-c,n)^2$ by interchanging the two factors and acts trivially on $\HH_{n-c,n-c}^n$.

\begin{definition} \label{gamma} For each $1\leq i \leq c$ define an incidence variety
\begin{equation*}
\Gamma_{i} = \{(\Lambda,\Lambda'): \mathrm{codim}_{\PP^n} (\Lambda \cap \Lambda') \leq d-1+i \} \subseteq \mathbf{Gr}(n-c,n) \times \mathbf{Gr}(n-d,n).
\end{equation*}
\end{definition}
Note that $\Xi$ is defined on the open subset where the two planes meet transversely. If $X$ spans $\mathbf{P}^n$ (when $n\geq c+d-1$) then this open set is precisely the complement of $\Gamma_c$. Moreover, in this case, $\Xi$ is also defined on the complement of $\Gamma_{c-1}$ (Lemma \ref{extend1}). By explicitly resolving $\Xi$ and studying the induced morphism, we obtain

\newtheorem*{thm:main}{Theorem \ref{main}}
\begin{thm:main} Let $c \geq 2$ and $n \geq 2c-1$. The component $\HH_{n-c,n-c}^n$ is smooth and the map $\Xi$ induces an isomorphism
$$
\Bl_{\widebar{\Gamma}_{c-1}}\cdots\Bl_{\widebar{\Gamma}_{1}}\Sym^{2}\mathbf{Gr}(n-c,n)\longrightarrow \HH_{n-c,n-c}^n
$$
where  $\widebar{\Gamma}_i$  is the strict transform of $\Gamma_i/\mathfrak{S}_2$. 

If $n < 2c-1$, the morphism $\mathcal{H}_{n-c,n-c}^n \longrightarrow \mathbf{Gr}(2n-2c+1,n)$ that sends a scheme to its linear span is smooth; the fiber over a point $\Lambda$ is $\mathcal{H}_{n-c,n-c}(\Lambda)$. 
\end{thm:main}

\newtheorem*{thm:maintwo}{Theorem \ref{maintwo}}
\begin{thm:maintwo} 
Let $d > c \geq 2$ and $n \geq c+d-1$. The component $\mathcal{H}_{n-c,n-d}^n$ is smooth and $\Xi$ extends to an isomorphism
\begin{equation*}
\Xi: \Bl_{\Gamma_{c-1}} \cdots \Bl_{\Gamma_1} (\mathbf{Gr}(n-c,n) \times \mathbf{Gr}(n-d,n)) \longrightarrow \mathcal{H}_{n-c,n-d}^n. \footnote{By abuse of notation, we use $\Gamma_i$  to also denote the strict transform of $\Gamma_i$.
}
\end{equation*}
If $n < c+d -1$, the morphism $\mathcal{H}_{n-c,n-d}^n \longrightarrow \mathbf{Gr}(2n-c-d+1,n)$ that sends a scheme to its linear span is smooth; the fiber over a point $\Lambda$ is $\mathcal{H}_{n-c,n-d}(\Lambda)$. 
\end{thm:maintwo}

Historically, Harris \cite{joe} suggested that $\HH_{1,1}^3 \simeq \Bl_{\widebar{\Gamma}_{1}}\Sym^{2}\mathbf{Gr}(1,3)$ and that $\Hilb^{2t+2} \, \mathbf{P}^3$ is the union of $\HH_{1,1}^3$ and another smooth component meeting transversely. The authors of \cite{ccn} generalized this and proved that $\HH_{n-2,n-2}^n \simeq \Bl_{\widebar{\Gamma}_{1}}\Sym^{2}\mathbf{Gr}(n-2,n)$ is smooth and meets exactly one other component in $\Hilb^{P_{n-2,n-2}^n(t)}  \, \mathbf{P}^n$. 
A major step in the proof of these statements was a computation of an analytic neighbourhood of a point in the intersection of the two components using the tangent-obstruction theory for the Hilbert scheme \cite[Proposition 2.6]{ccn}. 
Unfortunately, for general $c,d$ there are many, sometimes singular, components meeting $\HH_{n-c,n-d}^n$ (Remark \ref{manycomponents}). Thus a description of a neighbourhood of a point in the intersection of all these components is most likely intractable. Our proof of Theorem \ref{main} circumvents this by using the explicit construction of $\Xi$ and studying the induced map on tangent spaces. 

In \cite{ritvik} we expounded on the philosophy that the complexity of a Hilbert scheme can be measured by their number of Borel fixed points. In line with this reasoning, we have the following result:

\newtheorem*{thm:uniqueborel}{Theorem \ref{uniqueborel2}, \ref{uniqueborel'}} 
\begin{thm:uniqueborel}  
The component $\HH_{n-c,n-d}^n$ has a unique Borel fixed point.
\end{thm:uniqueborel}

We also give a complete description of all the subschemes parameterized by $\HH_{n-c,n-d}^n$. In light of Theorem \ref{main}, \ref{maintwo} it is enough to consider the case $n \geq c+d-1 $. A \textbf{double structure} on an integral subscheme $Z \subseteq \mathbf{P}^n$ is a subscheme $Z' \subseteq \PP^n$ such that $Z'_{\text{red}} = Z$ and $\deg(Z') = 2\deg(Z)$. A double structure is said to be \textbf{pure} if it has no embedded components.

\newtheorem*{thm:primdecomp}{Theorem \ref{primdecomp}} 
\begin{thm:primdecomp} 
Let $n \geq 2c-1$. Let $Z$ be a subscheme parameterized by $\mathcal{H}_{n-c,n-c}^n$.  Then $Z$ is a pair of planes meeting transversely, or there exists a sequence of integers $1 \leq i_1 < \cdots < i_r \leq c$ and a flag of linear spaces $\Lambda^1 \subseteq \Lambda^{2} \subseteq \cdots \subseteq \Lambda^{r} \subseteq \PP^n$ with $\mathrm{codim}_{\PP^n}(\Lambda^{\ell})= (c+i_{\ell}-1)$ for each $\ell$, such that
\begin{enumerate}
\item If $i_1 >1$ then $Z$ is a union of two planes meeting along $\Lambda^1$ with embedded pure double structures on $\Lambda^{\ell}$ for each $1 \leq \ell \leq r$.
\item If $i_1 =1$ then $Z$ is a pure double structure on $\Lambda^1$ with embedded pure double structures on $\Lambda^{\ell}$ for each $2 \leq \ell \leq r$. 
\end{enumerate}
\end{thm:primdecomp}
The description when $c \ne d$ is similar and can be found in Theorem \ref{primdecomp2}. 

\newtheorem*{cor:2c}{Corollary \ref{2kpointss}, \ref{2cpointss}} 
\begin{cor:2c}  
Up to projective equivalence, there are exactly $2^c$ schemes parameterized by $\mathcal{H}_{n-c,n-d}^n$.
\end{cor:2c}

\smallskip
When $\text{char} \, \kk = 0$, we use our explicit description of $\Xi$ and the classification of ideals parameterized to study the effective and nef cones of $\HH_{n-c,n-d}^n$. As a consequence, we deduce that $\HH_{n-c,n-d}^n$ is always a Mori dream space.

\begin{definition} Let $Y$ be a smooth projective variety with $\Cl(Y)$ finitely generated. Then $Y$ is a \textbf{Mori dream space} if the Cox Ring of $Y$ is finitely generated over $\kk$. The Cox ring of $Y$ is defined to be 
$\bigoplus_{\mathbf{m} \in \ZZ^k}H^0(Y,\mathcal{O}_{Y}(\sum_i \mathbf{m}_iD_i))$
where  $D_1,\dots,D_k$ are chosen to generate $\Cl(Y)$.
\end{definition}

We also determine the pairs $(c,d)$ for which the component is Fano. For the rest of the introduction $\Lambda_m$ will always denote an $m$-dimensional linear subspace of $\PP^n$. We begin with a description of the divisors.

\begin{definition} \label{Di} Let $n \geq 2c-1$. For each $1 \leq i \leq c-1$ and a choice of a flag of linear spaces $\{\Lambda_{i-1} \subseteq \Lambda_{2c-1-i}\}$, let $D_i$ denote the divisor class of the locus of subschemes $Z \in \HH_{n-c,n-c}^n$, for which the linear span of $\Lambda_{i-1} \cup (Z\cap \Lambda_{2c-1-i})$ has dimension less than $2c-i-1$.  Let $D_c$ denote the divisor class of the locus of subschemes that meet a fixed $\Lambda_{c-1}$.
\end{definition}

\begin{definition} \label{Ni} Let $n \geq 2c-1$. Let $N_1$ denote the divisor class of the locus of generically non-reduced subschemes in $\HH_{n-c,n-c}^n$. For each $2 \leq i \leq c-1$, let $N_i$ denote the divisor class of the locus of subschemes with an embedded $(n-c+1-i)$-plane. If $n =2c-1$ let $N_c$ denote the divisor class of the locus of subschemes with an embedded point. If $n > 2c-1$ let $N_c$ denote the class of the closure of the locus of pairs of planes meeting transversely, where the intersection of the two planes meets a fixed $\Lambda_{2c-1}$.
\end{definition} 

Here are the results when $c=d$ and the pair of planes span $\PP^n$. 
 
\begin{customthm}{F} \label{theoremF}
Let $c\geq 2$ and $n \geq 2c-1$. The component $\HH_{n-c,n-c}^n$ is a Mori dream space and we have,
$$
\Eff(\HH_{n-c,n-c}^n) = \langle N_1,\dots,N_c \rangle \quad \text{and} \quad 
\Nef(\HH_{n-c,n-c}^n) = \langle D_1,\dots,D_c \rangle.
$$
Moreover, $\HH_{n-c,n-c}^n$ is Fano if and only if either $c=3$ and $n =5$, or $c \ne 3$ and $n \in \{2c-1,2c\}$. 
\end{customthm}

To state the results when the pair of planes do not span $\PP^n$, it is more convenient to use dimension instead of codimension to index the component. In particular, the component parameterizing subschemes that do not span $\PP^n$ are of the form $\HH_{c-1,d-1}^n$ with $n > c+d-1$.

\begin{definition} Let $n > 2c-1$. For each $1 \leq i \leq c-1$ and a choice of flag $\{\Lambda_{n-2c+i} \subseteq \Lambda_{n-i}\}$, let $D_{i}'$ denote the divisor class of the locus of subschemes $Z \in \HH_{c-1,c-1}^n$, for which the linear span of $\Lambda_{n-2c+i} \cup (\Lambda_{n-i} \cap Z)$ has dimension less than $n-i$. Let $D_c'$ denote the divisor class of the locus of subschemes meeting a fixed $\Lambda_{n-c}$. Let $F$ denote the divisor class of the locus of subschemes whose linear span meets a fixed $\Lambda_{n-2c}$.
\end{definition}

\begin{definition}  Let $n > 2c-1$. Let $N_1'$ denote the divisor class of the locus of generically non-reduced subschemes in $\HH_{c-1,c-1}^n$. For each $2 \leq i \leq c$, let $N_i'$ denote the divisor class of the locus of subschemes with an embedded $(c-i)$-plane.
\end{definition}

Here are the results when $c=d$ and the pair of planes do not span $\PP^n$.

\begin{customthm}{G} \label{theoremG}
Let $c \geq 2$ and $n > 2c-1$. The component $\HH_{c-1,c-1}^n$ is Fano and thus a Mori dream space. Moreover we have,
$$
\Eff(\HH_{c-1,c-1}^n) = \langle N_1',\dots,N_c',F \rangle \quad \text{and} \quad 
\Nef(\HH_{c-1,c-1}^n)  = \langle D_1',\dots,D_c',F \rangle.
$$
\end{customthm}

The precise results when $c \ne d$ can be found in Section \ref{last bro}. We conclude the introduction by describing the components that are Fano in this case; the results mirror the case of $c=d$.

\begin{customthm}{H}  \label{theoremH}
The component $\HH_{c-1,d-1}^n$ is Fano. The component $\HH_{n-c,n-d}^n$ is Fano if and only if  either $c=2$ and $n \in \{d+1,\dots,2d-1\}$, or $c \geq 3$ and $n \in \{c+d-1,c+d\}$.
\end{customthm}

\begin{organization*} In Section \ref{one} we construct the component $\HH_{n-c,n-d}^n$ and show that the rational map $\Xi$ is defined away from $\Gamma_{c-1}$. In Section \ref{two} we thoroughly study the case $c=d$.  We begin by explicitly constructing a morphism, also denoted $\Xi$, from a sequence of blowups to $\HH_{n-c,n-c}^n$ (Proposition \ref{mainone}, Proposition \ref{mainonee}). We then construct a Gr\"{o}bner basis for ideals in the image of $\Xi$ (Lemma \ref{GROBNER}), which is indispensable in showing $\Xi$ is bijective and proving Theorems \ref{uniqueborel2}, \ref{primdecomp}. By analyzing the differential of $\Xi$ at the Borel fixed point we deduce Theorem \ref{main}. In Section \ref{three} we explain how to carry out all of the proofs of Section \ref{two} with little to no modification for the case $c \ne d$. In Section \ref{Section Divisors} we study the divisors on $\HH_{n-c,n-c}^n$ and provide local equations for them. In Sections \ref{Section Part I}, \ref{Section Part II} we study the birational geometry of $\HH_{n-c,n-c}^n, \HH_{c-1,c-1}^n$ and prove Theorem \ref{theoremF} and Theorem \ref{theoremG}. More precisely, the cones are computed in Proposition \ref{cones} and Proposition \ref{lincomb2}. The fact that the components are Mori dream spaces is established in Theorem \ref{MORI}. In Section \ref{last bro} we explain how to carry out all of the proofs of Section \ref{Section Part I} and \ref{Section Part II} for the case $c \ne d$.
\end{organization*}


\section{Preliminaries} \label{one}
In this section we fix our notation, verify the existence of a component parameterizing a pair of linear spaces (Proposition \ref{component}) and describe some of its properties.

\medskip
\textbf{Notation:} Let $\mathbf{k}$ be an algebraically closed field. For Sections \ref{one} - \ref{three} we will assume $\mathrm{char}\, \kk \ne 2$ and for Sections \ref{four} - \ref{seven} we will assume $\mathrm{char} \, \kk = 0$. We use $S$ to denote the polynomial ring $\mathbf{k}[x_0,\dots,x_n]$ and $S_d$ to denote the subspace of monomials of degree $d$. For a homogenous ideal $I \subseteq S$ we use $I_d$ to denote the subspace of degree $d$ elements of $I$.  We use $[I]$ or $[X]$ to denote the $\kk$-point in the Hilbert scheme corresponding to $X=\Proj(S/I) \subseteq \mathbf{P}^n$ and we use $P_X(t)$ or $P_{S/I}(t)$ to denote its Hilbert polynomial. The ideal associated to a subscheme always refers to its saturated ideal. 

We use $\Gr(r,n)$ to denote the Grassmannian variety parameterizing $r$-dimensional linear spaces in $\PP^n$. The \textbf{span} of a subscheme $X \subseteq \PP^n$ is the linear subspace $V(H^0(\PP^n,I_X(1))) \subseteq \PP^n$.
The letters $c$ and $d$ are reserved for the codimension of linear spaces in $\mathbf{P}^n$; throughout the paper, we always assume $n \geq d \geq c \geq 2$. Similarly we reserve the letter $k = c = d$ for the case they are equal.

All the divisors we will consider are assumed to be Cartier. Given a smooth variety $Y$, we let $N^1(Y)$ denote the group of Cartier divisors modulo numerical equivalence. $\Nef(Y)$ and $\Eff(Y)$ denote the nef and effective cones of $Y$, respectively. We use $\langle D_1,\dots, D_l \rangle$ to denote the convex cone in $N^1(Y)\otimes \mathbf{R}$ generated by the divisors $D_i$. For more details we refer to \cite[Chapter 1]{debarre}.

\medskip
Let $X$ denote the union of an $(n-c)$-plane and $(n-d)$-plane meeting transversely in $\mathbf{P}^n$. It is clear that $X$ is parameterized by an open subset of $\mathbf{Gr}(n-c,n) \times \mathbf{Gr}(n-d,n)$ of dimension $c(n-c+1)+d(n-d+1)$. If we show that the tangent space to $[X]$ on its Hilbert scheme has dimension $c(n-c+1)+d(n-d+1)$, it will follow immediately that there is an  irreducible component of $\Hilb^{P_{n-c,n-d}^n(t)} \, \mathbf{P}^n$ whose general member parameterizes $X$ and whose natural scheme structure is reduced. 

Since $X$ is projectively equivalent to $Z = V(x_0,\dots,x_{c-1}) \cup V(x_{n-d+1},\dots,x_n)$; thus it suffices to compute the tangent space to $[Z]$ on its Hilbert scheme. For the rest of this section we fix $Z$ and $P(t) = P_{n-c,n-d}^n(t)$. 

If $Z \simeq \mathbf{P}^{n-c} \sqcup \mathbf{P}^{n-d}$ is a disjoint union of linear spaces, it is smooth; this occurs if and only if $n \leq c+d-1$. In this case we have a splitting of normals sheaves
$$
\mathscr{N}_{Z/\mathbf{P}^n} = \mathscr{N}_{\mathbf{P}^{n-c}/\mathbf{P}^n} \oplus \mathscr{N}_{\mathbf{P}^{n-d}/\mathbf{P}^n} \simeq \mathscr{O}_{\mathbf{P}^{n-c}}^c(1) \oplus \mathscr{O}_{\mathbf{P}^{n-d}}^d(1).
$$
Thus we obtain, $h^0(\mathbf{P}^n,\mathscr{N}_{Z/\mathbf{P}^n}) = c(n-c+1) + d(n-d+1)$ and $h^1(\mathbf{P}^n,\mathscr{N}_{Z/\mathbf{P}^n})= 0$. It follows that $[Z]$ is a smooth point on its Hilbert scheme \cite[Theorem 1.1c]{hdeform}. If $n>c+d-1$, we will explicitly compute the tangent space to $[Z]$ using the following result:

 \begin{thm}[Comparison Theorem \cite{ps}] \label{compare} Let $X\subseteq\mathbf{P}^{n}$ be a subscheme with ideal $I_X=(f_{1},\dots,f_{r}) \subseteq S$ where $\deg f_{i}=e_{i}$ satisfying, $(S/I_X)_{e}\simeq H^{0}(\PP^n,\mathscr{O}_{X}(e))$ for $e=e_{1},\dots,e_{r}$. Then there is an isomorphism between the universal deformation space of $I_X$ and that of $X$.
In particular, $T_{[X]}\, \Hilb^{P(t)} \, \mathbf{P}^n = H^0(\mathbf{P}^n,\mathscr{N}_{X/\mathbf{P}^n}) = \Hom_S(I_X,S/I_X)_0$.
\end{thm}

\begin{remark} \label{applycomparison}
With notation as in the above Theorem, consider the following exact sequence in local cohomology \cite[Corollary A1.12]{syzygies},
\begin{equation*}
0 \longrightarrow H_{\mathfrak{m}}^{0}(S/I_X) \longrightarrow S/I_X \longrightarrow 
				H^0_{\star}(\mathbf{P}^n,\mathscr{O}_{X}) \longrightarrow H_{\mathfrak{m}}^{1}(S/I_X) \longrightarrow 0.
\end{equation*}
If we show that $H_{\mathfrak{m}}^{i}(S/I_X)_{e} = 0$ for $e = e_1,\dots,e_r$ and $i=0,1$, then the Comparison theorem would apply. Here are two instances in which this is true
\begin{enumerate}
\item The depth of $S/I_X$ is at least $2$ \cite[Corollary A1.13]{syzygies}.
\item The Castlenuovo-Mumford regularity of the ideal $I_X$ is $\min{\{e_1,\dots,e_r\}}$ \cite[Proposition 4.16]{syzygies}. Note that $\mathrm{reg}(I_X) = \mathrm{reg}(S/I_X)+1$.
\end{enumerate}
\end{remark}

Since $n> c+d-1$, the depth of $S/I_Z$ is at least $2$. It follows from the previous Remark that the comparison theorem applies for $Z$.

\begin{lemma} \label{tan1} We have $\dim_{\kk} T_{[Z]} \Hilb^{P(t)} \mathbf{P}^n = c(n-c+1)+d(n-d+1) $.  
\end{lemma}
\begin{proof} 
We only need to consider the case $n > c+d-1$. Moreover, it suffices to show that the tangent space dimension is at most $c(n-c+1)+d(n-d+1)$. In particular it is enough to show that any $\varphi \in \Hom(I_Z,S/I_Z)_0$ can be written as
\begin{equation} \label{phi}
\begin{aligned}
\varphi(x_ix_j) = \sum_{\ell =0}^{n-d} a^{j}_{\ell}x_ix_{\ell} + \sum_{\ell = c}^{n} b^{i}_{\ell}x_jx_{\ell}
\end{aligned}
\end{equation}
for any $0 \leq i \leq c-1$ and $n-d+1 \leq j \leq n$ with some constants, $a_{\ell}^{i},b_{\ell}^i \in \kk$. 

Let us first show that $\varphi(x_ix_j)$ is supported on $\{x_ix_{0},\dots,x_ix_{n-d},x_jx_{c},\dots,x_jx_n\}$. Let $i,j$ be any integers satisfying $0 \leq i \leq c-1$ and $n-d+1 \leq j \leq n$. Choose $j'$ such that $n-d+1 \leq j' \leq n$ and $j \ne j'$.  Since $\varphi$ is an $S$-module homomorphism we have, $x_{j'}\varphi(x_{i}x_{j}) = x_{j}\varphi(x_{i}x_{j'})$. 
This implies that $x_{j}$ divides every non-zero monomial in $\varphi(x_{i}x_{j})$ that is not annihilated by $x_{j'}$ in $S/I_Z$. It follows that $\varphi(x_ix_j)$ is supported on 
$$
 \mathcal{C} = \{x_px_q:0 \leq p \leq c-1, 0 \leq  q \leq n-d\} \cup \{x_jx_c,\dots,x_jx_{n}\}.
$$

Similarly, choose $i'$ such that $0 \leq i' \leq c-1$ and $i' \ne i$. Then the equality $x_{i'}\varphi(x_{i}x_{j}) = x_{i}\varphi(x_{i'}x_{j})$ implies $x_{i}$ divides every monomial in $\varphi(x_{i}x_{j})$ that is not annihilated by $x_{i'}$. Once again we see that $\varphi(x_ix_j)$ is supported on 
$$
\mathcal{C}' = \{x_ix_0,\dots,x_ix_{n-d}\} \cup \{x_px_q:c \leq p \leq n, n-d+1 \leq  q \leq n\}. 
$$
Thus $\varphi(x_ix_j)$ is supported on 
$
\mathcal{C} \cap \mathcal{C}' = \{x_ix_{0},\dots,x_ix_{n-d},x_jx_{c},\dots,x_jx_n\}.
$ 

For any $i,j$, write $\varphi(x_ix_j) = \sum_{\ell=0}^{n-d}a^{i,j}_{\ell}x_ix_{\ell} + \sum_{\ell=c}^{n}b^{i,j}_{\ell}x_jx_{\ell}$ with $b^{ij}_{\ell},a^{ij}_{\ell} \in \kk$. Using the relation $x_{j'}\varphi(x_{i}x_{j}) = x_{j}\varphi(x_{i}x_{j'})$ we see that $b_{\ell}^{i,j} = b_{\ell}^{i,j'}$ for each $\ell$ and all $j,j'$. Using the relation $x_{i'}\varphi(x_{i}x_{j}) = x_{i}\varphi(x_{i}'x_{j})$ we obtain $a_{\ell}^{i,j} = a_{\ell}^{i',j}$ for each $\ell$ and all $i,i'$. Thus $\varphi$ is of the form described in (\ref{phi}).
\end{proof}

We immediately deduce the following.

\begin{prop} \label{component} There is an integral component of $\Hilb^{P(t)} \, \mathbf{P}^n$, denoted $\mathcal{H}_{n-c,n-d}^n$ or $\mathcal{H}_{n-c,n-d}(\mathbf{P}^n)$,  whose general point parameterizes an $(n-c)$-plane and an $(n-d)$-plane meeting transversely in $\mathbf{P}^n$.
\end{prop}

In the introduction we defined a rational map (\ref{XI})
\begin{equation*}
\Xi: \mathbf{Gr}(n-c,n) \times  \mathbf{Gr}(n-d,n) \dashrightarrow \mathcal{H}_{n-c,n-d}^n, \quad (\Lambda,\Lambda') \mapsto [I_{\Lambda}I_{\Lambda'}].
\end{equation*}
This map is well defined along the locus where $\Lambda, \Lambda'$ meet transversely, because in this situation $I_{\Lambda}I_{\Lambda'} = I_{\Lambda}\cap I_{\Lambda'}$. In many cases, $\Xi$ is in fact defined on a slightly larger open set.

\begin{lemma} \label{extend1} Let $n \geq c+d-1$. The rational map $\Xi$ extends to the complement of $\Gamma_{c-1}$.
\end{lemma}
\begin{proof}  We need to show that $\Xi$ is defined along $\Gamma_c \setminus\, \Gamma_{c-1}$. Up to projective equivalence, an element of $\Gamma_c \setminus \Gamma_{c-1}$ is  of the form $V(x_0,\dots,x_{c-1})\cup V(x_0,x_{c},\dots,x_{c+d-2})$. It suffices to show that $J=(x_0,\dots,x_{c-1})(x_0,x_{c},\dots,x_{c+d-2})$ has Hilbert polynomial $P(t)$. It follows by inspecting the minimal generators of $J$ that for any $t \geq 1$, $(S/J)_{t}$ is spanned by
$$
x_0\kk[x_{c+d-1},\dots,x_n]_{t-1} \oplus \bigoplus_{i=1}^{c-1}x_i\kk[x_i,\dots,x_{c-1},x_{c+d-1},\dots,x_n]_{t-1} \oplus \kk[x_{c},\dots,x_n]_t.
$$
Thus the Hilbert polynomial of $S/J$ is
$$
\begin{aligned}
\binom{n-c-d+t}{t-1} + \sum_{i=1}^{c-1} \binom{n-d-i+t}{t-1} + \binom{n-c+t}{t}.
\end{aligned}
$$ 
Using the "Hockey-Stick" identity this simplifies to
$$
\binom{n-c+t}{t} + \binom{n-d+t}{t} - \binom{n-c-d+t}{t} = P(t).
$$
\end{proof}

\begin{lemma} \label{twotoone} Let $n \geq c+d-1$ and consider the open set
$$
\mathcal{V} =  (\mathbf{Gr}(n-c,n)\times \mathbf{Gr}(n-d,n)) \setminus \, \Gamma_{c-1} \subseteq \mathbf{Gr}(n-c,n)\times \mathbf{Gr}(n-d,n).
$$ 
The morphism $\Xi|_{\mathcal{V}}: \mathcal{V} \longrightarrow \mathcal{H}_{n-c,n-d}^n$ is injective if $c \ne d$ and two-to-one if $c=d$. 
\end{lemma}
\begin{proof} 
Assume $\Xi|_{\mathcal{V}}(\Lambda,\Lambda') =  \Xi|_{\mathcal{V}}(\tilde{\Lambda},\tilde{\Lambda}') = [Y]$ for some scheme $Y$. Observe that $I_{\Lambda}I_{\Lambda'}$ is a saturated ideal. Indeed, up to projective equivalence, $\Lambda \cup \Lambda' =V(x_0,\dots,x_{c-1})\cup V(x_c,\dots,x_{c-d-2},x_i)$ with $i \in \{0,c-d-1\}$. In both cases, $I_{\Lambda}I_{\Lambda'}$ is clearly saturated. Thus we have $I_Y = I_{\Lambda}I_{\Lambda'}$ and taking  nilradicals we obtain
$$
I_{\Lambda \cup \Lambda'} =I_{\Lambda} \cap I_{\Lambda'} = \sqrt{I_{\Lambda} \cap I_{\Lambda'}} = \sqrt{I_{\Lambda}I_{\Lambda'}} = I_{Y_{\text{red}}}.
$$ 
Similarly, $I_{\tilde{\Lambda} \cup \tilde{\Lambda}'} = I_{Y_{\text{red}}}$. Equating the two expressions we have $\Lambda \cup \Lambda' = \tilde{\Lambda} \cup \tilde{\Lambda}'$. The conclusion now follows.
\end{proof}


\section{Structure of $\mathcal{H}_{n-k,n-k}^n$} \label{two}
This section is devoted to an analysis of $\HH_{n-k,n-k}^n$. The first major goal of this section is to prove that $\mathcal{H}_{n-k,n-k}^n$ is smooth. We start with the case when the pair of planes parameterized spans $\mathbf{P}^n$. We construct a bijective morphism from a non-singular variety to $\mathcal{H}_{n-k,n-k}^n$ and deduce this is an isomorphism by proving its differential is injective (Theorem \ref{main}). For the case where the pair of planes do not span $\mathbf{P}^n$, we construct a certain fibration to reduce to the case where they do span (Corollary \ref{mainthree}).

\smallskip
Let $n \geq 2k-1$ and $\mathcal{X}_0 = \Gr(n-k,n)^2$. For each $1 \leq v \leq k-1$, let $\mathcal{X}_v = \Bl_{\Gamma_v} \cdots \Bl_{\Gamma_1} \mathcal{X}_0$ and let $\pi_v:\mathcal{X}_v \longrightarrow \mathcal{X}_0$ be the blow-up morphism. 
The map (\ref{XI}) induces a  rational map
\begin{equation} 
\label{ximap} \Xi : \mathcal{X}_{k-1} = \Bl_{\Gamma_{k-1}}\cdots \Bl_{\Gamma_{1}}\mathbf{Gr}(n-k,n)^{2} \dashrightarrow \mathcal{H}_{n-k,n-k}^n
\end{equation}
defined away from the strict transforms of the exceptional divisors. In order to study the structure of $\HH_{n-k,n-k}^n$, we will begin by extending $\Xi$ to a morphism on $\XX_{k-1}$. 

\smallskip
For each ordered basis $\EE = \{e_0,\dots,e_{n}\}$ of $S_1$ we obtain an affine neighbourhood $U_{\EE} = \Spec \kk[a_{i,j},b_{i,j}]_{0 \leq i \leq k-1}^{k \leq j \leq n}$ of $\XX_0$ such that the $\kk$-points of $U_{\EE}$ correspond to
{\small
\begin{equation} \label{bfnotation}
(\Lambda(\mathbf{a}),\Lambda(\mathbf{b})) := (V(e_{0}+\sum_{j=k}^{n}a_{0,j}e_{j},\dots,e_{k-1}+\sum_{j=k}^{n}a_{k-1,j}e_{j}),V(e_{0}+\sum_{j=k}^{n}b_{0,j}e_{j},\dots,e_{k-1}+\sum_{j=k}^{n}b_{k-1,j}e_{j})).
\end{equation}}

It is clear that as $\EE$ ranges over all ordered basis of $S_1$, the set of $U_{\EE}$ cover $\XX_0$. In particular, it suffices to  extend $\Xi$ along each $\pi^{-1}_{k-1}(U_{\EE})$ in a compatible way. For notational convenience we may assume $\EE = \{x_0,\dots,x_n\}$ and let $U_0 = U_{\EE}$. Observe that the locus $\Gamma_v \cap U_0$ is cut out by the ideal generated by the $v\times v$ minors of the matrix
\begin{equation*}
M = 
\begin{pmatrix}
a_{0,k}-b_{0,k} & \cdots & a_{0,n}-b_{0,n}\\
\vdots &  & \vdots\\
a_{k-1,k}-b_{k-1,k} & \cdots & a_{k-1,n}-b_{k-1,n}
\end{pmatrix}.
\end{equation*}
Thus $\pi_{k-1}^{-1}(U_0)$ is obtained by blowing up $U_0$ along the strict transforms of the ideal generated by the $v\times v$ minors of $M$ for $v=1,\dots,k-1$, in that order.
\medskip

\begin{prop} \label{bigmatrix} For each $1 \leq v \leq k-1$, there exists non-singular affine open subsets $U_v \subseteq \XX_v$ such that the following hold.
\begin{enumerate}
\item We have $ U_v \subseteq \Bl_{\Gamma_v \cap U_{v-1}} U_{v-1} \subseteq \XX_v$.
\item   On the open set $U_v$, the matrix $\pi_v^{\star}(M)$ is row equivalent to the matrix
\begingroup\makeatletter\def\f@size{8}\check@mathfonts
\begin{equation*}
\begin{tikzpicture}[baseline=(current bounding box.center)]
\matrix (m) [matrix of math nodes, nodes in empty cells, left delimiter={(}, right delimiter={)}]{
\lambda_{1}\cdots\lambda_{v}(T_{0,k}^{(v)}-T_{0,n-v+1}^{(v)}T_{k-v,k}^{(v)}) & \cdots & \lambda_{1}\cdots\lambda_{v}(T_{0,n-v}^{(v)}-T_{0,n-v+1}^{(v)}T_{k-v,n-v}^{(v)}) & 0 & \cdots & 0 & 0\\
\vdots &  & \vdots & \vdots&   &\vdots & \vdots  \\
\lambda_{1}\cdots\lambda_{v}(T_{k-v-1,k}^{(v)}-T_{k-v-1,n-v+1}^{(v)}T_{k-v,k}^{(v)}) & \cdots & \lambda_{1}\cdots\lambda_{v}(T_{k-v-1,n-v}^{(v)}-T_{k-v-1,n-v+1}^{(v)}T_{k-v,n-v}^{(v)}) & 0& \cdots & 0 & 0\\
\lambda_{1}\cdots\lambda_{v}T_{k-v,k}^{(v)} & \cdots  & \lambda_{1}\cdots\lambda_{v}T_{k-v,n-v}^{(v)} & \lambda_{1}\cdots\lambda_{v} & \ddots & \vdots  & \vdots\\
\vdots &  & \vdots & \vdots & \ddots & 0 & \vdots \\
\lambda_{1}\lambda_2T_{k-2,k}^{(2)} & \cdots & \lambda_{1}\lambda_2T_{k-2,n-v}^{(2)}  & \lambda_{1}\lambda_2T_{k-2,n-v+1}^{(2)} & \cdots   & \lambda_1\lambda_2 & 0\\
\lambda_{1}T_{k-1,k}^{(1)} & \cdots & \lambda_{1}T_{k-1,n-v}^{(1)} & \lambda_{1}T_{k-1,n-v+1}^{(1)}   & \cdots & \lambda_{1}T_{k-1,n-1}^{(1)} & \lambda_{1} \\
} ;
\end{tikzpicture}
\end{equation*}
\endgroup
where $$\lambda_1 = a_{k-1,n}-b_{k-1,n} \text{  and  }  \lambda_{i}=T_{k-i,n-i+1}^{(i-1)}-T_{k-i,n-i+2}^{(i-1)}T_{k-i+1,n-i+1}^{(i-1)} \, \text{  for each  } 2 \leq i \leq k-1.$$
\item The strict transform of $\Gamma_{v+1}$ on $U_v$ is cut out by $$(T_{i,j}^{(v)}-T^{(v)}_{i,n-v+1}T_{k-v,j}^{(v)})_{k\leq j\leq n-v}^{0\leq i\leq k-v-1}.$$
\item $\Gamma_{v+1} \cap U_v$ is non-singular and the blowup along this locus is given by 
$$
\Bl_{\Gamma_{v+1} \cap U_v} U_v := \Proj \kk[U_v][T_{i,j}^{(v+1)}]_{i,j}/(\text{Koszul Relations}).
$$

\end{enumerate}
\end{prop}
\begin{proof} We begin with the definition of $U_1$. Since $\Gamma_1$ is cut out by $(a_{i,j} - b_{i,j})_{i,j}$ on $U_0$, it is a non-singular subscheme and we have $\Bl_{\Gamma_1 \cap U_0 } U_0 = \Proj \kk[U_0][T^{(1)}_{i,j}]_{i,j}/(\text{Koszul relations})$. We define $U_1 = D(T^{(1)}_{k-1,n})$.

Let $M_{v}$ denote the matrix appearing in item (ii). We will prove items (i) - (iv) inductively starting with $v=1$. Item (i) is true for $v=1$ by construction. On the open set $U_1$, the Koszul relations simplify to $a_{i,j}-b_{i,j} = \lambda_1T_{i,j}^{(1)}$; here we have set $T_{k-1,n}^{(1)} =1$. Substituting this into the matrix $\pi_{1}^{\star}(M)$ and subtracting appropriate multiples of the bottom row from every other row, we obtain the matrix
$$
M_1 = \begin{pmatrix}
\lambda_1(T_{0,k}^{(1)}-T_{0,n}^{(1)}T_{k-1,k}^{(1)}) & \cdots & \lambda_1(T_{0,n-1}^{(1)}-T_{0,n}^{(1)}T_{k-1,n-1}^{(1)}) & 0\\
\vdots &  & \vdots & \vdots\\
\lambda_1(T_{k-2,k}^{(1)}-T_{k-2,n}^{(1)}T_{k-1,k}^{(1)}) &  & \lambda_1(T_{k-2,n-1}^{(1)}-T_{k-2,n}^{(1)}T_{k-1,n-1}^{(1)}) & 0\\
\lambda_1T_{k-1,k}^{(1)} & \cdots & \lambda_1T_{k-1,n-1}^{(1)} & \lambda_1
\end{pmatrix}.
$$
This proves item (ii) for $v=1$.  The ideal generated by the $2\times2$ minors of $M_1$ is $\lambda_1^2(T_{i,j}^{(1)}-T_{i,n}^{(1)}T^{(1)}_{k-1,j})_{0 \leq j \leq n-1}^{0 \leq i \leq k-2}$. Thus the ideal of the strict transform of $\Gamma_{2}$ is $(T_{i,j}^{(1)}-T_{i,n}^{(1)}T^{(1)}_{k-1,j})_{0 \leq j \leq n-1}^{0 \leq i \leq k-2}$. Since this ideal is generated by a  regular sequence, the blowup along it is non-singular and equal to $\Bl_{\Gamma_{2}\cap U_{1}} U_{1} :=  \Proj\,\kk[U_{1}][T_{i,j}^{(2)}]_{i,j}/(\text{Koszul relations})$. This proves item (iii) and (iv) for $v=1$.

Now assume items (i) - (iv) have been proved for some $1 \leq  v \leq k-2$. Define $U_{v+1} = D(T_{k-v-1,n-v}^{(v+1)})$; equivalently let  $T_{k-v-1,n-v}^{(v+1)}= 1$. Then the Koszul relations on this open simplify to $T_{i,j}^{(v)}-T^{(v)}_{i,n-v+1}T_{k-v,j}^{(v)} = \lambda_{v+1}T_{i,j}^{(v+1)}$. Once we substitute this into the matrix $M_v$, it is straightforward to row reduce the matrix so that it becomes $M_{v+1}$. Items (i) - (iv) will follow immediately as explained in the previous paragraph.
\end{proof}


\begin{remark} \label{coordinates} It follows from Proposition \ref{bigmatrix} that a set of algebraically independent coordinates on $ U_{k-1}$ is 
\begin{equation*}
\{b_{i,j}\}_{0 \leq i \leq k-1}^{k \leq j \leq n} \cup\{T^{(j)}_{i,n-j+1}\}_{1 \leq j \leq k-1}^{0 \leq i \leq k-1-j} \cup \{\lambda_1,\dots,\lambda_{k-1}\} \cup \{T_{k-i,j}^{(i)}\}_{k \leq j \leq n-i}^{1 \leq i \leq k-1} \cup \{T_{0,j}^{(k)}\}_{k \leq j \leq n-k+1}
\end{equation*}
with $T^{(k)}_{0,j} = T_{0,j}^{(k-1)}-T^{(k-1)}_{0,n-k+2}T_{1,j}^{(k-1)}$ for all $j$. 
\end{remark}

\begin{prop} \label{mainone} Let $n \geq 2k-1$. The rational map $\Xi$ (\ref{ximap}) extends to a morphism $U_{k-1} \longrightarrow \HH_{n-k,n-k}^n$.
\end{prop}
\begin{proof} We will use $\mathbf{a}$ to denote the tuple $(a_{i,j})_{i,j}$ and similarly use $\mathbf{b}$ and $\mathbf{T}^{(v)}$ to denote their corresponding tuples. Moreover, we will use $\Lambda(\mathbf{a})$ to denote the $(n-k)$-plane corresponding to $\mathbf{a}$ as in (\ref{bfnotation}).
For each $0 \leq i \leq k-1$ let $y_i = x_i + \sum_{j=k}^nb_{i,j}x_j$. At the moment, $\Xi$ maps 
{\small
\begin{align}
\label{originalXi} (\mathbf{a}, \mathbf{b},\mathbf{T}^{(1)},\dots,\mathbf{T}^{(k)})
		&\mapsto  \left[ I_{\Lambda(\mathbf{a})}I_{\Lambda(\mathbf{b})} \right] \\ 
		& = \left[ (y_{0}+\sum_{j=k}^{n}(a_{0,j}-b_{0,j})x_{j},\dots,y_{k-1}+\sum_{j=k}^{n}(a_{k-1,j}-b_{k-1,j})x_{j}) 
		(y_0,\dots,y_{k-1}) \right] \notag 
\end{align}
}
and this is undefined along the strict transforms of the exceptional divisors. Although we may express $\mathbf{a}$ in terms of $\mathbf{b}$ and $\{\mathbf{T}^{(v)}\}_v$, we will still describe formulas in terms of $\mathbf{a}$ as it simplifies the exposition.

Observe that a minimal set of generators for $I_{\Lambda(\mathbf{a})}$ is given by the rows of $\begin{bmatrix} \text{Id}_{k\times k} \, \vert \,  M\end{bmatrix} \boldsymbol{z}^{T}$ where $\boldsymbol{z} = \begin{bmatrix} y_0 & \cdots & y_{k-1} & x_k & \cdots & x_n \end{bmatrix}$ is a row vector. Applying row operations to $\begin{bmatrix} \text{Id}_{k\times k} \, \vert \,  M\end{bmatrix}$ will produce different minimal sets of generators. In particular, applying the row operations we did to $M$ to get $M_{k-1}$ (Proposition \ref{bigmatrix} (ii)) to the matrix $\begin{bmatrix} \text{Id}_{k\times k} \, \vert \,  M\end{bmatrix}$
we obtain a new set of generators $\alpha_0,\dots,\alpha_{k-1}$ of $I_{\Lambda(\mathbf{a})}$ where
\begin{equation*}
\alpha_{p}= y_{p}-\sum_{j=1}^{k-1-p}T_{p,n-j+1}^{(j)}y_{k-j}+\sum_{j=k}^{n-(k-1-p)}\lambda_{1}\cdots\lambda_{k-p}T_{p,j}^{(k-p)}x_{j} \quad \text{for} \quad  0 < p \leq k-1
\end{equation*}
and
\begin{equation*}
\alpha_{0}= y_{0}-\sum_{j=1}^{k-1}T_{0,n-j+1}^{(j)}y_{k-j}+\sum_{j=k}^{n-(k-1)}\lambda_{1}\cdots\lambda_{k-1}T_{0,j}^{(k)}x_j
\end{equation*}
with $T^{(k)}_{0,j} = T_{0,j}^{(k-1)}-T^{(k-1)}_{0,n-k+2}T_{1,j}^{(k-1)}$ for all $j$. By construction, $T_{k-v,n-v+1}^{(v)}= 1$  for all $1 \leq v \leq k-1$. 

\smallskip
For $0\leq  p<q\leq k-1$ define the following "cross terms"
\begin{equation*}
\begin{aligned}
\beta_{p,q}= \left(y_{p}-\sum_{j=1}^{k_p}T_{p,n-j+1}^{(j)}y_{k-j}\right)\left(\sum_{j=k}^{n-k_q}T_{q,j}^{(k-q)}x_{j}\right)- \lambda_{p,q}\left(y_{q}- \sum_{j=1}^{k_q}T_{q,n-j+1}^{(j)}y_{k-j}\right)\left(\sum_{j=k}^{n-k_p}T_{p,j}^{(k-p)}x_{j}\right),
\end{aligned}
\end{equation*}
where $k_p = k-1-p$ for all $p$ and 
$\lambda_{p,q}  =  \begin{cases} \lambda_{k-q+1}\cdots\lambda_{k-p} \, \text{ if } p > 0
	\\ \lambda_{k-q+1}\cdots\lambda_{k-1} \,\, \text{ if } p =0.  \end{cases}$ 

\smallskip	
Note that our convention implies $\lambda_{0,1}=1$. Extend $\Xi$ to $U_{k-1}$ by mapping
{\small
\begin{equation} \label{finalmap1}
\begin{aligned}
(\mathbf{a}, \mathbf{b},\mathbf{T}^{(1)},\dots,\mathbf{T}^{(k)}) & \mapsto \left[ I_{\Lambda(\mathbf{a})}(y_{0},\dots,y_{k-1})+(\beta_{p,q})_{0\leq p<q\leq k-1} \right] \\
& = \left[ \left(x_i+\sum_{j =k}^na_{i,j}\right)_{0 \leq i \leq k-1}\left(x_i+\sum_{j =k}^nb_{i,j}\right)_{0 \leq i \leq k-1} + \left(\beta_{p,q}\right)_{0\leq p<q\leq k-1}\right].
\end{aligned}
\end{equation}
}

Note that (\ref{finalmap1}) extends the original rational map (\ref{originalXi}). Indeed, (\ref{originalXi}) is defined away from the strict transform of all the the exceptional divisors; this is the locus where $\lambda_{1},\dots,\lambda_{k-1}\ne0$. In this case we have
\begin{equation} 
\label{simplifygens}
(y_0,\dots,y_{k-1})I_{\Lambda(\mathbf{a})} \ni \left(y_{p}-\sum_{j=1}^{k_p}T_{p,n-j+1}^{(j)}y_{k-j}\right)\alpha_{q}-\left(y_{q}-\sum_{j=1}^{k_q}T_{q,n-j+1}^{(j)}y_{k-j}\right)\alpha_{p}  = \lambda_1 \cdots \lambda_{k-q} \beta_{p,q}.
\end{equation}
Thus $\beta_{p,q} \in I_{\Lambda(\mathbf{a})}(y_0,\dots,y_{k-1})$ and (\ref{originalXi}) and (\ref{finalmap1}) coincide.

To show that the image of (\ref{finalmap1}) is well defined, it is enough to show that the Hilbert polynomial of an ideal $J = I_{\Lambda(\aaa)}I_{\Lambda(\bb)}+(\beta_{p,q})_{0\leq p<q\leq k-1}$ in this image is $P_{n-k,n-k}^n(t)$. In Lemma \ref{GROBNER} we define a term order $>$ on $S$ for which 
$$
\inn_{>} J = (x_0,\dots,x_{k-1})^2+(x_px_{n-k_q})_{0 \leq p < q \leq k-1}.
$$
Since there is a flat degeneration from $J$ to $\inn_{>} J$ it suffices to show $\inn_{>} J$ has the desired Hilbert polynomial. It is easy to see that $(S/\inn_{>}J)_{t}$ is spanned by
$$
\bigoplus_{i=0}^{k-1}x_i\kk[x_k,\dots,x_{n-k+i+1}]_{t-1}\oplus \kk[x_{k},\dots,x_n]_t.
$$
Using this and the Hockey-Stick identity we deduce that Hilbert polynomial of $S/\inn_{>} J$ is
$$
\binom{n-k+t}{t} + \sum_{i=0}^{k-1}\binom{n-2k+i+t}{t-1} = \binom{n-k+t}{t} + \binom{n-k+t}{t} - \binom{n-2k+t}{t} =P_{n-k,n-k}^n(t).
$$
\end{proof}

Prior to proving Lemma \ref{GROBNER} we need the following auxiliary result.

\begin{lemma} \label{projequiv} The ideal $I_{\Lambda(\mathbf{a})}I_{\Lambda(\mathbf{b})}+(\beta_{p,q})_{0\leq p<q\leq k-1}$ in the image of Equation (\ref{finalmap1}) is projectively equivalent to an ideal of the form
\begin{equation} \label{projideal}
(x_{p}+\mu_{p,k}x_{n-k_p})_{0\leq p\leq k-1}(x_0,\dots,x_{k-1}) + (x_px_{n-k_q} -\mu_{p,q} x_q x_{n-k_p})_{0\leq p <q \leq k-1},
\end{equation}
with $\mu_{i} \in \kk$ and $\mu_{p,q} = \mu_{k-q+1}\cdots \mu_{k-p}$ for any $0 \leq p < q \leq k$.
\end{lemma}
\begin{proof} Applying the projective transformation that maps $ x_i \mapsto x_{i}-\sum_{j\geq k}b_{i,j}x_{j}$ if $ i \leq k-1$ and fixes the other $x_i$, we may assume $\mathbf{b} = \mathbf{0}$. For each $0 \leq i \leq k-1$ let $\tau_i$ denote the map that sends $x_i \mapsto x_{i}+\sum_{j=1}^{k-i-1}T_{i,n-j+1}^{(j)}x_{k-j}$ and fixes the other $i$. It is clear that $\tau_{k-1} \circ \cdots \circ \tau_{0}(I)$ equals,
{\small
\begin{equation*}
\begin{aligned}
\left(x_p +\sum_{j=k}^{n-k_p}\lambda_1\cdots \lambda_{k-p}T_{p,j}^{(k-p)}x_{j}\right)_{0 \leq p \leq k-1}(x_0,\dots,x_{k-1}) + \left(x_p\left(\sum_{j=k}^{n-k_q}T_{q,j}^{(k-q)}x_{j}\right) -\lambda_{p,q} x_q\left(\sum_{j=k}^{n-k_p}T_{p,j}^{(k-p)}x_{j}\right)\right)_{p<q}
\end{aligned}
\end{equation*}
}

For each $0 \leq i \leq k-1$ let $\mu_i = \lambda_i$. If $T^{(k)}_{0,j} = 0$ for all $j$ then let $\mu_k = 0$. If not, choose the largest index $\ell$ for which $T^{(k)}_{0,\ell} \ne 0$ and let $\mu_k =  T^{(k)}_{0,\ell}$.

For each $1 \leq i \leq k-1$ consider the map $\tau_{n-k_i}$, that maps $x_{n-k_i} \mapsto x_{n-k_i} - \sum_{j=k}^{n-k_i-1}T_{i,j}^{(k-i)}x_{j}$ and fixes the other $x_i$. As we range over all $i$, we obtain maps $\tau_{n},\dots, \tau_{n-(k-2)}$. If $\mu_k=0$ let $\tau_{n-(k-1)} $ be the identity; else let $\tau_{n-(k-1)}$ denote the map that sends $x_{\ell} \mapsto x_{n-k_0}- \frac{1}{\mu_k}\sum_{j=k}^{\ell-1}T^{(k)}_{0,j}$, $x_{n-k_0} \mapsto x_{\ell}$ if $\ell < n-k_0$, and fixes the other $x_i$.

Using the fact that $T^{(k-i)}_{i,n-k_i} =1$ on the open set $U_{k-1}$, it is straightforward to check that $\tau_{n-(k-1)} \circ \cdots \tau_{n} \circ \tau_{k-1} \circ \cdots \circ \tau_{0}(I)$ is of the desired form.
\end{proof}

\begin{lemma} \label{GROBNER}   \label{realin}  Let $>$ denote the lexicographic ordering on $S$ with terms ordered by $x_{0} >x_{1}> \cdots >x_{k-1} > x_n >x_{n-1}> \cdots >x_k$. Let $J = I_{\Lambda(\mathbf{a})}I_{\Lambda(\bb)}+(\beta_{p,q})_{0 \leq p < q \leq k-1}$ denote the ideal in the image of Equation (\ref{finalmap1}). Then we have
\begin{equation*} 
\inn_{>} J = (x_0,\dots,x_{k-1})^2 + (x_px_{n-k_q})_{0 \leq p < q\leq k-1}
\end{equation*}
\end{lemma}
\begin{proof} Let $J'$ denote the ideal in (\ref{projideal}). We will first show that 
\begin{equation}
\label{projgrobner}
\inn_{>}J' = (x_0,\dots,x_{k-1})^2 + (x_px_{n-k_q})_{0 \leq p <q \leq k-1}.
\end{equation}

Let $\gamma_{p,q} = (x_p+\mu_{p,k}x_{n-k_p})x_q$ for $0\leq p \leq q\leq k-1$ and $\delta_{p,q} = x_px_{n-k_q} -\mu_{p,q} x_q x_{n-k_p}$ for $0\leq p < q \leq k-1$. Since $\inn_{>} \gamma_{p,q} = x_px_q$ and $\inn_{>} \delta_{p,q}= x_px_{n-k_q}$, to prove (\ref{projgrobner}) it is enough to show that $G = \{\gamma_{p,q},\delta_{p,q}\}_{p,q}$ is a Gr\"{o}bner basis for $J'$. Note that $G$ generates $J'$ because for $p< q$ we have
\begin{align}
\label{swapeq} (x_q+ \mu_{q,k}x_{n-k_q})x_p  &= (x_p+\mu_{p,k}x_{n-k_p})x_q  +  \mu_{q,k}(x_px_{n-k_q} - \mu_{p,q}x_qx_{n-k_p})\\
\nonumber & = \gamma_{p,q} +\mu_{q,k}\delta_{p,q} \in (G).
\end{align} 
Notice that $\mu_{p,q}\mu_{q,k} = \mu_{p,k}$ and this will be used repeatedly in the rest of the proof.

Given $a,b \in S$ we denote their \textit{S-pair} by $R(a,b) = (\frac{\inn_{>}b}{h})a - (\frac{\inn_{>}a}{h})b$ with $h = \mathrm{gcd}(\inn_{>}(a),\inn_{>}(b))$. To show that $G$ forms a Gr\"{o}bner basis we need to show that there is a \textit{standard expression} for the S-pairs in terms of elements of $G$ with no remainder \cite[Section 2.2-2.3]{monomialideals}.

\smallskip
\textbf{Case 1.} The standard expression of $R(\gamma_{p_1,q_1},\gamma_{p_2,q_2})$: Let $h = \mathrm{gcd}(\inn_>\gamma_{p_1,q_1},\inn_>\gamma_{p_2,q_2})$ and we may assume $p_1 \leq p_2$. If $h=1$ then $p_1 < p_2$ and we have
\begin{align*}
R(\gamma_{p_1,q_1},\gamma_{p_2,q_2})  &= x_{p_2}x_{q_2}\gamma_{p_1,q_1} - x_{p_1}x_{q_1}\gamma_{p_2,q_2}  \\
								&=  \mu_{p_1,k}x_{p_2}x_{q_2}x_{n-k_{p_1}}x_{q_1} 
								-  \mu_{p_2,k}x_{p_1}x_{q_1}x_{n-k_{p_2}}x_{q_2} \\
								 &= -\mu_{p_2,k}x_{q_1}x_{q_2}\delta_{p_1,p_2}.				
\end{align*}
This is obviously a standard expression with no remainder. If $h=x_{p_1}$ then $p_1 = p_2$ or $p_1 = q_2$; in the latter case we still have $p_1 =p_2$ as our assumptions imply $p_1 \leq p_2 \leq q_2$. Thus in both the situations we obtain $R(\gamma_{p_1,q_1},\gamma_{p_2,q_2}) = x_{q_2}\gamma_{p_1,q_1} - x_{q_1}\gamma_{p_1,q_2} = 0$. If $h = x_{q_1}$ we have either $q_1 = q_2$ or $q_1 = p_2$. If $q_1 = q_2$ then as shown above we obtain
\begin{align*}
R(\gamma_{p_1,q_1},\gamma_{p_2,q_2}) = x_{p_2}\gamma_{p_1,q_1} - x_{p_1}\gamma_{p_2,q_1}  
								= \mu_{p_1,k}x_{p_2}x_{n-k_{p_1}}x_{q_1} 
								- \mu_{p_2,k}x_{p_1}x_{n-k_{p_2}}x_{q_1}  
								 = -\mu_{p_2,k}x_{q_1}\delta_{p_1,p_2} .
\end{align*}
Similarly, if $q_1 = p_2$ we obtain
$R(\gamma_{p_1,q_1},\gamma_{p_2,q_2}) = x_{q_2}\gamma_{p_1,p_2} - x_{p_1}\gamma_{p_2,q_2} 
								 = -\mu_{p_2,k}x_{q_2}\delta_{p_1,p_2}
$
(if $p_1 = p_2$ this is just $0$). If $h = x_{p_1}x_{q_1}$ then we have $p_1 = q_1 = p_2 = q_2$ or $p_1 =p_2 < q_1 = q_2$; in either case $R(\gamma_{p_1,q_1},\gamma_{p_2,q_2}) = 0$. 

\smallskip
\textbf{Case 2.} The standard expression of $R(\delta_{p_1,q_1},\delta_{p_2,q_2})$: Let $h = \mathrm{gcd}(\inn_>\delta_{p_1,q_1},\inn_>\delta_{p_2,q_2})$ and assume $p_1 \leq p_2$. If $h = 1$ we have $p_1 < p_2$ and $q_1 \ne q_2$. Then we obtain
\begin{align*}
R(\delta_{p_1,q_1},\delta_{p_2,q_2})  &= x_{p_2}x_{n-k_{q_2}}\delta_{p_1,q_1} - x_{p_1}x_{n-k_{q_1}}\delta_{p_2,q_2}  \\
								&=  -\mu_{p_1,q_1}x_{p_2}x_{n-k_{q_2}}x_{q_1}x_{n-k_{p_1}} 
								+  \mu_{p_2,q_2}x_{p_1}x_{n-k_{q_1}} x_{q_2}x_{n-k_{p_2}}\\
								 &=  \mu_{p_2,q_2}x_{q_2}x_{n-k_{q_1}}\delta_{p_1,p_2}
								 	-x_{p_2}x_{n-k_{p_1}}(\mu_{p_1,q_1}x_{q_1}x_{n-k_{q_2}} 
								 	- \mu_{p_1,p_2}\mu_{p_2,q_2}x_{q_2}x_{n-k_{q_1}})		\\
								& = \begin{cases}
									 \mu_{p_2,q_2}x_{q_2}x_{n-k_{q_1}}\delta_{p_1,p_2} 
											- \mu_{p_1,q_1}x_{p_2}x_{n-k_{p_1}}\delta_{q_1,q_2} 
													\quad\quad \text{ if }  q_1 < q_2 \\
									 \mu_{p_2,q_2}x_{q_2}x_{n-k_{q_1}}\delta_{p_1,p_2} 
											+ \mu_{p_1,q_2}x_{p_2}x_{n-k_{p_1}}\delta_{q_2,q_1} 
												\quad\quad\text{ if }  q_2 < q_1.\\
								 	\end{cases}		
\end{align*}
Each of the above cases is a standard expression in terms of $G$ with no remainder \footnote{If $\mu_{p_2,q_2} \ne 0$ then $\inn_{>}R(\delta_{p_1,q_1},\delta_{p_2,q_2}) =  \mu_{p_2,q_2}x_{p_1}x_{n-k_{q_1}} x_{q_2}x_{n-k_{p_2}} $. This is greater or equal to $\inn_{>}(x_{q_2}x_{n-k_{q_1}}\delta_{p_1,p_2})$ and $\inn_{>}(x_{p_2}x_{n-k_{p_1}}\delta_{q_1,q_2})$.}. If $h = x_{n-{k_{q_1}}}$ we have $q_1 = q_2$ and $p_1 < p_2$. Then we obtain
\begin{equation*}
\begin{aligned}
R(\delta_{p_1,q_1},\delta_{p_2,q_2}) &= x_{p_2}\delta_{p_1,q_2} -x_{p_1}\delta_{p_2,q_2}  \\
			& = -\mu_{p_1,q_2}x_{p_2}x_{q_2}x_{n-k_{p_1}} + \mu_{p_2,q_2}x_{p_1}x_{q_2}x_{n-k_{p_2}} \\
			& = \mu_{p_2,q_2}x_{q_2}\delta_{p_1,p_2}.
\end{aligned}
\end{equation*}
If $h =x_{p_1}$ we have $p_1 = p_2$ and wlog we may assume $q_1 < q_2$. Then we have 
\begin{equation*}
\begin{aligned}
R(\delta_{p_1,q_1},\delta_{p_2,q_2}) &= x_{n-k_{q_2}}\delta_{p_1,q_1} -x_{n-k_{q_1}}\delta_{p_1,q_2}  \\
			& = -\mu_{p_1,q_1}x_{n-q_2}x_{q_1}x_{n-k_{p_1}} + \mu_{p_1,q_2}x_{n-k_{q_1}}x_{q_2}x_{n-k_{p_1}} \\
			& = -\mu_{p_1,q_1}x_{n-k_{p_1}}\delta_{q_1,q_2}.
\end{aligned}
\end{equation*}
Finally if $h=x_{p_1}x_{n-k_{q_1}}$ we have $p_1 = p_2 < q_1 = q_2$ and thus $R(\delta_{p_1,q_1},\delta_{p_2,q_2})  = 0$. 

\smallskip
\textbf{Case 3.} The standard expression of $R(\gamma_{p_1,q_1},\delta_{p_2,q_2})$: Let $h = \mathrm{gcd}(\inn_>\gamma_{p_1,q_1},\inn_>\delta_{p_2,q_2})$ and note that $h \in \{1,x_{p_1},x_{q_1}\}$. If $h= x_{p_1}$ we have $p_1 = p_2 $ and using (\ref{swapeq}) we obtain
\begin{align*}
R(\gamma_{p_1,q_1},\delta_{p_2,q_2}) &= x_{n-k_{q_2}}\gamma_{p_1,q_1} -x_{q_1}\delta_{p_1,q_2}  \\
			& = \mu_{p_1,k}x_{n-k_{q_2}}x_{n-k_{p_1}}x_{q_1} + \mu_{p_1,q_2}x_{q_1}x_{n-k_{p_1}}x_{q_2} \\
			& = \begin{cases}
					\mu_{p_1,q_2}x_{n-k_{p_1}}\gamma_{q_2,q_1} 
						\quad\quad\quad\quad\quad\quad\quad\quad\,\,\text{ if } q_1 \geq q_2 \\
					 \mu_{p_1,q_2}x_{n-k_{p_1}}\gamma_{q_1,q_2} 
					 		+\mu_{p_1,k}x_{n-k_{p_1}}\delta_{q_1,q_2}\, \text{ if } q_1 < q_2.
				\end{cases}
\end{align*}
Both these cases are standard expressions with no remainder. If $h=x_{q_1}$ then $q_1 = p_2$ and we obtain,
\begin{align*}
R(\gamma_{p_1,q_1},\delta_{p_2,q_2}) &= x_{n-k_{q_2}}\gamma_{p_1,p_2} -x_{p_1}\delta_{p_2,q_2}  \\
			& = \mu_{p_1,k}x_{n-k_{q_2}}x_{n-k_{p_1}}x_{p_2} + \mu_{p_2,q_2}x_{p_1}x_{n-k_{p_2}}x_{q_2} \\
			& = x_{n-k_{q_2}}\gamma_{p_1,q_2} - x_{p_1}\delta_{p_2,q_2}.
\end{align*}
Finally consider the case $h=1$. If we further assume $p_2 < p_1$ and $q_2 < p_1$ we have
\begin{align*}
R(\gamma_{p_1,q_1},\delta_{p_2,q_2})  &= x_{p_2}x_{n-k_{q_2}}\gamma_{p_1,q_1} - x_{p_1}x_{q_1}\delta_{p_2,q_2}  \\
								&=  \mu_{p_1,k}x_{p_2}x_{n-k_{q_2}}x_{n-k_{p_1}}x_{q_1}
								+ \mu_{p_2,q_2}x_{p_1}x_{q_1}x_{q_2}x_{n-k_{p_2}} \\
								& = \mu_{p_1,k}x_{n-k_{q_2}}x_{q_1}\delta_{p_2,p_1} 
								+\mu_{p_2,k}x_{n-k_{q_2}}x_{q_1}x_{p_1}x_{n-k_{p_2}} 
								+ \mu_{p_2,q_2}x_{p_1}x_{q_1}x_{q_2}x_{n-k_{p_2}} \\
								& = \mu_{p_1,k}x_{n-k_{q_2}}x_{p_2}\delta_{p_2,p_1}  
								+ \mu_{p_2,q_2}x_{n-k_{p_2}}x_{q_1}\gamma_{q_2,p_1}
\end{align*}
This is a standard expression with no remainder. We omit the other cases as their proofs are very similar (use Equation \ref{swapeq}). We have now shown that $G$ is a Gr\"obner basis for $J'$.

Since $J'$ and $\inn_{>} J'$ have the same Hilbert function (as graded $S$-modules) and $J$ is projectively equivalent to $J'$, $J$ and $\inn_{>}J'$ have the same Hilbert function. On the other hand, $(x_0,\dots,x_{k-1})^2 \subseteq \inn_{>}J$ and $x_px_{n-k_q} = \inn_{>}(\beta_{p,q}) \in \inn_{>}J$. Thus $\inn_{>}J \supseteq \inn_{>}J'$. Since these ideals have the same Hilbert function they must be equal, completing the proof. 
\end{proof}

\begin{remark} \label{knotation} For the rest of the paper, $>$ will always denote the term order from Lemma \ref{GROBNER} and $k_p$ will always denote $k-1-p$.
\end{remark}

The following Lemma sheds some light on the structure of the subschemes in the image of the morphism, $U_{k-1} \longrightarrow \HH_{n-k,n-k}^n$.

\begin{lemma} \label{anothermingens}  \label{saturated} 
Let $J = I_{\Lambda(\mathbf{a})}I_{\Lambda(\bb)}+(\beta_{p,q})_{0 \leq p < q \leq k-1}$ denote the ideal in the image of the morphism (\ref{finalmap1}). Then the following statements are true
\begin{enumerate}
\item The ideal $J$ is saturated. 
\item If all the $\lambda_i$ are non-zero and $\mathbf{T}^{(k)} \ne \mathbf{0}$ then $J$ is the ideal of a pair of $(n-k)$-planes meeting transversely.
\item If all the $\lambda_i$ are non-zero and $\mathbf{T}^{(k)} = \mathbf{0}$ then $\sqrt{J}$ is the ideal of a pair of $(n-k)$-planes meeting along an $(n-2k+1)$-plane.
\item Let $\ell$ be the smallest index for which $\lambda_{\ell} =0$. Then we have
$$
J = I_{\Lambda(\mathbf{a})}I_{\Lambda(\mathbf{b})}+(\beta_{p,q})_{0 \leq p < q \leq k-\ell}
$$
and $\sqrt{J}$ is the ideal of a pair of $(n-k)$-planes meeting along an $(n-k+1-\ell)$-plane.
\end{enumerate}
\end{lemma}

\begin{proof} Item (i) follows from the fact that  $\mathrm{depth}_{\mathfrak{m}}(S/J) \geq \mathrm{depth}_{\mathfrak{m}}(S/\inn_{>}J)  \geq 1$ where $\mathfrak{m} = (x_0,\dots,x_n)$. The first inequality is \cite[Theorem 3.3.4]{monomialideals} and the second inequality is true because $x_k$ is a non-zero divisor on $S/\inn_{>} J$.

Notice that $\Lambda(\aaa)$ and $\Lambda(\bb)$ meet along a $(n-k+1-\ell)$-plane precisely when the matrix $M$ (Proposition \ref{bigmatrix} (ii)) has rank $\ell-1$. As a consequence items (ii), (iii) and the second half of (iv) follow immediately. The other half of item (iv) follows from  Equation \ref{simplifygens} as it shows $\beta_{p,q} \in I_{\Lambda(\mathbf{a})}I_{\Lambda(\bb)}$ for any $q > k-\ell$.
\end{proof}

\begin{prop} \label{mainonee} Let $n \geq 2k-1$. Then $\Xi$ induces a surjective, $\GL(n+1)$-equivariant morphism
\begin{equation*} 
\widebar{\Xi}: \XX_{k-1}/\mathfrak{S}_2 \simeq \Bl_{\Gamma_{k-1}} \cdots \Bl_{\Gamma_1} \Sym^2 \mathbf{Gr}(n-k,n) \longrightarrow \mathcal{H}_{n-k,n-k}^n.
\end{equation*}
Moreover, the quotient $\XX_{k-1}/\mathfrak{S}_2 $ is non-singular.
\end{prop}
\begin{proof} In Proposition \ref{mainone} we showed that $\Xi$ extends to a map from $U_{k-1}$. We will now explain how the same argument gives a morphism on all of $\pi_{k-1}^{-1}(U_0)$. Consider a pair 
$$
\bg=(\bg^1,\bg^2) = ((\bg^1_1,\dots,\bg^1_k),(\bg^2_1,\dots,\bg^2_{k-1}))
$$ 
with $\bg^1$ an ordered $k$-subset of $\{0,\dots,k-1\}$ and $\bg^2$ an ordered $(k-1)$-subset of $\{k,\dots,n\}$. For any such $\bg$ we can define a sequence of open sets $U_{1}^{\bg},\dots, U_{k-1}^{\bg}$ such that 
\begin{enumerate}
\item[(1)] $U_1^{\bg} = D(T^{(1)}_{\bg^1_1,\bg^2_1}) \subseteq \Bl_{\Gamma_1 \cap U_0 } U_0 $ and let $T^{\bg,(1)}_{i,j} = T^{(1)}_{i,j}$.
\item[(2)] For $v \geq 1$, the strict transform of $\Gamma_{v+1}$ on $U_{v}^{\bg}$ is cut out by
$$
\left(T_{i,j}^{\bg,(v)}-T^{\bg,(v)}_{i,\bg^2_v}T_{\bg^1_v,j}^{\bg,(v)}\right)^{i \in \{0,\dots,k-1\} \setminus \{\bg_1^1,\dots,\bg_{v}^1\}}_{j \in \{k,\dots,n\} \setminus \{\bg_1^2,\dots,\bg_{v}^2\}}
$$

\item[(3)] For $v \geq 1$, the locus $\Gamma_{v+1} \cap U_v^{\bg}$ is non-singular and 
$$
\Bl_{\Gamma_{v+1}\cap U_v^{\bg}} U_v^{\bg} \simeq \Proj \kk[U_v^{\bg}][T^{\bg,(v)}_{i,j}]_{i,j}/(\text{Koszul Relations}).
$$
\item[(4)] For $v \geq 1$, we have $U_v^{\bg} = D(T^{\bg,(v)}_{\bg^1_v,\bg^2_v})  \subseteq \Bl_{\Gamma_v \cap U_{v-1}^{\bg}} U_{v-1}^{\bg}$.
\end{enumerate}
Due to symmetry, the proof of Proposition \ref{bigmatrix} also establishes the above statements (note that $U_{k-1} = U^{\bg}_{k-1}$ with $\bg^1 = (k-1,k-2,\dots,0)$ and  $\bg^2 = (n,n-1,\dots,n-k+2)$). It follows that  $\{U_{k-1}^{\bg}\}_{\bg}$ is an affine cover of $\pi^{-1}_{k-1}(U_0)$ with the natural gluing maps. We omit an explicit description of the gluing maps as they will never be used. 

\smallskip
To construct the $U_{v}^{\bg}$ and verify statement (2), we would have to row reduce $M$ in a way analogous to Proposition \ref{bigmatrix} (each $\bg$ corresponds to a different sequence of row redutions). We will omit an explicit description of the matrix, but the corresponding lambdas are
$$
\lambda_{1}^{\bg} = a_{\bg^1_1,\bg^2_1} -b_{\bg^1_1,\bg^2_1} \quad \text{and} 
	\quad \lambda^{\bg}_i = T^{\bg,(i-1)}_{\bg^1_i,\bg^2_i}- T^{\bg,(i-1)}_{\bg^1_i,\bg^2_{i-1}}T^{\bg,(i-1)}_{\bg^1_{i-1},\bg^2_{i}} 
	\quad \text{for each } \, 2 \leq i \leq k-1.
$$
As in the proof of Proposition \ref{mainone} we can choose a minimal generating set, $\alpha^{\bg}_0,\dots,\alpha^{\bg}_{k-1}$ of  $I_{\Lambda(\mathbf{a})}$ where
\begin{align*}
\alpha_{p}^{\bg} = y_{\bg^1_{k-p}}-\sum_{j=1}^{k-1-p} T_{\bg^1_{k-p},\bg^2_j}^{\bg,(j)}y_{\bg^1_{j}}
			+	\sum_{j \in \{k,\dots,n\}\setminus \{\bg^2_1,\dots,\bg^2_{k-1-p}\}}
									\lambda_{1}^{\bg}\cdots\lambda^{\bg}_{k-p}T_{\bg^1_{k-p},j}^{\bg,(k-p)}x_{j} 
\end{align*}
for $0 < p \leq k-1$ and
\begin{equation*}
\alpha_{0}^{\bg} = y_{\bg^1_{k}}-\sum_{j=1}^{k-1} T_{\bg^1_{k},\bg^2_j}^{\bg,(j)}y_{\bg^1_{j}}
			+	\sum_{j \in \{k,\dots,n\}\setminus \{\bg^2_1,\dots,\bg^2_{k-1}\}}
									\lambda_{1}^{\bg}\cdots\lambda^{\bg}_{k-1}T_{\bg^1_{k},j}^{\bg,(k)}x_{j} 
\end{equation*}
with $T^{\bg,(k)}_{\bg^1_k,j} = T_{\bg^1_k,j}^{\bg,(k-1)}-T_{\bg^1_k,\bg^2_{k-1}}^{\bg,(k-1)}T_{\bg^1_{k-1},j}^{\bg,(k-1)}$. 

\smallskip
For $0\leq  p<q\leq k-1$  we may define analogous "cross terms" 
\begin{align*}
\beta_{p,q}^{\bg} &= 
	\left(y_{\bg^1_{k-p}}-\sum_{j=1}^{k-1-p} T_{\bg^1_{k-p},\bg^2_j}^{\bg,(j)}y_{\bg^1_{j}}\right)
		\left(\sum_{j \in \{k,\dots,n\}\setminus \{\bg^2_1,\dots,\bg^2_{k-1-q}\}}T_{\bg^1_{k-q},j}^{\bg,(k-q)}x_{j}\right)
		\\ 
		& \quad\quad\quad - \lambda^{\bg}_{p,q}
			\left(y_{\bg^1_{k-q}}-\sum_{j=1}^{k-1-q} T_{\bg^1_{k-q},\bg^2_j}^{\bg,(j)}y_{\bg^1_{j}}\right)
		\left(\sum_{j \in \{k,\dots,n\}\setminus \{\bg^2_1,\dots,\bg^2_{k-1-p}\}}T_{\bg^1_{k-p},j}^{\bg,(k-p)}x_{j}\right).
\end{align*}

Thus we obtain a morphism 
\begin{equation}
\label{finalmap2}
\Xi_{U^{\bg}_{k-1}}:(\mathbf{a},\mathbf{b},\TT^{\bg,(1)},\dots,\TT^{\bg,(k)}) \mapsto \left[ I_{\Lambda(\mathbf{a})}I_{\Lambda(\mathbf{b})}+ (\beta_{p,q}^{\bg})_{0 \leq p < q \leq k-1} \right].
\end{equation}
This is well defined as any ideal in the image of $\Xi_{U^{\bg}_{k-1}}$ is still projectively equivalent to an ideal in (\ref{projideal}) (the proof of Lemma \ref{projequiv} works with straightforward modifications). As explained in Proposition \ref{mainone}, $\Xi_{U^{\bg}_{k-1}}$ will also extend the original rational map (\ref{originalXi}) for each $\bg$. Thus  for any $\bg,\bg'$, $\Xi_{U^{\bg}_{k-1}}$ and $\Xi_{U^{\bg'}_{k-1}}$ agree on an open subset of $U^{\bg}_{k-1} \cap U^{\bg'}_{k-1}$. By uniqueness of extensions, they will agree on all of $U^{\bg}_{k-1} \cap U^{\bg'}_{k-1}$ . Gluing all these maps gives us a morphism $\pi^{-1}_{k-1}(U_0) \longrightarrow \HH_{n-k,n-k}^n$. 

As mentioned in the beginning of the section, $\Gr(n-k,n)^2$ is covered by open sets of the form $U_{\mathcal{E}}$ where $\mathcal{E}$ ranges over all ordered bases of $S_1$. Since assuming $\mathcal{E} = \{x_0,\dots,x_{n}\}$ was purely notational, all the discussion in this section applies verbatim to $\pi_{k-1}^{-1}(U_{\mathcal{E}})$. In particular, we obtain a morphism on each $\pi_{k-1}^{-1}(U_{\mathcal{E}})$ that extends the original rational map (\ref{originalXi}). Thus we can glue all these maps to obtain a morphism $\Xi: \XX_{k-1} \longrightarrow \HH_{n-k,n-k}^n$.

Let $\mathfrak{S}_2=\{1,g\}$ be the group on two elements and consider its natural on $\mathbf{Gr}(n-k,n)^2$ given by interchanging the two factors. Since each of the $\Gamma_i$ are $\mathfrak{S}_2$ stable, the action extends to the blowup $\XX_{k-1}$. If we consider the trivial action of $\mathfrak{S}_2$ on $\HH_{n-k,n-k}^n$, then our construction shows that $\Xi$ is $\mathfrak{S}_2$-equivariant. Thus, we get an induced morphism $\widebar{\Xi}: \XX_{k-1}/\mathfrak{S}_2 \longrightarrow \HH_{n-k,n-k}^n$.

Since $\mathrm{char} \, \kk \ne 2$ and $g$ fixes a divisor (the strict transform of the exceptional divisor of $\XX_1$), the Chevalley-Shephard-Todd theorem \cite[Theorem 7.14]{invarianttheory} implies that the quotient is non-singular. Note that
$$
\XX_{k-1}/\mathfrak{S}_2 
			= (\Bl_{\Gamma_{k-1}} \cdots \Bl_{\Gamma_1} \mathbf{Gr}(n-k,n)^2)/\mathfrak{S}_2 
			\simeq \Bl_{\widebar{\Gamma}_{k-1}} \cdots \Bl_{\widebar{\Gamma}_1} \Sym^2 \mathbf{Gr}(n-k,n).
$$
Since $\Xi$ is dominant and $\XX_{k-1}$ is projective, $\widebar{\Xi}$ is surjective. 

The natural action of $\text{GL}(n+1)$ on $\mathbf{P}^n$ induces an action on $\Gr(n-k,n)^2$ and on $\HH_{n-k,n-k}^n$. Since the $\Gamma_i$ are stable under this action, it extends to an action on $\XX_{k-1}$. To show that $\Xi$ is $\GL(n+1)$-equivariant we need to show that for any $g\in \text{GL}(n+1)$ the two morphisms, $\Xi \circ g: \XX_{k-1} \to \HH_{n-k,n-k}^n $ given by $w \mapsto \Xi(gw)$  and $g \circ \Xi : \XX_{k-1} \to \HH_{n-k,n-k}^n $ given by $w \mapsto g\Xi(w)$ are identical.  For any $(\Lambda,\Lambda')$ in the open set  $\Gr(n-k,n)^2 \setminus \Gamma_k \subseteq \XX_{k-1}$ we have 
$$
(\Xi \circ g)(\Lambda,\Lambda') = \Xi(g(\Lambda), g(\Lambda'))= g(\Lambda) \cup g(\Lambda') = g(\Lambda \cup \Lambda') =  (g \circ \Xi)(\Lambda,\Lambda').
$$
Thus $\Xi \circ g$ and $g \circ \Xi$ must agree on all of $\XX_{k-1}$. It follows that $\widebar{\Xi}$ is also $\text{GL}(n+1)$-equivariant.
\end{proof}

\begin{cor} \label{2kgens} Let $n \geq 2k-1$. Any subscheme parameterized by $\mathcal{H}_{n-k,n-k}^n$ is minimally cut out by $k^2$ quadrics.
\end{cor}
\begin{proof} 
By the discussion in Proposition \ref{mainonee} we may reduce to considering subschemes cut out by ideals in the image of morphism (\ref{finalmap1}). Let $J$ denote any such ideal and note that $J$, as presented, is generated by quadrics. By Lemma \ref{saturated} (i), $J$ is saturated and thus is the ideal of its corresponding subscheme. Therefore it suffices to show that $\dim_{\kk} J_2 = k^2$. Since $S/J$ and $S/\inn_{>} J$ have the same Hilbert function we have $ \dim_{\kk} J_2 = \dim_{\kk} (\inn_{>} J)_2 = k^2$ (Lemma \ref{GROBNER}).
\end{proof}

\begin{remark} \label{reducetoUk} The analogue of Lemma \ref{anothermingens} holds verbatim for ideals in the image of Equation (\ref{finalmap2}). The analogue of Lemma \ref{GROBNER} is as follows: Let $J$ be any ideal in the image of Equation (\ref{finalmap2}) and let $>_{\bg}$ denote a lexicographic ordering on $S$ for which 
$$
x_{\bg_{k}^1} > x_{\bg_{k-1}^1}> \cdots > x_{\bg^1_1} > x_{\bg^2_1} > \cdots >x_{\bg^2_{k-1}} > x_{h_1} > \cdots > x_{h_{n-2k+2}}.
$$
We may choose any $h_i$ so that $\{h_1,\dots,h_{n-2k+2}\} = \{k,\dots,n\} \setminus\{\bg^2_1,\dots,\bg^2_{k-1}\}$. Then we have
$$
\inn_{>_{\bg}} J = (x_0,\dots,x_{k-1})^2 + (x_{\bg^1_{k-p}}x_{\bg^2_{k-q}})_{0 \leq p < q \leq k-1}
$$
\end{remark}

We split the proof of the injectivity of $\widebar{\Xi}$ into two steps. Here is the first step.

\begin{lemma} \label{injective} For any $\bg$, the restriction $\widebar{\Xi}: U_{k-1}^{\bg}/\mathfrak{S}_2 \longrightarrow \HH_{n-k,n-k}^n$ is injective.
\end{lemma}
\begin{proof} It is evident from our construction that $U_{k-1}^{\bg}$  is $\mathfrak{S}_2$-stable and thus the quotient $U^{\bg}_{k-1}/\mathfrak{S}_2$ is well defined. Without loss of generality we may assume $U_{k-1}^{\bg}=U_{k-1} $. To prove the Lemma it suffices to show that for any $\tZ,\hZ \in U_{k-1}$ satisfying $\Xi(\tZ) = \Xi(\hZ)$, we have $\tZ = \hZ$ or $g(\tZ) = \hZ$ where where $g$ is the non-identity of $\mathfrak{S}_2$. Let $\tZ = (\tilde{\mathbf{a}}, \tilde{\mathbf{b}},\tilde{\mathbf{T}}^{(1)},\dots,\tilde{\mathbf{T}}^{(k)}) $ and $\hZ = (\hat{\mathbf{a}},\hat{\mathbf{b}},\hat{\mathbf{T}}^{(1)},\dots,\hat{\mathbf{T}}^{(k)})$ be their coordinates on $U_{k-1}$. The "betas" and "lambdas" corresponding to $\tZ$ are denoted by $\tilde{\beta}_{i,j}$ and $\tilde{\lambda}_i$ respectively, and the ones corresponding to $\hZ$ are denoted by $\hat{\beta}_{i,j}$ and $\hat{\lambda}_i$.

We have 
$
\Lambda(\tilde{\mathbf{a}}) \cup \Lambda(\tilde{\mathbf{b}})= \Xi(\tZ)_{\text{red}}=\Xi(\hZ)_{\text{red}} = \Lambda(\hat{\mathbf{a}}) \cup \Lambda(\hat{\mathbf{b}})$. 
After possibly replacing $\tZ,\hZ$ by $g(\tZ), g(\hZ)$ respectively, we may assume  
$\tilde{\mathbf{a}} = \hat{\mathbf{a}}$ and $\tilde{\mathbf{b}}= \hat{\mathbf{b}}$. Thus to prove that $\widebar{\Xi}$ is injective, we need to now show that $\tZ= \hZ$. Since $\Xi$ is $\GL(n+1)$-equivariant we may apply a projective transformation and assume $\tilde{\mathbf{b}}= \hat{\mathbf{b}} = \mathbf{0}$. For simplicity we let $\aaa := \tilde{\mathbf{a}}= \hat{\mathbf{a}}$.

By Lemma \ref{anothermingens}, $\Xi(\tZ)_{\text{red}}=\Xi(\hZ)_{\text{red}}$ is a pair of $(n-k$)-planes meeting along an $(n-k+1-\ell)$-plane for some $1 \leq \ell \leq k+1$. If $\ell \in \{k,k+1\}$ then $\widetilde{Z},\widehat{Z}$ lie in an open set along which $\Xi$ was already shown to be two-to-one (Lemma \ref{twotoone}). Thus we may assume $\ell \leq k-1$. By Lemma \ref{anothermingens} it is also the smallest index for which $\tilde{\lambda}_{\ell} =0$ and, symmetrically, the smallest index for which $\hat{\lambda}_{\ell} = 0$.

Using Lemma \ref{anothermingens} (iv) we get $\Xi(\tZ) = [I_{\Lambda(\mathbf{a})}I_{\Lambda(\mathbf{0})}+(\tilde{\beta}_{p,q})_{0\leq p<q\leq k-\ell}]$ and $\Xi(\hZ) = [I_{\Lambda(\mathbf{a})}I_{\Lambda(\mathbf{0})}+(\hat{\beta}_{p,q})_{0\leq p<q\leq k-\ell}]$. Using Lemma \ref{anothermingens} (i) we have the equality
$$
I_{\Lambda(\mathbf{a})}I_{\Lambda(\mathbf{0})}+(\tilde{\beta}_{p,q})_{0\leq p<q\leq k-\ell}
	= I_{\Lambda(\mathbf{a})}I_{\Lambda(\mathbf{0})}+(\hat{\beta}_{p,q})_{0\leq p<q\leq k-\ell}.
$$
I claim that $(\tilde{\beta}_{p,q})_{0 \leq p < q \leq k- \ell} = (\hat{\beta}_{p,q})_{0 \leq p < q \leq k- \ell}$. Assume $\tilde{\beta}_{p,q} = \alpha + \omega$ with $\alpha \in I_{\Lambda(\mathbf{a})}I_{\Lambda(\mathbf{0})}$ and $\omega \in (\hat{\beta}_{p,q})_{0 \leq p < q \leq k- \ell}$ such that $\alpha,\omega$ are linearly independent and homogenous of degree $2$. Since $\hat{\lambda}_{\ell} = \tilde{\lambda}_{\ell} = 0$, the construction in Proposition \ref{mainone} implies
\begin{align*}
I_{\Lambda(\mathbf{a})}I_{\Lambda(\mathbf{0})}  =  (\alpha_0,\dots,\alpha_{k-1})(x_0,\dots,x_{k-1})
					\subseteq   (x_0,\dots,x_{k-1},x_{n-\ell+2},\dots,x_n)(x_0,\dots,x_{k-1}) 
\end{align*}
and
$$
(\tilde{\beta}_{p,q})_{0 \leq p < q \leq k- \ell} , (\hat{\beta}_{p,q})_{0 \leq p < q \leq k- \ell} \subseteq (x_0,\dots,x_{k-1})(x_k,\dots,x_{n-\ell+1}).
$$
This implies $\alpha = 0$ and we obtain $B = (\tilde{\beta}_{p,q})_{0 \leq p < q \leq k- \ell} = (\hat{\beta}_{p,q})_{0 \leq p < q \leq k- \ell}$. The proof will be complete once we the show that the coordinates from Remark \ref{coordinates} of $\widetilde{Z}$ coincide with those of $\widehat{Z}$.

\smallskip
 It follows from the proof of Proposition \ref{bigmatrix} that the coordinate $T^{(v)}_{i,j}$ admits a formal expression
\begin{equation} \label{Ts}
T^{(v)}_{i,j} = \frac{A_{i,j,v}(\mathbf{a},\bb,\lambda_1,\dots,\lambda_v)}{\lambda_{1}^{\epsilon_1}\cdots\lambda_{v}^{\epsilon_{v}}}
\end{equation}
with $A_{i,j,v}$ a polynomial in $\mathbf{a},\bb,\lambda_1,\dots,\lambda_v$ and $\epsilon_1,\dots,\epsilon_{v} \geq 1$. Similarly, each $\lambda_v$ admits a formal expression
\begin{equation} \label{lambdas}
\lambda_{v} = \frac{B_{i,j,v}(\mathbf{a},\bb,\lambda_1,\dots,\lambda_{v-1})}{\lambda_{1}^{\epsilon_1}\cdots\lambda_{v-1}^{\epsilon_{v-1}}}
\end{equation}
with  $B_{i,j,v}$ a polynomial in $\mathbf{a},\bb,\lambda_1,\dots,\lambda_{v-1}$ and $\epsilon_1,\dots,\epsilon_{v-1} \geq 1$.
\begin{enumerate}
\item $\hat{\lambda}_i = \tilde{\lambda}_i$ for all $i \leq \ell$: We clearly have $\hat{\lambda}_1 = a_{k-1,n} = \tilde{\lambda}_1$. Since $\hat{\lambda}_v \ne 0$ for all $v \leq \ell -1$ we can inductively apply (\ref{lambdas}) to obtain
\begin{equation*}
\hat{\lambda}_{v} 
			= \frac{B_{i,j,v}(\mathbf{a},\mathbf{0},\hat{\lambda}_1,\dots,\hat{\lambda}_{v-1})}{\hat{\lambda}_{1}^{\epsilon_1}\cdots\hat{\lambda}_{v-1}^{\epsilon_{v-1}}} 
			= \frac{B_{i,j,v}(\mathbf{a},\mathbf{0},\tilde{\lambda}_1,\dots,\tilde{\lambda}_{v-1})}{\tilde{\lambda}_{1}^{\epsilon_1}\cdots\tilde{\lambda}_{v-1}^{\epsilon_{v-1}}} 
			= \tilde{\lambda}_{v}.
\end{equation*}
\item $\hat{T}^{(v)}_{i,j} = \tilde{T}^{(v)}_{i,j}$ for all $v \leq \ell-1$  and all $i,j$: Analogous to item (i) above, where we instead use (\ref{Ts}) to conclude 
$$\hat{T}^{(v)}_{i,j} = 
					\frac{A_{i,j,v}(\mathbf{a},\mathbf{0},\hat{\lambda}_1,\dots,\hat{\lambda}_v)}{\hat{\lambda}_1^{\epsilon_1}\cdots\hat{\lambda}_{v}^{\epsilon_{v}}} 
					= \frac{A_{i,j,v}(\mathbf{a},\mathbf{0},\tilde{\lambda}_1,\dots,\tilde{\lambda}_v)}{\tilde{\lambda}_{1}^{\epsilon_1}\cdots\tilde{\lambda}_{v}^{\epsilon_{v}}}  
					= \tilde{T}^{(v)}_{i,j}.
$$
\item $\hat{T}^{(v)}_{i,j} = \tilde{T}^{(v)}_{i,j}$ for all $k-1 \geq v \geq \ell$ and all relevant $i,j$ (those appearing as coordinates in Remark \ref{coordinates}): Let $r,s$ be any integers such that $0 \leq r < s \leq k-\ell$ and assume $\hat{\beta}_{r,s} = \sum_{0 \leq p < q \leq k-\ell}c_{p,q}\tilde{\beta}_{p,q}$ for some constants $c_{p,q} \in \kk$. 
Let $p' = \min\{p:c_{p,q} \ne 0\}$ and $q' = \max\{q:c_{p',q} \ne 0\}$. 
Then 
$$
x_rx_{n-k_s} = \inn_{>}(\hat{\beta}_{r,s}) = \inn_{>}\left(\sum_{0 \leq p < q \leq k-\ell}c_{p,q}\tilde{\beta}_{p,q}\right) = c_{p',q'}x_{p'}x_{n-k_{q'}}.$$ 
It follows that $\tilde{\beta}_{r,s} = \hat{\beta}_{r,s}$. Equating the terms supported on $x_r$ we obtain
$$
\sum_{j=k}^{n-k_s}\hat{T}^{(k-s)}_{s,j}x_j = \sum_{j=k}^{n-k_s}\tilde{T}^{(k-s)}_{s,j}x_j.
$$ 
It follows that $\hat{T}_{s,j}^{(k-s)} = \tilde{T}_{s,j}^{(k-s)}$ for all $k \leq j < n-k_s$. Similarly, equating the terms supported on $x_{n-k_s}$ we obtain $\hat{T}^{(j)}_{r,n-j+1} = \tilde{T}^{(j)}_{r,n-j+1}$ for all $1 \leq j \leq k_r.$

\smallskip
\item $\hat{T}^{(k)}_{0,j} = \tilde{T}^{(k)}_{0,j}$ for all $k \leq j \leq n-k+1$: Combining $\hat{\beta}_{0,1} = \tilde{\beta}_{0,1}$ and the equality of coordinates in (iii) we obtain
{\small
$$
\hat{\lambda}_{0,1}\left(x_{1}- \sum_{j=1}^{k-2}\hat{T}_{1,n-j+1}^{(j)}x_{k-j}\right)\left(\sum_{j=k}^{n-(k-1)}\hat{T}_{0,j}^{(k)}x_{j}\right)
=\tilde{\lambda}_{0,1}\left(x_{1}- \sum_{j=1}^{k-2}\tilde{T}_{1,n-j+1}^{(j)}x_{k-j}\right)\left(\sum_{j=k}^{n-(k-1)}\tilde{T}_{0,j}^{(k)}x_{j}\right).
$$
}
Since $\hat{\lambda}_{0,1} = 1 = \tilde{\lambda}_{0,1}$, equating the coefficients of the monomials containing $x_1$ gives the desired result.
\smallskip
\item $\hat{\lambda}_i = \tilde{\lambda}_i$ for all $i \geq \ell+1$: For each $\ell+1 \leq i \leq k-1$ we have $\tilde{\beta}_{k-i,k-i+1} =   \hat{\beta}_{k-i,k-i+1}$. Note that $\hat{\lambda}_{k-i,k-i+1} = \hat{\lambda}_i$ and $\tilde{\lambda}_{k-i,k-i+1} = \tilde{\lambda}_i$. Using the equality of coordinates in (iii), the expression $\tilde{\beta}_{k-i,k-i+1} =   \hat{\beta}_{k-i,k-i+1}$ reduces to
{\small
\begin{equation*} 
\hat{\lambda}_{i}\left(x_{k-i+1}- \sum_{j=1}^{i-2}\hat{T}_{k-i+1,n-j+1}^{(j)}x_{k-j}\right)\left(\sum_{j=k}^{n-i+1}\hat{T}_{k-i,j}^{(i)}x_{j}\right) =
\tilde{\lambda}_{i}\left(x_{k-i+1}- \sum_{j=1}^{i-2}\tilde{T}_{k-i+1,n-j+1}^{(j)}x_{k-j}\right)\left(\sum_{j=k}^{n-i+1}\tilde{T}_{k-i,j}^{(i)}x_{j}\right).
\end{equation*}
}
Equating the coefficients of $x_{k-i+1}x_{n-i+1}$ gives the desired result.\qedhere
\end{enumerate} 
\end{proof}

\begin{lemma} \label{Binjective} The fiber of $\Xi$ over the point  $[(x_0,\dots,x_{k-1})^2 +(x_px_{n-k_q})_{0 < p <q \leq k-1}]$ consists of a single element.
\end{lemma}
\begin{proof} Let $J$ denote the ideal $(x_0,\dots,x_{k-1})^2 +(x_px_{n-k_q})_{0 < p <q \leq k-1}$. Let $X \in U_{k-1}$ be the point with all the coordinates of Remark \ref{coordinates} equal to $0$. We clearly have $\Xi(X) = [J]$. Now assume $Z \in \XX_{k-1}$ such that $\Xi(Z) = [J]$. Since $J_{\text{red}} = (x_0,\dots,x_{k-1})$, we must have $Z \in \pi^{-1}_{k-1}(U_0)$. In particular, $Z \in U_{k-1}^{\bg}$ for some $\bg$. By Remark \ref{reducetoUk} we have
$$
(x_0,\dots,x_{k-1})^2 + (x_{\bg^1_{k-p}}x_{\bg^2_{k-q}})_{0 \leq p < q \leq k-1} = \inn_{>_{\bg}} \Xi(Z) = \inn_{>_{\bg}}J = J.
$$
Comparing the monomial generators of the two ideals we deduce that $\bg^1_{k-p} = p$ for all $0 \leq p \leq k-2$; this forces $\bg^1_1 = k-1$. But then we also obtain $\bg^2_{k-q} = n-k_q= n - (k-q)+1$ for all $1 \leq q \leq k-1$. Thus $U_{k-1}^{\bg} = U_{k-1}$ and by Lemma \ref{injective}, $Z = X$ or $g(Z) = X$ for the non-identity $g \in \mathfrak{S}_2$. Since $\Xi(Z)_{\text{red}} = \Xi(X)_{\text{red}} = V(x_0,\dots,x_{k-1})$ we must have $g(Z)=Z$; thus $Z= X$.
\end{proof}

\begin{prop} \label{INJECTIVE} Let $n \geq 2k-1$. The morphism $\widebar{\Xi}: \XX_{k-1}/\mathfrak{S}_2 \longrightarrow \HH_{n-k,n-k}^n$ is injective.
\end{prop}
\begin{proof} Let $Y,Z \in \mathcal{X}_{k-1}$ such that $\Xi(Y) = \Xi(Z)$. Since $\Xi(Y)_{\text{red}}=\Xi(Z)_{\text{red}}$ we may assume wlog that $Y,Z \in \pi^{-1}_{k-1}(U_0)$. We may also assume wlog that $Y \in U_{k-1}$. By Lemma \ref{injective} we only need to show that $Z \in U_{k-1}$. Let $\ell \geq 1$ be the maximal value such that $Z \in U_{k-1}^{\bg}$ with $\bg^1_i = k-i$ and $\bg^2_i = n-i+1$ for all $i < \ell$. We need to show that $\ell = k$ (then automatically, $\bg^1_k = 0$). For the sake of a contradiction, assume that $\ell < k$.  Our method is to compare certain initial ideal degenerations of $\Xi(Z)$ and $\Xi(Y)$.

\smallskip
Let $\mathbf{w}$ be any integral weight order corresponding to $>$ \cite[Section 15]{eisenbud}. For any $t \in \kk^{\star}$ let $g_t \in \GL(n+1)$ denote the automorphism that maps $x_i \mapsto t^{-\mathbf{w}(i)}x_i$.  Since each $g_t$ just scales the coordinates the following facts are immediate
\begin{itemize}
\item[(1)] $g_t$ induces an action on $\XX_0$ and extends to all the blowups $\XX_v$.
\item[(2)] $g_t$ fixes $U_{\ell}^{\bg}$ and also fixes any closed subset of the form $V(T^{\bg,(\ell)}_{i,j})$.
\item[(3)] For each $\ell$ let $\psi_{\ell}: \XX_{k-1} \longrightarrow \XX_{\ell}$ denote the blowdown map. Then $\psi_{\ell}$ is $\GL(n+1)$-equivariant and thus $\psi_{\ell}(g_t) = g_t(\psi_{\ell})$.
\end{itemize}

Let $Y_0 = \lim_{t\to 0} g_t(Y)$ and $Z_0 = \lim_{t\to 0} g_t(Z)$. Using \cite[Theorem 15.17]{eisenbud} and Lemma \ref{GROBNER} we obtain 
$$
 \Xi(Y_0) = \lim_{t \to 0} g_t(\Xi(Y)) = \inn_{>} \Xi(Y) = (x_0,\dots,x_{k-1})^2+(x_px_{n-k_q})_{0<p< q \leq  k-1}.
$$
Similarly, $\Xi(Z_0) = (x_0,\dots,x_{k-1})^2+(x_px_{n-k_q})_{0<p< q \leq  k-1}= \Xi(Y_0)$. By Lemma \ref{Binjective}, $Z_0 = Y_0$.

\smallskip
Using the notation in item (3) and our assumption on $\ell$, $\psi_{\ell}(Z)$ and $\psi_{\ell}(Y)$ are $\kk$-points of $\Proj \kk[U_{\ell-1}][T^{(\ell)}_{i,j}]/(\text{Koszul}) \subseteq \XX_{\ell}$. By maximality of $\ell$ we have $T^{(\ell)}_{k-\ell,n-\ell+1}(\psi_{\ell}(Z)) = 0$ i.e. $\psi_{\ell}(Z)$ lies in $V(T^{(\ell)}_{k-\ell,n-\ell+1})$. Then by item (2) we still have $\psi_{\ell}(g_t(Z)) = g_t(\psi_{\ell}(Z)) \in V(T^{(\ell)}_{k-\ell,n-\ell+1})$. Thus the limit $\psi_{\ell}(Z_0)$ also lies in there. But this contradicts the fact that $T^{(\ell)}_{k-\ell,n-\ell+1}(\psi_{\ell}(Y_0)) = T^{(\ell)}_{k-\ell,n-\ell+1}(Y_0) \ne 0$ (since $Y_0$ lies in $U_{k-1}$). Thus $\ell = k$ and we have $Z,Y \in U_{k-1}$, as required.
\end{proof}

\begin{remark} \label{uniquepoint} It follows that the preimage $\Xi^{-1}(Z)$ is a single point precisely when $Z_{\text{red}}$ is an $(n-k)$-plane. This occurs precisely when $Z$ is generically non-reduced, c.f. Theorem \ref{primdecomp} \footnote{If the reader is only interested in the classification of subschemes parameterized by $\HH_{n-k,n-k}^n$ they can directly skip to Lemma \ref{CM}}.
\end{remark}


\smallskip
The group $\GL(n+1)$ acts on $S$ and thus on $\Hilb^{P(t)} \, \mathbf{P}^n$ by a change of coordinates.  An ideal of $S$ or its corresponding point on the Hilbert scheme is said to be \textbf{Borel fixed} if it is fixed by the Borel subgroup of $\GL(n+1)$ consisting of upper triangular matrices. Since a Borel fixed ideal is fixed by the subgroup of diagonal matrices, it is generated by monomials. We will now show that $\HH_{n-k,n-k}^n$ has a unique Borel fixed point. We begin with a combinatorial characterization of the Borel fixed ideals, see \cite[Section 15]{eisenbud} for details.

\begin{definition} \label{boreldef} Let $I \subseteq S$ be a monomial ideal and $p$ a prime number. The ideal $I$ is said to be $0$-\textbf{Borel fixed} if for any monomial generator $m\in I$ divisible by $x_j$, we have $\frac{x_i}{x_j}m \in I$ for all $i<j$. The ideal $I$ is said to be $p$-\textbf{Borel fixed} if for any monomial generator $m\in I$ divisible by $x_j^\beta$ but no higher power of $x_j$, we have $(\frac{x_i}{x_j})^{\alpha}m \in I$ for all $i<j$ and $\alpha \preceq_{p} \beta$ (this means that each digit in the $p$-base expansion of $\alpha$ is less than or equal to each digit in the $p$-base expansion of $\beta$).
\end{definition}

Note that a $0$-Borel fixed ideal is always $p$-Borel fixed for any $p$.

\begin{prop}\emph{\cite[Theorem 15.23]{eisenbud}} Let  $ \mathrm{char} \, \kk = p\geq 0$. Then $I \subseteq S$ is Borel fixed if and only if it $I$ is $p$-Borel.
\end{prop}

In our situation, $\mathrm{char} \, \kk = p \geq 0$ with $p \ne 2$. Let $I$ be a saturated $p$-Borel fixed ideal parameterized by $\HH_{n-k,n-k}^n$. Since $I$ is a monomial ideal generated by quadrics (Corollary \ref{2kgens}) and $p \ne 2$, the condition $\alpha \preceq_{p} \beta$ in Definition \ref{boreldef} reduces to the condition $\alpha \leq \beta$. In particular, $I$ is always $0$-Borel.

\begin{prop} {\label{uniqueborel}} Let $n\geq 2k-1$. Consider the ideal
$$
I_{n-k,n-k}^n = \sum_{i=0}^{k-1}x_i(x_i,\dots,x_{2k-2-i}) = (x_0,\dots,x_{k-1})^2 + (x_px_{2k-1-q})_{0 \leq p < q \leq k-1}.
$$
Then $[I_{n-k,n-k}^n]$ is the unique Borel fixed point on $\mathcal{H}_{n-k,n-k}^n$.
\end{prop}
\begin{proof} As noted above, Borel fixed ideals in $\HH_{n-k,n-k}^n$ are the same as $0$-Borel fixed ideals. Since $I_{n-k,n-k}^n$ is projectively equivalent to $(x_0,\dots,x_{k-1})^2 + (x_px_{n-k_q})_{0\leq p < q \leq k-1}$, it lies in $\mathcal{H}_{n-k,n-k}^n$.  It also clear that $I_{n-k,n-k}^n$ is Borel fixed.   Let $B$ be any saturated $0$-Borel fixed ideal on $\mathcal{H}_{n-k,n-k}^n$. Then it is of the form $B = \sum_{i=0}^{\epsilon}x_i(x_i,\dots,x_{a_i})$ with $n-1 \geq a_0 \geq a_1\geq \cdots \geq a_\epsilon \geq \epsilon$. Since $\sqrt{B} = (x_0,\dots,x_{\epsilon})$ has codimension $k$, we obtain $\epsilon=k-1$.

Arguing as in the end of the proof of Proposition \ref{mainone} we see that the Hilbert polynomial of $B$ is $\binom{n-k+t}{t} + \sum_{i=0}^{k-1}\binom{t+n-a_i-2}{t-1}$. Equating this with the Hilbert polynomial of $I_{n-k,n-k}^n$ we have
$$
\sum_{i=0}^{k-1}\binom{n-2k+i+t}{t-1} = \sum_{i=0}^{k-1}\binom{t+n-a_i-2}{t-1}.
$$
Since the set $\{\binom{t-1+a}{a}\}_{a\in \mathbf{N}}$ is a $\mathbf{Q}$-basis for $\mathbf{Q}[t]$, we obtain $a_i = 2k-i-2$ for all $i$; therefore $B = I_{n-k,n-k}^n$.
\end{proof}

\begin{lemma} \label{regularitythm} Let $I$ be a (saturated) ideal parameterized by $\HH_{n-k,n-k}^n$. Then the Castelnuovo-Mumford regularity of $I$ is $2$ and $T_{[I]}\Hilb^{P_{n-k,n-k}^n(t)} \, \mathbf{P}^n = \Hom_S(I,S/I)_0$.
\end{lemma}
\begin{proof} Since $I$ is generated by quadrics, the regularity is at least $2$. Up to projective equivalence, we may assume $I$ is of the form (\ref{projideal}). By \cite[Theorem 3.3.4]{monomialideals} we have also $\mathrm{reg}(I) \leq \mathrm{reg}(\inn_{>}I)$. Note that $\inn_{>}I$ is projectively equivalent to $I_{n-k,n-k}^n$ and the regularity of a $0$-Borel ideal is the highest degree of a minimal monomial generator \cite[Corollary 7.2.3]{monomialideals}. Thus $\mathrm{reg}(I) \leq \mathrm{reg}(I_{n-k,n-k}^n) = 2$, as required.  The description of the tangent space follows from Remark \ref{applycomparison} and Theorem \ref{compare}.
\end{proof}

\begin{definition} \label{zetazeta} Let $\zeta$ denote the pre-image of $[I_{n-k,n-k}^n]$ in $\XX_{k-1}$ (Remark \ref{uniquepoint}) and let $\bar{\zeta}$ denote the image of $\zeta$ in $\XX_{k-1}/\mathfrak{S}_2$.
\end{definition}

By constructing curves passing through $\zeta$ and $\bar{\zeta}$ we will now show that the differential $d\widebar{\Xi}_{\bar{\zeta}}$ is injective. This is a major portion of the proof of Theorem \ref{main}.

\begin{lemma} \label{dacurves} Let $n \geq 2k-1$. The differential $d\widebar{\Xi}_{\bar{\zeta}}:T_{\bar{\zeta}} (\XX_{k-1}/\mathfrak{S}_2) \longrightarrow T_{[I_{n-k,n-k}^n]} \HH_{n-k,n-k}^n$ is injective.
\end{lemma}
\begin{proof}
Note that we have a factorization
$$
\begin{tikzcd}
T_{\zeta} \XX_{k-1} \arrow[r] \arrow[dr]
& T_{\bar{\zeta}} (\XX_{k-1}/\mathfrak{S}_2) \arrow[d]\\
& T_{[I_{n-k,n-k}^n]} \HH_{n-k,n-k}^n
\end{tikzcd}
$$
By non-singularity we also have $\dim_{\kk}T_{\zeta} \XX_{k-1} = \dim_{\kk}T_{\bar{\zeta}} (\XX_{k-1}/\mathfrak{S}_2)$. Thus to show that $d\widebar{\Xi}_{\bar{\zeta}}$ is injective it suffices to establish the following two facts
\begin{itemize} \itemsep0.2em
\item[(1)] $d\Xi_{\zeta}:T_{\zeta} \XX_{k-1} \longrightarrow T_{[I_{n-k,n-k}^n]} \HH_{n-k,n-k}^n$ has a 1 dimensional kernel
\item[(2)] The exists $\omega \in T_{\bar{\zeta}} (\XX_{k-1}/\mathfrak{S}_2)$ for which $d\widebar{\Xi}_{\bar{\zeta}}(\omega)$ does not lie in the image of $d\Xi_{\zeta}$.
\end{itemize}

We begin with item (1). Let $\bg^1 = (k-1,k-2,\dots,0)$ and $\bg^2= (k,k+1,\dots,2k-2)$. Then $\zeta$ is the point $\mathbf{0}$ on $U^{\bg}_{k-1}$  (Proposition \ref{mainonee}). As in Remark \ref{coordinates} a set of coordinates on $U^{\bg}_{k-1}$ is $\mathcal{N} = \mathcal{N}_1 \cup \cdots\cup \mathcal{N}_5$ where
\begin{align*}
\mathcal{N}_1 = \{b_{i,j}\}_{0\leq i \leq k-1}^{k \leq j \leq n}, \quad 
\mathcal{N}_2 =  \{T^{\bg,(j)}_{i,k-1+j}\}^{0 \leq i \leq k-1-j}_{1 \leq j \leq k-1}, \quad	
\mathcal{N}_3 =   \{T_{k-i,j}^{\bg,(i)}\}_{k+i \leq j \leq n}^{1 \leq i \leq k-1}, \\
\mathcal{N}_4 = \{\lambda_1^{\bg},\dots,\lambda_{k-1}^{\bg}\} , \quad	
\mathcal{N}_5 =   \{T_{0,j}^{\bg,(k)}\}_{2k-1 \leq j \leq n}.
\end{align*}
For each $\eta \in \mathcal{N} $ we define a curve $D_{\eta}: \Spec \kk[t] \longrightarrow U_{k-1}^{\bg}$, passing through $\mathbf{0}$, by setting $\eta = t$ and all the other coordinates in $\mathcal{N}$ to $0$. 

Let $ \iota: \Spec \kk[t]/(t^2) \longrightarrow \Spec \kk[t]$ be a first order deformation of the origin. Since $\XX_{k-1}$ is non-singular the set $\{D_{\eta}\circ \iota \}_{\eta \in \mathcal{N}}$ is a basis for $T_{\mathbf{0}} U_{k-1}^{\bg} = T_{\zeta}\XX_{k-1}$. 
We need to study the dimension of $\{d\Xi_{\zeta}(D_{\eta}\circ \iota )\}_{\eta}$. Since $d\Xi_{\zeta}(D_{\eta}\circ \iota) =(\Xi \circ D_{\eta}) \circ \iota$ we begin with an explicit description of each  $\Xi \circ D_{\eta}$. The items below follow directly from the construction of the map (\ref{finalmap2}).

\begin{enumerate}
\item If $\eta = b_{i,j} \in \mathcal{N}_1$ then $\Xi \circ D_{\eta}(t)$ is
\begin{align*}
(x_0,\dots,x_{i-1},x_i+tx_j,x_{i+1},\dots,x_{k-1})^2 +(x_px_{2k-1-q})^{0 \leq p < q \leq k-1}_{p \ne i}  \\
				+ (x_i+tx_j)(x_k,\dots,x_{2k-2-i}).
\end{align*}
\item If $\eta = T^{\bg,(j)}_{i,k-1+j} \in \mathcal{N}_2$ then $\Xi \circ D_{\eta}(t)$ is 
$$
(x_0,\dots,x_{k-1})^2 +(x_px_{2k-1-q})^{0 \leq p < q \leq k-1}_{p \ne i} 
				+ (x_i-tx_{k-j})(x_k,\dots,x_{2k-2-i}).
$$
\item If $\eta = T^{\bg,(i)}_{k-i,j} \in \mathcal{N}_3$ then $\Xi \circ D_{\eta}(t)$ is
$$
(x_0,\dots,x_{k-1})^2  + (x_px_{2k-1-q})^{0 \leq p < q \leq k-1}_{q \ne k-i } 
				+ (x_0,\dots,x_{k-i-1})(x_{k-1+i}+tx_j).
$$
\item If $\eta = \lambda_i^{\bg}$ with $i >1$ then $\Xi \circ D_{\eta}(t)$ is 
$$
(x_0,\dots,x_{k-1})^2 + (x_px_{2k-1-q})_{(p,q) \ne (k-i,k-i+1)}^{0 \leq p < q \leq k-1} 
				+ (x_{k-i}x_{k+i-2}-tx_{k-i+1}x_{k+i-1}).
$$
\item If $\eta = \lambda_1^{\bg}$ then $\Xi \circ D_{\eta}(t)$ is 
$$
(x_0,\dots,x_{k-2})(x_0,\dots,x_{k-1}) +(x_{k-1}+tx_k)x_{k-1} + (x_px_{2k-1-q})_{0 \leq p < q \leq k-1}. 
$$
\item If $\eta = T_{0,j}^{\bg,(k)} \in \mathcal{N}_5$ then $\Xi \circ D_{\eta}(t)$ is
$$
(x_0,\dots,x_{k-1})^2 + (x_px_{2k-1-q})^{0 \leq p < q \leq k-1}_{(p,q) \ne (0,1)} 
				+ (x_0x_{2k-2}-tx_1x_j).
$$
\end{enumerate}

Let $I = I_{n-k,n-k}^n$ and under the inclusion $\mathcal{H}_{n-k,n-k}^n \subseteq \Hilb^{P_{n-k,n-k}^n(t)} \mathbf{P}^n$, we may identify $T_{[I]}\mathcal{H}_{n-k,n-k}^n$ with a subspace of $\Hom(I,S/I)_{0}$ (Lemma \ref{regularitythm}). We can explicitly describe this identification using  \cite[Proposition 2.3]{hdeform}. In particular, by re-indexing, we obtain
\begin{align*}
\text{span}_{\kk}\{d\Xi_{\zeta}(D_{\eta}\circ \iota)\}_{\eta \in \mathcal{N}_1 \cup \mathcal{N}_2 \cup \mathcal{N}_3} 
&= 	\text{span}_{\kk}\left(\left\{-x_j\frac{\partial}{\partial x_{i}}\right\}^{k \leq j \leq n}_{0 \leq i \leq k-1} \cup
	\left\{x_j\frac{\partial}{\partial x_{i}}\right\}^{i+1 \leq j \leq k-1}_{0 \leq i \leq k-2} 
	 \cup \left\{-x_j\frac{\partial}{\partial x_{i}}\right\}^{i+1 \leq j \leq n}_{k \leq i \leq 2k-2}\right) \\
&= 	 \text{span}_{\kk}\left\{x_j\frac{\partial}{\partial x_{i}}\right\}_{0 \leq i \leq 2k-2}^{i+1 \leq j \leq n}.
\end{align*}
These are the \emph{trivial deformations} i.e. the ones induced by a change of coordinates. For $i \in \{1,\dots,k-2\}$ let $\Delta_i$ be the derivation that maps $x_ix_{2k-2-i} \mapsto x_{i+1}x_{2k-1-i}$ and other generators to $0$. Let $\Delta_{k-1}$ denote the derivation that maps $x_{k-1}^2 \mapsto x_{k-1}x_k$ and the other generators to $0$.
For $i \in \{2k-1,\dots,n\}$ let $\Delta_i$ to the derivation that maps $x_0x_{2k-2} \mapsto x_{1}x_{i}$. Then we have 
$$
\text{span}_{\kk}\{d\Xi_{\zeta}(D_{\eta}\circ \iota)\}_{\eta \in \mathcal{N}_4 \cup \mathcal{N}_5} = 
	\text{span}_{\kk}(\{\Delta_i\}_{1 \leq i \leq k-1} \cup \{\Delta_i\}_{2k-1 \leq i \leq n}).
$$
Notice that the derivation $\Delta_{k-1}$ is a scalar multiple of $x_k\frac{\partial}{\partial x_{k-1}}$. Thus to prove (1) it suffices to show that the set $\{x_j\frac{\partial}{\partial x_{i}}\}_{0 \leq i \leq 2k-2}^{i+1 \leq j \leq n} \cup \{\Delta_i\}_{1 \leq i \leq k-2} \cup \{\Delta_i\}_{2k-1 \leq i \leq n}$ is linearly independent. 

Assume we had a linear combination 
\begin{equation} \label{derivationbro}
\sum_{\substack{0 \leq i \leq 2k-2 \\ i+1 \leq j \leq n}} \epsilon_{i,j}x_j\frac{\partial}{\partial x_i}
	+ \sum_{\substack{1\leq i \leq k-2 \\ 2k-1 \leq i \leq n}} \epsilon_{i}\Delta_i \equiv 0 \bmod I
\end{equation}
with some constants $\epsilon_{i,j},\epsilon_i \in \kk$. 
Assume $\epsilon_{p,q} \ne 0$ for some $p<q$. Since $x_px_{2k-2-p} \in I$ we may  evaluate (\ref{derivationbro}) at $x_px_{2k-2-p}$ to obtain
\begin{equation} \label{derivationbro1}
\sum_{p+1 \leq j \leq n} \epsilon_{p,j}x_jx_{2k-2-p} + \sum_{2k-1-p \leq j \leq n} \epsilon_{2k-2-p,j}x_jx_p + Q \equiv 0 \bmod I
\end{equation} 
where 
\begin{align*}
Q = \begin{cases} 
	\sum_{i=2k-1}^n \epsilon_ix_1x_i 	& \text{if } p = 0, \, 2k-2, \\
	 \epsilon_px_{p+1}x_{2k-1-p} 		& \text{if } 1 \leq p \leq k-2, \\
	 0							& \text{if } p = k-1 \\
	  \epsilon_{2k-2-p}x_{2k-1-p}x_{p+1} & \text{if } k \leq p \leq 2k-3.
	  \end{cases}
\end{align*} 
Observe that the monomial $x_qx_{2k-2-p}$ does not appear in the support of $Q$. Thus, in the left hand side of (\ref{derivationbro1}), the monomial $x_qx_{2k-2-p}$ appears with a coefficient of $\epsilon_{p,q}$ if $p \ne k-1$ and a coefficient of $2\epsilon_{p,q}$ if $p=k-1$.  In either case, the coefficient is non-zero. But this is a contradiction as $x_{q}x_{2k-2-p} \notin I$. Thus we have $\epsilon_{p,q} = 0$ for all $p,q$. Evaluating (\ref{derivationbro}) at $x_px_{2k-2-p}$ we see that $\epsilon_p = 0$ for every $p \in \{1,\dots,k-2\}$. Finally, evaluating (\ref{derivationbro}) at $x_0x_{2k-2}$ we  obtain $\sum_{i=2k-1}^n\epsilon_i x_1x_i \equiv 0 \bmod I$. Since $x_1x_i \notin I$ for all $i \geq 2k-1$, we must have that $\epsilon_i = 0$ for all $i$. This completes the proof of item (1).

Let $\Delta \in \Hom(I,S/I)_0$ denote the derivation that maps $x_{k-1}x_k \mapsto x_k^2$ and all the other generators to $0$. By evaluating at $x_{k-1}x_k$ it is easy to see that $\Delta$ does not lie in the span of $\{x_j\frac{\partial}{\partial x_{i}}\}_{0 \leq i \leq 2k-2}^{i+1 \leq j \leq n} \cup \{\Delta_i\}_{1 \leq i \leq k-2} \cup \{\Delta_i\}_{2k-1 \leq i \leq n}$. Consider the curve $C : \Spec \kk[t] \to \HH_{n-k,n-k}^n$ given by 
$$
t \mapsto (x_0,\dots,x_{k-2})(x_0,\dots,x_{k-1}) + (x_{k-1}^2 - tx_k^2) + (x_px_{2k-1-q})_{0 \leq p < q \leq k-1}
$$
This is well defined because for any given $s \in \kk$,  $C(s)$ is the point in  $U_{k-1}^{\bg}$ with $\lambda_1^{\bg} = -2\sqrt{s}$, $b_{k-1,k} = \sqrt{s}$ and all other coordinates equal $0$. It is also clear that $C \circ \iota$ corresponds to the derivation $\Delta$. Thus to prove item (2) it suffices to find a curve $C':  \Spec \kk[t] \to \XX_{k-1}/\mathfrak{S}_2$ passing through $\bar{\zeta}$ for which $d_{\bar{\zeta}}\widebar{\Xi} (C' \circ \iota) = C \circ \iota$.

Let  $Z$ denote the image of $C$ and let $Z'$ denote the pullback $\widebar{\Xi}^{-1}(Z) \subseteq \XX_{k-1}/\mathfrak{S}_2$. I claim that $\widebar{\Xi}|_{Z'}:Z' \to Z$ is an isomorphism.  Since $Z$ is non-singular, $Z'$ is Cohen-Macaulay and $\widebar{\Xi}$ is bijective, the morphism $\widebar{\Xi}|_{Z'}$ is flat. It is clear that a finite flat degree $1$ morphism is an isomorphism. Thus $C' = \widebar{\Xi}|_{Z'}^{-1} \circ C: \Spec \kk[t] \to \XX_{k-1}/\mathfrak{S}_2$  is the desired curve. 
\end{proof}

We are now ready to prove the main Theorem.

\begin{customthm}{A} \label{main} Let $n \geq 2k-1$. The component $\mathcal{H}_{n-k,n-k}^{n}$ is smooth and isomorphic to 
$$
\XX_{k-1}/\mathfrak{S}_2 = \Bl_{\widebar{\Gamma}_{k-1}}\cdots\Bl_{\widebar{\Gamma}_{1}}\Sym^{2}\Gr(n-k,n).
$$ 
\end{customthm}
\begin{proof} Proposition \ref{mainonee} and \ref{INJECTIVE} together show that $\widebar{\Xi}$ is bijective and $\XX_{k-1}/\mathfrak{S}_2$ is non-singular. Since $\widebar{\Xi}$ is $\GL(n+1)$-equivariant, $\bar{\zeta}$ (Definition \ref{zetazeta}) is the unique Borel fixed point on $\XX_{k-1}/\mathfrak{S}_2$. By Borel's fixed point theorem, the closure of the Borel orbit of any point in $\XX_{k-1}/\mathfrak{S}_2$ contains $\bar{\zeta}$. Thus to show that $\widebar{\Xi}$ is an isomorphism, it suffices to show that it is an isomorphism in a neighbourhood of $\bar{\zeta}$. By the proof of \cite[Theorem 14.9]{harris}, this is equivalent to showing that $d\widebar{\Xi}_{\bar{\zeta}}: T_{\bar{\zeta}} (\XX_{k-1}/\mathfrak{S}_2) \longrightarrow T_{[I_{n-k,n-k}^n]} \HH_{n-k,n-k}^n$ is injective. This is precisely the content of Lemma \ref{dacurves}.
\end{proof}

When the pair of planes do not span $\PP^n$, we obtain the following fibration

\begin{cor} \label{mainthree} Let $n < 2k-1$. The morphism $\rho:\mathcal{H}_{n-k,n-k}^n \longrightarrow \mathbf{Gr}(2n-2k+1,n)$ that sends a scheme to its linear span is smooth; the fiber over a point $\Lambda$ is $\mathcal{H}_{n-k,n-k}(\Lambda)$. 
\end{cor}
\begin{proof} Recall that the linear span of a subscheme $Z \subseteq \mathbf{P}^n$ is the linear space $V(H^0(\mathbf{P}^n,I_Z(1)))  \subseteq \mathbf{P}^n$.
Let $\mathcal{Y} \longrightarrow \mathbf{A}^1$ be a flat family such that for $t\ne 0$, $\mathcal{Y}_t$ is a disjoint pair of $(n-k)$-planes. It is clear that for any $t \ne 0$, the linear span of $\mathcal{Y}_t$ is a $(2n-2k+1)$-plane. By upper semicontunity, the limit $\mathcal{Y}_0$ also lies in a $(2n-2k+1)$-plane, which we denote by $\Lambda$. Thus $\mathcal{Y}_0$ defines a point in $\mathcal{H}_{n-k,n-k}^n(\Lambda)$ and by Corollary \ref{2kgens}, we see that the linear span of $\mathcal{Y}_0$ is all of $\Lambda$. It follows that the linear span of any subscheme parameterized by $\mathcal{H}_{n-k,n-k}(\PP^n)$ is of dimension $2n-2k+1$.

\smallskip
For each ordered basis $\EE=\{e_0,\dots,e_n\}$ of $S_1$ we obtain an open neighbourhood $U_{\EE} = \Spec \kk[f_{i,j}]_{0 \leq i \leq 2k-2-n}^{2k-1-n\leq j \leq n}$ of $\Lambda_{\EE} = V(e_0,\dots,e_{2k-2-n})$ in $\Gr(2n-2k+1,n)$. The $\kk$-point $\mathbf{f} = (f_{i,j})_{i,j}$ is identified with
$$
V(e_0 + \sum_{j=2k-1-n}^n f_{0,j}e_j,\dots,e_{2k-2-n} + \sum_{j=2k-1-n}^n f_{2k-2-n,j}e_j).
$$
Let $\EE = \{e_i\}_i, \EE' = \{e_i'\}_i$ be ordered bases of $S_1$. The isomorphism $\Lambda_{\EE} \to \Lambda_{\EE'}$ given by mapping $e_i \mapsto e'_i$ for all $i$ induces an an isomorphism $\psi_{\EE,\EE'}: \HH_{n-k,n-k}(\Lambda_{\EE}) \longrightarrow \HH_{n-k,n-k}(\Lambda_{\EE'})$. Define the following
\begin{itemize}
\item $\XX_{\EE} = \HH_{n-k,n-k}(\Lambda_{\EE}) \times U_{\EE}$, 
\item $\XX_{\EE,\EE'} = \HH_{n-k,n-k}(\Lambda_{\EE}) \times (U_{\EE} \cap U_{\EE'}) \subseteq \XX_{\EE}$,
\item $\varphi_{\EE,\EE'} = \psi_{\EE,\EE'} \times \text{id}  :\XX_{\EE,\EE'} \longrightarrow \XX_{\EE',\EE}. $ 
\end{itemize}
It is clear that 
$\varphi_{\EE,\EE'}^{-1}= \varphi_{\EE',\EE}$,
$ \varphi_{\EE',\EE''} \circ \varphi_{\EE,\EE'} = \varphi_{\EE,\EE''}$  on $\XX_{\EE,\EE'} \cap \XX_{\EE,\EE''}$ and
$\varphi_{\EE,\EE'}(\XX_{\EE,\EE'} \cap \XX_{\EE,\EE''}) = \XX_{\EE',\EE} \cap \XX_{\EE',\EE''}$.
Thus the set of schemes $\{X_{\EE}\}_{\EE}$ glue to a smooth scheme $\XX$ (Theorem \ref{main}). 

For each $\EE$ we obtain a natural morphism $g_{\EE}: U_{\EE} \longrightarrow \text{GL}(n+1)$ such that for any $\mathbf{f}$, $g_{\EE}(\mathbf{f})$ is the map that sends $e_i \mapsto e_i + \sum_{j=2k-1-n}^n f_{i,j}e_j$ if $i \leq 2k-2-n$ and fixes the other coordinates. Thus we may define a map
$$
\HH_{n-k,n-k}(\Lambda_{\EE}) \times U_{\EE} \longrightarrow \HH_{n-k,n-k}(\PP^n), 
			\quad (X,\mathbf{f}) \mapsto g_{\EE}(\mathbf{f})(X).
$$
These maps glue to a morphism $ \Pi: \XX \longrightarrow \HH_{n-k,n-k}^n$. By the first paragraph, $\Pi$ is a bijective morphism. It is also clear that the differential to $\Pi$ is injective at all points. As noted in Theorem \ref{main}, this implies that $\Pi$ is an isomorphism. By construction, there is a smooth fibration $\rho:\XX \longrightarrow \Gr(2n-2k+1,n)$ of the desired form.
\end{proof}

\begin{customthm}{C} \label{uniqueborel2} $\mathcal{H}_{n-k,n-k}^n$ has a unique Borel fixed point.
\end{customthm}
\begin{proof} By Proposition \ref{uniqueborel} we my assume $n < 2k-1$. If $X$ is Borel fixed then its linear span $V((I_{X})_1)$ is also Borel fixed. Thus $X$ lies in the fiber $\rho^{-1}(V(x_0,\dots,x_{2k-2-n})) \simeq \mathcal{H}_{n-k,n-k}^{2n-2k+1}$. Moreover, the Borel action on $\mathcal{H}_{n-k,n-k}^n$ restricts to the Borel action on this fiber. By Proposition \ref{uniqueborel} this fiber has a unique Borel fixed point; thus $X$ is unique.
\end{proof}

We now turn our attention to the subschemes parameterized by $\mathcal{H}_{n-k,n-k}^n$. Since we are going to describe these subschemes up to projective equivalence, we may assume $n \geq 2k-1$ (Corollary \ref{mainthree}). We begin with two Lemmas that will aid in the proof of Theorem \ref{primdecomp}.

\begin{lemma} \label{CM} Let $J = (x_0,\dots,x_{k-1})^2 + (x_{p}x_{n-k_q}-\mu_{p,q}x_{q}x_{n-k_p})_{0 \leq p < q \leq k-1}$ with $\mu_{i} \in \kk$ and $\mu_{p,q} = \mu_{k-q+1}\cdots \mu_{k-p}$ for any $0 \leq p < q \leq k$. If all the $\mu_{i}$ are non-zero then the subscheme defined by $J$ is Cohen-Macaulay; in particular, it has no embedded components. Moreover, the subscheme defined by $J$ is double structure on $V(x_0,\dots,x_{k-1})$.

\end{lemma}
\begin{proof} Applying the change of coordinates that maps $x_p \mapsto \mu_{p,k}x_p$ for all $p \leq k-1$ and fixing the other coordinates, we may assume $\mu_{p,q}=1$ for all $p,q$. If $n > 2k-1$, the variables $x_{k},\dots,x_{n-k}$ form a regular sequence as they do not appear in the support of the generators of $J$. Thus we may quotient by the ideal $(x_k,\dots,x_{n-k})$  to reduce to the case $n=2k-1$; in this case $n-k_p =k+p$. Since $\Proj(S/J)$ is supported on $V(x_0,\dots,x_{k-1})$, it suffices to verify the Cohen-Macaulayness on the open sets $D(x_k),\dots,D(x_{2k-1})$. 

On the open set $W= D(x_{k})$ we may set $x_k = 1$. Then for all $j\ne 0$ we have $x_{j} - x_0x_{k+j} =-(x_0x_{k+j} - x_{j}x_k) \in J\vert_W$ and this implies $J\vert_W = (x_0^2, x_1-x_0x_{k+1},\dots,x_{k-1}-x_0x_{2k-1})$. Since $x_k,\dots,x_{2k-1}$ forms a regular sequence on $(S/J)\vert_W$, $\Proj(S/J)|_{W}$ is a  Cohen-Macaulay subscheme of dimension $k-1$. The argument for the other open sets is the same.

Since the Hilbert polynomial of $\Proj(S/J)$ is $P_{n-k,n-k}^n(t)$, its degree is $2$; thus it is a double structure on the linear space $V(x_0,\dots,x_{k-1})$
\end{proof}

\begin{remark} \label{CM'} More generally, $(x_{\epsilon_1},\dots,x_{\epsilon_2})^2+(x_{p}x_{n-k_q}-\mu_{p,q}x_{q}x_{n-k_p})_{\epsilon_1 \leq p < q \leq \epsilon_2}$ is Cohen-Macaulay for any $0 \leq \epsilon_1 \leq \epsilon_2 \leq k-1$, assuming $\mu_i \ne 0$ for all $i$. 
\end{remark}

\begin{lemma} \label{monodecomp} Let $0 \leq \epsilon_1 \leq \epsilon_2 \leq k-1$ and let $J(\epsilon_1,\epsilon_2) = (x_{\epsilon_1},\dots,x_{\epsilon_2})^2 +(x_px_{n-k_q})_{\epsilon_1 \leq p < q \leq \epsilon_2}$. Then we have a primary decomposition
$$
J(\epsilon_1,\epsilon_2)= 
\bigcap_{j = \epsilon_1}^{\epsilon_2} (x_{\epsilon_1},\dots,x_{j-1},x_j^2,x_{j+1},\dots,x_{\epsilon_2},x_{n-k_{j+1}},\dots,x_{n-k_{\epsilon_2}}).
$$
\end{lemma}
\begin{proof} For the first statement we proceed by induction on $\epsilon_2$. The base case $\epsilon_2 = \epsilon_1$ is vacuous and by induction we may assume
\begin{align*}
J(\epsilon_1,\epsilon_2+1)=\left[(x_{\epsilon_1},\dots,x_{\epsilon_2})^2 +(x_px_{n-k_q})_{\epsilon_1 \leq p < q \leq \epsilon_2} + (x_{\epsilon_2+1},x_{n-k_{\epsilon_2+1}})\right] \cap (x_{\epsilon_1},\dots,x_{\epsilon_2},x_{\epsilon_2+1}^2). 
\end{align*}
The conclusion now follows from the fact that if $I_1=(m_1,\dots,m_{i_1}),I_2 = (m_1,\dots,m_{i_2})$ are monomial ideals then $I_1 \cap I_2 = (\mathrm{lcm}(m_im_j): 1 \leq i \leq i_1, 1 \leq j \leq i_2)$.
\end{proof}

\begin{customthm}{D} \label{primdecomp} Let $n\geq 2k-1$. Let $Z$ be a subscheme parameterized by $\mathcal{H}_{n-k,n-k}^n$.  Then $Z$ is a pair of planes meeting transversely, or there exists a sequence of integers $1 \leq i_1 < \cdots < i_r \leq k$ and a flag of linear spaces $\Lambda^1 \subseteq \Lambda^{2} \subseteq \cdots \subseteq \Lambda^{r} \subseteq \mathbf{P}^n$ with $\mathrm{codim}_{\PP^n}(\Lambda^{\ell})= (k+i_{\ell}-1)$ for each $\ell$, such that
\begin{enumerate}
\item If $i_1 >1$ then $Z$ is a union of two planes meeting along $\Lambda^1$ with embedded pure double structures on $\Lambda^{\ell}$ for each $1 \leq \ell \leq r$.
\item If $i_1 =1$ then $Z$ is a pure double structure on $\Lambda^1$ with embedded pure double structures on $\Lambda^{\ell}$ for each $2 \leq \ell \leq r$. 
\end{enumerate}
\end{customthm}
\begin{proof} It suffices to compute a primary decomposition of the ideal
$$
J = (x_{p}+\mu_{p,k}x_{n-k_p})_{0\leq p\leq k-1}(x_0,\dots,x_{k-1}) + (x_px_{n-k_q} -\mu_{p,q} x_q x_{n-k_p})_{0\leq p <q \leq k-1}
$$
in (\ref{projideal}). Let $\mathfrak{P}_0 = (x_p+\mu_{p,k}x_{n-k_p})_{0 \leq p \leq k-1}$, $\mathfrak{P}_1 = (x_0,\dots,x_{k-1})$ and $\delta_{p,q} =x_px_{n-k_q} -\mu_{p,q} x_q x_{n-k_p}$ for each $0 \leq p < q \leq k-1$. Lemma \ref{anothermingens} (ii) implies that all the $\mu_i$ are non-zero if and only if $J$ is the ideal of a pair of $(n-k)$-planes meeting transversely. So we may assume some of the $\mu_i$ are zero. Let $i_1 < \cdots < i_r$ be all the indices $i$ for which $\mu_{i} = 0$. Set $i_0 =0 $ and $i_{r+1} = k+1$.  Lemma \ref{anothermingens} (iv) implies $\sqrt{J} = \mathfrak{P}_0 \cap \mathfrak{P}_1$ and  $J = \mathfrak{P}_0\mathfrak{P}_1 +  (\delta_{p,q})_{0\leq p <q \leq k-i_1}$. For each $2 \leq \ell \leq r+1$ define
\begin{align*}
\mathfrak{P}_{\ell} = (x_0,\dots,x_{k-i_{\ell}}) + (x_{k-i_\ell+1},\dots,x_{k-i_{\ell-1}})^2
 			+(\delta_{p,q})_{k-i_{\ell}+1 \leq p < q \leq k -i_{\ell-1}} +\\ 
			(x_{k-i_{\ell-1}+1},\dots,x_{k-1},x_{n-i_{\ell-1}+2},\dots,x_n).
\end{align*}

I claim that $J = \mathfrak{P}_0 \cap \mathfrak{P}_1 \cap \cdots \cap \mathfrak{P}_{r+1}$ (note that if $\mu_1 = 0$ then $\mathfrak{P}_0 =\mathfrak{P}_1$). We begin with the inclusion, $J \subseteq \mathfrak{P}_0 \cap \cdots \cap \mathfrak{P}_{r+1}$. It is enough to show that $\mathfrak{P}_0\mathfrak{P}_1$ and $\delta_{p,q}$ lie in $\mathfrak{P}_0 \cap \cdots \cap \mathfrak{P}_{r+1}$ for $0 \leq p < q \leq k-i_1$. Observe that
\begin{align*}
\mathfrak{P}_{0}\mathfrak{P}_{1} & = 
			((x_0,\dots,x_{k-i_1}) +(x_p + \mu_{p,k}x_{n-k_p})_{k-i_1 +1 \leq p \leq k-1})(x_0,\dots,x_{k-1})
\end{align*}
Clearly, $(x_0,\dots,x_{k-i_1})(x_0,\dots,x_{k-1}) \subseteq \mathfrak{P}_j$ for all $j$. We also have, $x_p, x_{n-k_p} \in \mathfrak{P}_j$ for all $k-i_1 +1 \leq p \leq k-1$ and all $j$. Thus $\mathfrak{P}_0\mathfrak{P}_1 \subseteq \mathfrak{P}_0 \cap \cdots \cap \mathfrak{P}_{r+1}$. It is clear that $\delta_{p,q} \in \mathfrak{P}_0 \cap \cdots \cap \mathfrak{P}_{r+1}$ if there is some $\ell$ such that $k-i_{\ell}+1 \leq p < q \leq k -i_{\ell-1}$. If this was not the case, then there is some $\ell$ such that $p \leq k-i_{\ell} < q$. This implies $\delta_{p,q} = x_{p}x_{n-k_q}$ and this lies in $(x_0,\dots,x_{k-i_j})$ if $j \leq \ell$ or in $(x_{n-i_{j-1}+2},\dots,x_n)$ if $j > \ell$; in either case, $\delta_{p,q} \in \mathfrak{P}_j$. Thus $\delta_{p,q} \in \mathfrak{P}_0 \cap \cdots \cap \mathfrak{P}_{r+1}$ and we have the desired containment. 

To get the other containment it suffices to show that $\mathfrak{P}_0 \cap \cdots \cap \mathfrak{P}_{r+1}$ has the same Hilbert function as $J$. We have
\begin{equation} \label{innnn}
\inn_{>} J \subseteq \inn_{>}\left(\mathfrak{P}_0 \cap \cdots \cap \mathfrak{P}_{r+1}\right) \subseteq
\inn_{>} (\mathfrak{P}_0 \cap \mathfrak{P}_1)  \cap \inn_{>}\mathfrak{P}_2 \cap \cdots \cap \inn_{>} \mathfrak{P}_{r+1}.
\end{equation}
Our goal is to show all these containments are equalities. Using Equation (\ref{swapeq}) we have 
\begin{align*}
\mathfrak{P}_{0}\cap \mathfrak{P}_{1} 
	&=((x_0,\dots,x_{k-i_1}) + (x_p+\mu_{p,k}x_{n-k_p})_{k-i_1+1\leq p\leq k-1})\cap (x_0,\dots,x_{k-1}) \\
	& = (x_0,\dots,x_{k-i_1}) + (x_p+\mu_{p,k}x_{n-k_p})_{k-i_1 +1 \leq p \leq k-1}\cap(x_{k-i_1+1},\dots,x_{k-1}) \\
	& = (x_0,\dots,x_{k-i_1}) +(x_p+\mu_{p,k}x_{n-k_p})_{k-i_1 +1 \leq p \leq k-1}(x_{k-i_1+1},\dots,x_{k-1})\\
	&= (x_0,\dots,x_{k-i_1}) + ((x_p+\mu_{p,k}x_{n-k_p})x_q)_{k-i_1 +1 \leq p \leq q \leq k-1} + (\delta_{p,q})_{k-i_1+1 \leq p <q \leq k-1}.
\end{align*}
Then the proof of Lemma \ref{GROBNER} immediately implies 
$$
\inn_{>} (\mathfrak{P}_{0}\cap \mathfrak{P}_{1}) = (x_0,\dots,x_{k-i_1}) +(x_{k-i_1+1},\dots,x_{k-1})^2 + (x_px_{n-k_q})_{k-i_1+1 \leq p <q \leq k-1}.
$$
Similarly for $\ell \geq 2$
\begin{align*}
\inn_{>} \mathfrak{P}_{\ell} = & \, \, (x_0,\dots,x_{k-i_{\ell}}) + (x_{k-i_\ell+1},\dots,x_{k-i_{\ell-1}})^2
 			+(x_px_{n-k_q})_{k-i_{\ell}+1 \leq p < q \leq k -i_{\ell-1}} +\\ 
			& \quad \quad (x_{k-i_{\ell-1}+1},\dots,x_{k-1},x_{n-i_{\ell-1}+2},\dots,x_n).
\end{align*}
Using Lemma \ref{monodecomp} we see that $\inn_{>} (\mathfrak{P}_0 \cap \mathfrak{P}_1) \cap \inn_{>} \mathfrak{P}_2 \cap \cdots \cap \inn_{>} \mathfrak{P}_{r+1}$ equals
$$
\bigcap_{\ell=1}^{r+1} \bigcap_{j=k-i_{\ell}+1}^{k-i_{\ell-1}}(x_0,\dots,x_{j-1},x_j^2,x_{j+1},\dots,x_{k-1},x_{n-k_{j+1}},\dots,x_n).\footnote{If $j =k$ the ideal $(x_0,\dots,x_{j-1},x_j^2,x_{j+1},\dots,x_{k-1},x_{n-k_{j+1}},\dots,x_n)$ is equal to $(x_0,\dots,x_{k-1})$.}
$$
Applying Lemma \ref{monodecomp} once again we see that this intersection is just $J(0,k-1) \cap (x_0,\dots,x_{k-1})$. But this ideal is precisely $(x_0,\dots,x_{k-1})^2 + (x_px_{n-k_q})_{0 <p <q \leq k-1} = \inn_{>} J$. Thus all the containments in (\ref{innnn}) are equalities and this shows that $J$ has the same Hilbert function as $\mathfrak{P}_{0}\cap \cdots \cap \mathfrak{P}_{r}$.

We are left with showing $\mathfrak{P}_{\ell}$ is a primary component for all $\ell \geq 2$. Going modulo the linear forms it suffices to show that
$(x_{k-i_\ell+1},\dots,x_{k-i_{\ell-1}})^2+(\beta_{p,q})_{k-i_{\ell}+1 \leq p < q \leq k -i_{\ell-1}}$ 
is a primary component. This is the content of Lemma \ref{CM} and Remark \ref{CM'}.
\end{proof}

\begin{customcor}{E} \label{2kpointss} Up to projective equivalence, there are exactly $2^k$ schemes parameterized by $\mathcal{H}_{n-k,n-k}^n$.
\end{customcor}
\begin{proof} By Corollary \ref{mainthree} we may assume $n \geq 2k-1$. It suffices to consider ideals $J$ of the form (\ref{projideal}). Let $\varphi$ denote the projective transformation that maps $x_p \mapsto \mu_{p,k}x_p$ if $\mu_{p,k} \ne 0$ and $0 \leq p \leq k-1$ and fixes the other coordinates. For a fixed $p$, note that if $\mu_{p,k} \ne 0$ then $\mu_{q,k} \ne 0$ and $\mu_{p,q} \ne 0$ for all $p<q$. Thus after applying $\varphi$ we may assume that the non-zero $\mu_i$ are equal to $1$. In particular, for each subset $W \subseteq \{1,\dots,k\}$ we obtain an ideal parameterized by $\mathcal{H}_{n-k,n-k}^n$ by setting $\mu_i = 0$ if $i \in W$ and $1$ otherwise; this gives at most $2^k$ distinct ideals. On the other hand, since projective transformations preserve the dimensions of the embedded structures, each of the $2^k$ ideals are projectively inequivalent.
\end{proof}


\begin{example} \label{specializations} We can now determine when there is a specialization $Z\rightsquigarrow Z'$ in $\HH_{n-k,n-k}^n$. For any subscheme $Z \in \HH_{n-k,n-k}^n$ let $W_Z = \{\epsilon_1,\dots,\epsilon_r\}$ be the set of dimensions of the embedded components of $Z$; if $Z$ is generically non-reduced include $n-k$ in that set. Then there is a specialization $Z\rightsquigarrow Z'$  if and only if $W_Z \subseteq W_{Z'}$

Here is a diagram of specializations for $\HH_{2,2}^5$. The non-reduced structures on points, lines and planes are represented by shadings. 

\begin{tikzpicture}[scale=1.4] \label{myprecious}

\draw (-2.8,-1.4) -- (-2.8,-0.4);
\draw (-2.8,-0.4) -- (-1.8,-0.4) node[opaque, above left =2pt] {(i)};
\draw (-1.8,-0.4) -- (-1.8,-1.4);
\draw[dashed] (-2.8,-1.4) -- (-1.8,-1.4);

\draw[rotate around={45:(-3,-3)}] (-2,-3)   -- ++(1,0) -- ++(0,1)   -- ++(-1,0) -- cycle ;
\draw[->][color=blue][thick] (-1.3,-1.3) -- (-0.7,-1.3);
\draw[->][color=blue][thick] (-1.3,-1.3) -- (-0.5,0.5);
\draw[->][color=blue][thick] (-1.3,-1.3) -- (-0.3,-3.0);


\draw (0,1) -- (1,1);
\draw (0,0) -- (0,1) node[below left=2pt] {(ii)};
\draw (0,0) -- (1,0);
\draw[dashed] (1,0) -- (1,1);

\draw[shift={+(-0.2,-0.9)}][rotate=45] (1.5,-0.5) -- (2.5,-0.5);
\draw[shift={+(-0.2,-0.9)}][rotate=45]   (1.5,-0.5)  -- (1.5,0.5);
\draw[shift={+(-0.2,-0.9)}][rotate=45] (1.5,0.5) -- (2.5,0.5);
\draw[shift={+(-0.2,-0.9)}][rotate=45]  (2.5,-0.5) -- (2.5,0.5);
\node[fill,circle,scale=.5] at (0.51,0.51) {};

\draw[->][color=blue][thick] (2.2,0.7) -- (3.3,0.7) ;
\draw[->][color=blue][thick] (2.2,0.7)  -- (3.3,-1) ;


\draw (0,-1) -- (0.5,-1);
\draw[dashed] (0.5,-1) -- (1,-1);
\draw (0,-2) -- (0,-1) node[below left=2pt] {(iii)};
\draw (0,-2) -- (1,-2);
\draw (1,-2) -- (1,-1.5);
\draw[dashed] (1,-1.5) -- (1,-1);

\draw[line width=3pt] (0.5,-2) -- (0.5,-1);
\draw[rotate around={45:(0.5,-2)}] (0.5,-2) -- (1.5,-2);
\draw[rotate around={45:(0.5,-1)}] (0.5,-1) -- (1.5,-1);
\draw (1.207,-0.293) -- (1.207,-1.293);

\draw[->][color=blue][thick] (1.5,-1.5) -- (3.3,0) ;
\draw[->][color=blue][thick] (1.5,-1.5)  -- (3.3,-3.0) ;

\fill (0,-4)  -- ++(1,0)  -- ++(0,1)  -- ++(-1,0) node[below left=2pt] {(v)}  -- cycle ;

\draw[->][color=blue][thick] (1.5,-3.5) -- (3.3,-1.8) ;
\draw[->][color=blue][thick] (1.5,-3.5)  -- (3.3,-3.5) ;


\draw (4,1) -- (4.5,1);
\draw[dashed] (4.5,1) -- (5,1);
\draw (4,0) -- (4,1) node[below left=2pt] {(iv)} ;
\draw (4,0) -- (5,0);
\draw (5,0) -- (5,0.5);
\draw[dashed] (5,0.5) -- (5,1);
\node[fill,circle,scale=.5] at (4.51,0.51) {};

\draw[line width=3pt][nearly opaque] (4.5,0) -- (4.5,1);
\draw[rotate around={45:(4.5,0)}] (4.5,0) -- (5.5,0);
\draw[rotate around={45:(4.5,1)}] (4.5,1) -- (5.5,1);
\draw (5.207,1.707) -- (5.207,0.707);

\draw[->][color=blue][thick] (5.5,0.5)  -- (7,-0.5) ;

\fill[nearly opaque] (4,-2)  -- ++(1,0)  -- ++(0,1)  -- ++(-1,0) node[opaque, below left=2pt] {(vi)}  -- cycle ;
\node[fill,circle,scale=.5] at (4.51,-1.49) {};

\draw[->][color=blue][thick] (5.5,-1.5)  -- (6.8,-1.5) ;

\fill[nearly transparent] (4,-4)  -- ++(1,0)  -- ++(0,1)  -- ++(-1,0) node[opaque, below left=2pt] {(vii)}  -- cycle ;
\draw[line width=3pt] (4.5,-4) -- (4.5,-3);

\draw[->][color=blue][thick] (5.5,-3.5) -- (7,-2.3) ;

\fill[nearly transparent] (7,-2)  -- ++(1,0)  -- ++(0,1) node[opaque, above left=2pt] {(viii)}  -- ++(-1,0) -- cycle ;
\draw[line width=3pt][semitransparent] (7.5,-2) -- (7.5,-1);
\node[fill,circle,scale=.5] at (7.51,-1.49) {} ;

\end{tikzpicture}
\end{example}

\begin{remark} In \cite{vain}, Vainsencher uses the map $\Xi: \Bl_{\Gamma_2}\Bl_{\Gamma_1}\mathbf{Gr}(2,5)^2 \to \HH_{2,2}^5$ to compute the degree of a family of rational cubic fourfolds in $\mathbf{P}^5$. However, he does not prove the smoothness of $\HH_{2,2}^5$.
\end{remark}

\section{Structure of $\mathcal{H}_{n-c,n-d}^n$} \label{three}
In this short section we explain how the proofs of the previous section carry over, almost identically, to the case when the pair of planes are of different dimension. We begin by explaining the special case of $c=1$ that we have omitted.

\begin{remark} If $c =1$ then $\Hilb^{P_{n-1,n-d}^n(t)} \mathbf{P}^n$ parameterizes ideals of codimension $1$. Using the decomposition in \cite[Proposition 2.4]{ritvik} we obtain
$$
\Hilb^{P_{n-1,n-d}^n(t)} \mathbf{P}^n = \Hilb^{\binom{n-1+t}{t}} \mathbf{P}^n \times \Hilb^{\binom{n-d+t}{t}} \mathbf{P}^n \simeq \Gr(n-1,n) \times \Gr(n-d,n).
$$
Thus $\HH_{n-1,n-d}^n$ is smooth and isomorphic to the full Hilbert scheme.  Alternatively, we can deduce this from the proof of Lemma \ref{twotoone} and a computation of the tangent space to the unique Borel fixed ideal on $\Hilb^{P_{n-1,n-d}^n(t)} \mathbf{P}^n$.
\end{remark}

Let $d>c \geq 2$ and assume $n \geq c + d -1$. Let $\mathcal{X}_{c-1} = \Bl_{\Gamma_{c-1}} \cdots \Bl_{\Gamma_1} (\mathbf{Gr}(n-c,n) \times \mathbf{Gr}(n-d,n))$ and let $\pi_{c-1}:\mathcal{X}_{c-1} \longrightarrow \mathbf{Gr}(n-c,n) \times \mathbf{Gr}(n-d,n) $ be the blow up. 

We have shown in Lemma \ref{twotoone} that the rational map $\Xi: \mathcal{X}_{c-1} \dashrightarrow \mathcal{H}_{n-c,n-d}^n$ is defined and one-to-one on the open set $\Gr(n-c,n)\times \Gr(n-d,n) \setminus \Gamma_1 \cup \cdots \cup \Gamma_{c-1}$. To extend $\Xi$ to $\XX_{c-1}$ we proceed as in Section \ref{two}. We first extend $\Xi$ to $\pi_{c-1}^{-1}(U_0)$ where $U_0= \Spec\kk[a_{i,j},b_{i,j}]_{i,j}$  is an open subset of $\Gr(n-c,n)\times \Gr(n-d,n)$ such that its $\kk$-points correspond to
$$
(\Lambda(\aaa),\Lambda(\bb)) = (V(x_{0}+\sum_{j=c}^{n}a_{0,j}x_{j},\dots,x_{c-1}+\sum_{j=c}^{n}a_{c-1,j}x_{j}),V(x_{0}+\sum_{j=d}^{n}b_{0,j}x_{j},\dots,x_{d-1}+\sum_{j=d}^{n}b_{d-1,j}x_{j})).
$$

We will now perform a few substitutions and obtain a different minimal set of generators for $I_{\Lambda(\aaa)}$ and $I_{\Lambda(\bb)}$. From these new presentations of $I_{\Lambda(\aaa)}$ and $I_{\Lambda(\bb)}$, it will be apparent how one has to mimic the arguments of Section \ref{two} to extend $\Xi$ to $\pi_{c-1}^{-1}(U_0)$, and thus all of $\XX_{c-1}$. For every $0 \leq i \leq c-1$, $0 \leq j \leq d-1$ and $d \leq p \leq n$ let
$$ 
y'_j = x_j+\sum_{\ell =d}^n b_{j,\ell}x_\ell, \quad y_i = y_i'+ \sum_{\ell = c}^{d-1}a_{i,\ell}y'_{\ell}, \quad \widetilde{b}_{i,p} = b_{i,p} + \sum_{\ell=c}^{d-1}a_{i,\ell}b_{\ell,p}.
$$
For any $0 \leq i \leq c-1$ we obtain
\begin{align*}
x_i + \sum_{j=c}^na_{i,j}x_j & = y_i' + \sum_{j=d}^n(a_{i,j}-b_{i,j})x_j + \sum_{j=c}^{d-1}a_{i,j}x_j \\
			& = y_i' + \sum_{j=d}^n(a_{i,j}-b_{i,j})x_j + \sum_{j=c}^{d-1}a_{i,j}\left(y_j' -\sum_{\ell = d}^nb_{j,\ell}x_{\ell}\right)\\
			& = y_i' + \sum_{j=c}^{d-1}a_{i,j}y_j' + \sum_{j=d}^n\left(a_{i,j}-b_{i,j} - \sum_{\ell=c}^{d-1}a_{i,\ell}b_{\ell,j}\right)x_j\\
			& = y_i + \sum_{j=d}^n(a_{i,j}-\widetilde{b}_{i,j})x_j.
\end{align*}
Thus we have 
\begin{align} \label{aaa}
I_{\Lambda(\aaa)} = (y_0+\sum_{j=d}^n(a_{0,j}-\widetilde{b}_{0,j})x_j,\dots,
y_{c-1}+\sum_{j=d}^n(a_{c-1,j}-\widetilde{b}_{c-1,j})x_j)
\end{align}
and 
\begin{align} \label{bb}
I_{\Lambda(\bb)} = (y_0',\dots,y_{d-1}') 
		&= (y_0,\dots,y_{c-1},y_c',\dots,y_d') \nonumber \\
		&= (y_0,\dots,y_{c-1},x_c+\sum_{j =d}^n b_{c,j}x_j,\dots,x_{d-1}+\sum_{j =d}^n b_{d-1,j}x_j).
\end{align}
From these descriptions of $I_{\Lambda(\aaa)}$ and $I_{\Lambda(\bb)}$ it follows that $\Gamma_v \cap U_0$ is cut out by the ideal generated by the $v\times v$ minors of the matrix

\begin{equation*}
\begin{pmatrix}
a_{0,d} -\widetilde{b}_{0,d} & \cdots & a_{0,n} -\widetilde{b}_{0,n}  \\
 \vdots &  & \vdots \\ 
a_{c-1,d} -\widetilde{b}_{c-1,d} & \cdots &  a_{c-1,n} -\widetilde{b}_{c-1,n}
\end{pmatrix}.
\end{equation*}
We can now prove an analogue of Proposition \ref{bigmatrix}. Moreover, using the presentations in (\ref{aaa}) and (\ref{bb}) and arguing as in Proposition \ref{mainone}, \ref{mainonee} we can construct a morphism $\pi^{-1}_{c-1}(U_0) \longrightarrow \HH_{n-c,n-d}^n$ extending the rational map $\Xi$. An argument identical to the one given for Proposition \ref{INJECTIVE} will show that this extends to a bijective morphism $\Xi:\XX_{c-1} \longrightarrow \mathcal{H}_{n-c,n-d}^n$. In a similar manner we may deduce the following results

\begin{customthm}{C'} \label{uniqueborel'} Let $d >c \geq 2$. The component $\HH_{n-c,n-d}^n$ has a unique Borel fixed point. If $n \geq c+d-1 $ the point
\begin{equation*}
I_{n-c,n-d}^n = (x_0,\dots,x_{c-1})(x_0,\dots,x_{d-1})+\sum_{i=0}^{c-1}x_i(x_d,\dots,x_{c+d-2-i})
\end{equation*}
is the unique Borel fixed point on $\mathcal{H}_{n-c,n-d}^n$.
\end{customthm}

Arguing as in Lemma \ref{dacurves}, Theorem \ref{main} and Corollary \ref{mainthree} we obtain
\begin{customthm}{B} \label{maintwo} Let $d > c \geq 2$ and $n \geq c+d-1$. The component $\mathcal{H}_{n-c,n-d}^n$ is smooth and there is an isomorphism
\begin{equation*}
\Xi: \Bl_{\Gamma_{c-1}} \cdots \Bl_{\Gamma_1} (\mathbf{Gr}(n-c,n) \times \mathbf{Gr}(n-d,n)) \longrightarrow \mathcal{H}_{n-c,n-d}^n.
\end{equation*}

If $n < c+d -1$, the morphism $\mathcal{H}_{n-c,n-d}^n \longrightarrow \mathbf{Gr}(2n-c-d+1,n)$ that sends a scheme to its linear span is smooth; the fiber over a point $\Lambda$ is $\mathcal{H}_{n-c,n-d}(\Lambda)$. 
\end{customthm}

\begin{customthm}{D'} \label{primdecomp2} Let $n \geq c+d-1$ and let $Z$ be a subscheme parameterized by $\mathcal{H}_{n-c,n-d}^n$.  Then $Z$ is a pair of planes meeting transversely, or there exists a sequence of integers $1 \leq i_1 < \cdots < i_r \leq c$ and a flag of linear spaces $\Lambda^1 \subseteq \Lambda^{2} \subseteq \cdots \subseteq \Lambda^{r} \subseteq \PP^n$ with $\mathrm{codim}_{\PP^n}(\Lambda^{\ell})= (d+i_{\ell}-1)$ for each $\ell$, such that
\begin{enumerate}
\item If $i_1 >1$ then $Z$ is a union of two planes meeting along $\Lambda^1$ with embedded pure double structures on $\Lambda^{\ell}$ for each $1 \leq \ell \leq r$.
\item If $i_1 =1$ then $Z$ is a codimension $c$-plane with embedded pure double structures on $\Lambda^{\ell}$ for each $1 \leq \ell \leq r$. 
\end{enumerate}
\end{customthm}

\begin{customcor}{E'} \label{2cpointss} Up to projective equivalence, there are exactly $2^c$ subschemes parameterized by $\mathcal{H}_{n-c,n-d}^n$. \end{customcor}

\begin{remark} \label{manycomponents} In \cite{ccn} it was shown that $\mathcal{H}_{n-2,n-2}^n$ meets exactly one other component in $\Hilb^{P_{n-2,n-2}^n(t)} \mathbf{P}^n$ and that this component is smooth. We will give two examples that show these statements are false in general.

The component $\mathcal{H}_{2,2}^{5}$ will meet the component whose general member parameterizes a pair of $2$-planes meeting at a point union an isolated point. It will also meet the component whose general member parameterizes a quadric union an isolated line. 

In \cite[Theorem 3.16]{ritvik} we show that $\Hilb^{P_{n-2,1}^n(t)} \mathbf{P}^n$ is a union of $\HH_{n-2,1}^n$ and a component $\mathcal{Y}$, whose general point parameterizes a line meeting an $(n-2)$-plane union an isolated point. We show that $\mathcal{Y}$ is singular; its singularity is a cone over the Segre embedding of $\PP^1 \times \PP^{n-2} \hookrightarrow \PP^{2(n-1)-1}$. 
\end{remark}

This completes the discussion of the local structure of $\HH_{n-c,n-d}^n$. The next four sections will pertain to its global geometry. As we did in Section \ref{two}, we begin studying divisors on $\HH_{n-c,n-d}^n$ with $c =d = k$ and $n \geq 2k-1$.


\section{Divisors on $\HH_{n-k,n-k}^n$} \label{Section Divisors} \label{four}
In this section we study the Picard group of $\HH_{n-k,n-k}^n$ for $n \geq 2k-1$. We give an explicit description of the divisors $D_i,N_i$ (Remark \ref{newDi}, \ref{newNi}) and describe equations for their pullback along $\Xi|_{U_{k-1}}$.

\begin{notation} \label{lambdak} We will use $\lambda_k$ to denote the coordinate $T^{(k)}_{0,n-k+1}$ on $U_{k-1}$ from Remark \ref{coordinates}. This convention will simplify the formulas for the equations we will obtain.
\end{notation}

The proofs of Theorem \ref{primdecomp} and Lemma \ref{projideal} give explicit equations for the various loci of embedded structures. 

\begin{lemma} \label{embeddedloci} Let $n \geq 2k-1$ and let $Z$ be a subscheme parameterized by $\Xi(U_{k-1})$. Then
\begin{enumerate}
\item $Z$ is a pair of planes meeting transversely if and only if $\lambda_1,\dots,\lambda_{k-1},\mathbf{T}^{(k)}\ne 0$.
\item $Z$ has an embedded $(n-2k+1)$-plane if and only if $\mathbf{T}^{(k)} = 0$.
\item For each $2 \leq i \leq k-1$, $Z$ has an embedded $(n-k+1-i)$-plane if and only if $\lambda_i = 0$.
\item $Z$ is generically non-reduced if and only if $\lambda_1 = 0$.
\end{enumerate}
\end{lemma}

\begin{definition} Consider the sequence of blowups 
$
\XX_{k-1} \overset{\psi_{k-1}}{\longrightarrow} \XX_{k-2} \overset{\psi_{k-2}}{\longrightarrow} \cdots  \overset{\psi_1}{\longrightarrow} \XX_0.
$
For each $i$ let $E_i$ denote the strict transform in $\mathcal{X}_{k-1}$ of the exceptional divisor of $\psi_i$. Let $E_k$ denote the strict transform of $\Gamma_k$.
\end{definition}

\begin{lemma} \label{picardgroup1} Let $n \geq 2k-1$. Then $N^1(\HH_{n-k,n-k}^{n}) = \Cl(\HH_{n-k,n-k}^n) =  \ZZ^k$. In particular, linear equivalence and numerical equivalence for divisors coincide.
\end{lemma}
\begin{proof} 
Since $\HH_{n-k,n-k}^n = \XX_{k-1}/\sss_2$ is a  smooth rational variety, its class group is torsion free. In particular, $N^1(\XX_{k-1}/\sss_2) = \Cl (\XX_{k-1}/\sss_2)$. Thus it suffices to prove that $\Cl(\XX_{k-1}/\sss_2)_{\QQ} := \Cl (\XX_{k-1}/\sss_2)\otimes \QQ$ is isomorphic to $\QQ^k$. By \cite[Example 1.7.6]{fulton} we have $\Cl(\XX_{k-1}/\sss_2)_{\QQ} =\Cl(\XX_{k-1})_{\QQ}^{\sss_2}$. Let $E_{1,0}$ and $E_{0,1}$ be the strict transform, in $\XX_{k-1}$, of $\mathcal{O}_{\XX_0}(1,0)$ and $\mathcal{O}_{\XX_0}(0,1)$, respectively. By \cite[Theorem 8.5]{hartshorne}, $\Cl(\XX_{k-1})_{\QQ}$ is freely generated by $E_1,\dots,E_{k-1},E_{1,0},E_{0,1}$. Since $\sss_2$ fixes $E_i$ and interchanges $E_{1,0}$ with $E_{0,1}$, it follows that
$$
\Cl(\XX_{k-1})^{\sss_2}_{\QQ} = \spn_{\QQ}\{E_1,\dots,E_{k-1},E_{1,0} + E_{0,1}\} \simeq \QQ^k. \qedhere
$$ 
\end{proof}

\begin{definition} Let $(\XX_0)^{\text{trv}} = \XX_0 \setminus \Gamma_k $ denote the open subset of $\XX_0$ consisting of pairs of $(n-k)$-planes such that the two planes in the pair meet transversely. We say that a pair of $(n-k)$-planes meets another plane $\Lambda$ transversely, if each plane in the pair meets $\Lambda$ transversely. 
\end{definition}

We now describe $D_i$ as a scheme theoretic image under $\Xi$.

\begin{remark} \label{newDi} For each $1 \leq i \leq k-1$ consider a flag $\mathcal{F}_i = \{\Lambda_{i-1} \subseteq \Lambda_{2k-1-i}\}$. Let $W_i \subseteq (\XX_0)^{\text{trv}} $ be the open subset consisting of pairs of planes that meet $\Lambda_{2k-1-i}$ transversely. Let $\hat{D}_i$ denote the (scheme theoretic) closure of 
$$
\{Z \in W_i: \dim_{\kk}\spn(\Lambda_{i-1} \cup (Z \cap \Lambda_{2k-1-i}))< 2k-1-i\}
$$ 
in $\XX_0$. Then $D_i$ is the image of the strict transform of $\hat{D}_i$ under the map $\Xi$.

\smallskip
Similarly, given a plane $\Lambda_{k-1}$, let $\hat{D}_k$ be the scheme theoretic closure of 
$$
\{Z \in (\XX_0)^{\text{trv}}: Z \cap \Lambda_{k-1} \ne \emptyset \}
$$
in $\XX_0$. Then $D_k$ is the image of the strict transform of  $\hat{D}_k$ under the map $\Xi$.
\end{remark}

\begin{lemma} \label{equationDi} The loci $D_i$ are divisorial. For $1 \leq i \leq k-1$ let $D_i$ be defined by the flag 
\begin{equation} \label{flagDi}
\Lambda_{i-1} = V(x_{i-1},x_{i+1},\dots,x_n) \subseteq \Lambda_{2k-i-1} = V(x_k,\dots, x_{n-k_{i-2}},x_{n-k_i}-x_{n-k_{i-1}}). 
\end{equation}
Then $\Xi^{\star}(D_i) \cap U_{k-1}$ is cut out by
$
T^{(k-i)}_{i-1,n-k_i} + T^{(k-i)}_{i-1,n-k_i}T^{(k-i)}_{i,n-k_{i-1}} + \lambda_{k-i+1}.
$
\end{lemma}
\begin{proof} Assume $1 \leq i \leq k-1$ and let $D_i$ be defined by the flag (\ref{flagDi}). To show that $D_i$ is a divisor, it suffices to show that $\hat{D}_i \cap W_i$ is a divisor in $W_i$ (notation from Remark \ref{newDi}). By symmetry, it is enough to show that $\hat{D}_i \cap W_i \cap U_0$ is a divisor in $W_i \cap U_0$.

Given a point $(\Lambda(\aaa),\Lambda(\bb)) \in W_i \cap U_0$ we have $(\Lambda(\aaa) \cup \Lambda(\bb)) \cap \Lambda_{2k-1-i} = P \cup Q$ for a pair of $(k-1-i)$-planes, $P$ and $Q$. For each $n-k_{i + 1} \leq j \leq n$ let $p_j$ (respectively $q_j$) denote the point in $P$ (respectively $Q$) obtained by setting $x_j = 1$ and $x_{\ell} = 0$ for all \textit{other} $\ell \geq k$ (there are no such points for $i = k-1$).
Explicitly, 
\begin{align*}
p_j =  (-a_{0,j}: \cdots: -a_{k-1,j} :0 : \cdots: 0: 1 : 0:  \cdots :0) \\
q_j =  (-b_{0,j}: \cdots: -b_{k-1,j} :0 : \cdots: 0: 1 : 0:  \cdots :0).
\end{align*}
Let $p_{n-k_i}$ (respectively $q_{n-k_i}$) denote the point in $P$ (respectively $Q$) obtained by setting $x_{n-k_i} = x_{n-k_{i-1}} =1$ and $x_{\ell} =0$ for all other $\ell \geq k$. Explicitly,
\begin{align*}
p_{n-k_i} =  (-a_{0,n-k_i}-a_{0,n-k_{i-1}}: \cdots: -a_{k-1,n-k_i}-a_{k-1,n-k_{i-1}} 
					:0 : \cdots: 0: 1:1 : 0:  \cdots :0) \\
q_{n-k_i} =  (-b_{0,n-k_i}-b_{0,n-k_{i-1}}: \cdots:-b_{k-1,n-k_i}-b_{k-1,n-k_{i-1}} 
					:0 : \cdots: 0: 1:1 : 0:  \cdots :0). 
\end{align*}
For each $\ell \in \{0,\dots,i-2,i\}$ let $r_{\ell} = V(x_0,\dots,x_{\ell-1},x_{\ell+1},\dots,x_n)$. 

By construction we have, $P = \mathrm{span}(p_{n-k_i},\dots,p_n) $,  $Q = \mathrm{span}(q_{n-k_i},\dots,q_n)$ and $\Lambda_{i-1} = \mathrm{span}(r_0,\dots,r_{i-2},r_i)$. It follows that the points in $\mathrm{span}(\Lambda_{i-1}\cup ((\Lambda(\aaa)\cup \Lambda(\bb)) \cap \Lambda_{2k-1-i}))$ are in the row span of the matrix 
$$
\begin{bmatrix}
q_{n-k_i} & \cdots & q_n & p_{n-k_i} & \cdots &p_n &r_0& \cdots& r_{i-2} &r_{i}
\end{bmatrix}^T.
$$
In particular, $\hat{D}_i \cap W_i\cap U_0$ is the locus where the matrix has rank less than $2k-i$. Let $\epsilon_{l,j} = a_{l,j} - b_{l,j}$ and apply the row operation
$$
\begin{pmatrix}
q_{n-k_i} \\ q_{n-k_{i+1}} \\  \vdots \\ q_n \\ p_{n-k_i} \\ \vdots  \\ p_n  \\ r_0 \\ \vdots \\ r_{i-2}  \\ r_{i}
\end{pmatrix}
\longrightarrow 
\begin{pmatrix}
q_{n-k_i} - p_{n-k_i} - \sum_l (\epsilon_{l,n-k_i} +\epsilon_{l,n-k_{i-1}}) r_l \vspace{2pt} \\ 
q_{n-k_{i+1}} - p_{n-k_{i+1}} - \sum_l \epsilon_{l,n-k_{i+1}}r_l  \vspace{2pt} \\ 
\vdots \\ 
q_{n} - p_{n} - \sum_l \epsilon_{l,n} r_l \vspace{2pt}\\ 
p_{n-k_i} \\ 
\vdots  \\ 
p_n  \\ 
r_0 \\ \vdots \\ r_{i-2}  \\ r_{i}
\end{pmatrix}
$$
and swap the $i$-th column and $(i-1)$-st column. It follows that the locus 
is cut out by the determinant of the submatrix
$$
\begin{pmatrix} 
\epsilon_{i-1,n-k_i} + \epsilon_{i-1,n-k_{i-1}} & \epsilon_{i+1,n-k_i} + \epsilon_{i+1,n-k_{i-1}} & \cdots & \epsilon_{k-1,n-k_i} +\epsilon_{k-1,n-k_{i-1}} & &  & \\
\epsilon_{i-1,n-k_{i+1}} & \epsilon_{i+1,n-k_{i+1}} & \cdots & \epsilon_{k-1,n-k_{i+1}} & &  & \\
\epsilon_{i-1,n-k_{i+2}} & \epsilon_{i+1,n-k_{i+2}} & \cdots & \epsilon_{k-1,n-k_{i+2}} & &  &  \\
\vdots & \vdots & & \vdots & & & \\
\epsilon_{i-1,n} & \epsilon_{i+1,n} & \cdots & \epsilon_{k-1,n} & &  & \\
\end{pmatrix}.
$$
Thus $\hat{D}_i  \cap W_i \cap U_0$ is a divisor and this determinant also cuts out $\hat{D}_i \cap U_0$. 

The strict transform of this determinant cuts out $\Xi^{\star}(D_i)\cap U_{k-1}$. Pulling back this matrix to $U_{k-1}$ and column reducing as in Proposition \ref{bigmatrix} we obtain
\begin{equation*}
\begin{pmatrix} 
\lambda_1\cdots\lambda_{k-i}(T^{(k-i)}_{i-1,n-k_i} + T^{(k-i)}_{i-1,n-k_{i-1}}) & \star & \cdots & \cdots & \star & \star \\\
0 & \lambda_1\cdots\lambda_{k-i-1} & \ddots&  &   & \vdots \\
0 & 0 & \ddots & \ddots &  & \vdots\\
\vdots & \cdots& \ddots &\ddots  & \star & \star \\
0 &\cdots &  & 0 & \lambda_{1}\lambda_2 & \star \\
0 &\cdots &  & 0 & 0 & \lambda_1 \\
\end{pmatrix}.
\end{equation*}
The strict transform of its determinant is $T^{(k-i)}_{i-1,n-k_i} + T^{(k-i)}_{i-1,n-k_{i-1}}$. 
\begin{itemize}[itemsep=0.5ex]
\item If $i >1$ we may use Proposition \ref{bigmatrix} (ii) to rewrite 
$T^{(k-i)}_{i-1,n-k_{i-1}} = \lambda_{k-i+1} + T^{(k-i)}_{i-1,n-k_i}T^{(k-i)}_{i,n-k_{i-1}}$.
\item If $i=1$ we may use Remark \ref{coordinates} to rewrite 
$T^{(k-1)}_{0,n-k+1} = \lambda_k + T_{0,n-k+2}^{(k-1)}T_{1,n-k+1}^{(k-1)} $.
\end{itemize}
In either case, $\Xi^{\star}(D_i) \cap U_{k-1}$ is cut out by the desired equation. Lastly,  $D_k$ is a divisor since $\hat{D}_k$ is the Weil divisor associated to $\mathcal{O}_{\XX_0}(1,1) \in \Pic \XX_0 \simeq \ZZ^2$.
\end{proof}

\begin{cor} \label{equationDi2} Let $0 \leq j < i$. For $1 \leq i \leq k-1$ let $D_i$ be defined by the flag 
\begin{align} \label{flagDi2}
\resizebox{.93\hsize}{!}{$
\Lambda_{i-1} = V(x_{j},x_{i+1},\dots,x_n) \subseteq \Lambda_{2k-i-1} = V(x_k,\dots, x_{n-k_{j-2}},x_{n-k_{j}}-x_{n-k_{j-1}},x_{n-k_{j+1}},\dots,x_{n-k_i})$\footnotemark} 
\end{align}
\footnotetext{if $j=0$ then $k_{j-1} = k_{-1} = k$ is still consistent with our convention, see Remark \ref{knotation}}
and let $D_k$ be defined by the plane
$$
\Lambda_{k-1} = V(x_{j}+x_{n-k_j},x_{k},\dots,x_{n-k_{j-1}},x_{n-k_{j+1}},\dots,x_n).
$$
Then $\Xi^{\star}(D_i) \cap U_{k-1}$ is cut out by a polynomial in the coordinates of Remark \ref{coordinates} that is linear in $\lambda_{k-j}$.
\end{cor}
\begin{proof} Assume $i \leq k-1$ and $j \ne 0$. Imitating the proof of Lemma \ref{equationDi} we see that $\Xi^{\star}(D_i) \cap U_{k-1}$ is cut out by  $T^{(k-i)}_{j,n-k_j} + T^{(k-i)}_{j,n-k_{j-1}}$. To express this in terms of our desired coordinates we will use the relation $T^{(\ell)}_{p,q} = T^{(\ell)}_{p,n-\ell+1}T^{(\ell)}_{k-\ell,q} + \lambda_{\ell+1}T^{(\ell+1)}_{p,q}$ which is true for any $q\leq n-k_p$ and any $p < k-\ell$ and $\ell < k-1$ (proof of Proposition \ref{bigmatrix}). 
Repeatedly applying this relation we obtain the following expressions
$$
T_{j,n-k_j}^{(k-i)} = \sum_{\ell = k-i}^{k-j-1} \lambda_{k-i+1}\cdots \lambda_{\ell} T^{(\ell)}_{j,n-\ell+1}T^{(\ell)}_{k-\ell,n-k_j} + \lambda_{k-i+1}\cdots\lambda_{k-j}
$$
and 
\begin{equation} \label{ijrelation} 
T_{j,q}^{(k-i)} = \sum_{\ell = k-i}^{k-j-1} \lambda_{k-i+1}\cdots \lambda_{\ell} T^{(\ell)}_{j,n-\ell+1}T^{(\ell)}_{k-\ell,q} 
				+ \lambda_{k-i+1}\cdots\lambda_{k-j}T^{(k-j)}_{j,q}
\end{equation}
for any $q < n-k_j$. Thus $T_{j,q}^{(k-i)}$, as a polynomial in the coordinates of Remark \ref{coordinates}, is linear in $\lambda_{k-j}$ for all $q \leq n-k_j$. This implies  $\Xi^{\star}(D_i) \cap U_{k-1}$ is linear in $\lambda_{k-j}$. 

\smallskip
Assume $i \leq k-1$ and $j =0$. Most of the argument from the previous paragraph still applies in this case. In particular, $\Xi^{\star}(D_i) \cap U_{k-1}$ is cut out by  $T^{(k-i)}_{0,n-k+1} + T^{(k-i)}_{0,n-k}$ and we have
\begin{equation} \label{ijrelation'} 
T_{0,q}^{(k-i)} = \sum_{\ell = k-i}^{k-2} \lambda_{k-i+1}\cdots \lambda_{\ell} T^{(\ell)}_{0,n-\ell+1}T^{(\ell)}_{k-\ell,q} + \lambda_{k-i+1}\cdots\lambda_{k-1}T^{(k-1)}_{0,q}
\end{equation}
for all $q \leq n-k+1 = n -k_0$. Notice that $T_{0,q}^{(k-1)} = T^{(k)}_{0,q} + T^{(k-1)}_{0,n-k+2}T_{1,q}^{(k-1)}$ for all $q \leq n-k+1$ and $T^{(k)}_{0,n-k+1} = \lambda_k$ (Remark \ref{coordinates}). Substituting this into (\ref{ijrelation'}) we see that $T^{(k-i)}_{0,n-k+1} + T^{(k-i)}_{0,n-k}$ is linear in $\lambda_k$.

\smallskip
Finally assume $i = k$. The locus of points $(\Lambda(\aaa),\Lambda(\bb)) \in U_0 $ meeting $\Lambda_{k-1}$ is clearly cut out by  $(a_{j,n-k_j}-1)(b_{j,n-k_j}-1)$. The pullback of this equation to $U_{k-1}$, which coincides with the strict transform, defines $\Xi^{\star}(D_k)$. If $j\ne 0$ we can use (\ref{ijrelation}) to deduce that
$$
(a_{j,n-k_j}-1)(b_{j,n-k_j} -1) =  \big(b_{j,n-k_j}+ \sum_{\ell = 1}^{k-j-1} \lambda_1\cdots \lambda_{\ell} T^{(\ell)}_{j,n-\ell+1}T^{(\ell)}_{k-\ell,n-k_j} + \lambda_{1}\cdots\lambda_{k-j}-1\big)(b_{j,n-k_j} -1).
$$
This expression is linear in $\lambda_{k-j}$. If $j =0$ we can argue in the previous paragraph and deduce linearity in $\lambda_k$. This completes the proof. 
\end{proof}

Here is an alternate description of $N_i$.

\begin{remark} \label{newNi} For each $1 \leq i \leq k-1$, let $N_i = \Xi(E_i)$. If $n=2k-1$ we let $N_k = \Xi(E_k)$. If $n > 2k-1$, let $\hat{N}_k$ denote the closure in $\XX_0$, of the locus of pairs of planes in $\XX_0^{\text{trv}}$ where the intersection of the two planes meets a fixed $\Lambda_{2k-1}$. Then $N_k$ is the image of the strict transform of $\hat{N}_k$ under $\Xi$. \end{remark}

In the next lemma we abuse notation and use "=" to mean equality as divisor classes.

\begin{lemma} \label{equationNi} Let $n \geq 2k-1$. The loci $N_i$ are divisorial. Moreover, we have 
\begin{enumerate}
\item $\Xi^{\star}(N_1) = 2E_1$.
\item $\Xi^{\star}(N_i) = E_i$ for $2 \leq i \leq k-1$.
\item If $n=2k-1$ then $\Xi^{\star}(N_k) = E_k$ and $\Xi^{\star}(N_k) \cap U_{k-1}$ is cut out by $\lambda_k$.
\item If $n > 2k-1$ let $\Lambda_{2k-1} = V(x_k,\dots,x_{n-k})$ be the plane defining $N_k$. Then $\Xi^{\star}(N_k) \cap U_{k-1}$ is cut out by $\lambda_k$.
\end{enumerate}
\end{lemma}
\begin{proof} Assume $1 \leq i \leq k-1$. Remark \ref{newNi} implies that the $N_i$ are divisors. Items (i), (ii) and the first half of (iii) follow from the fact that $\Xi$ is a finite, degree $2$ map branched along $N_1$ (although not phrased this way, it is part of the proof of Proposition \ref{mainonee}), see \cite[Chapter 1.7]{fulton}. The rest of item (iii) is a consequence of  Lemma \ref{embeddedloci} (ii).

Now assume $n>2k-1$ and let $\hat{N}_k$ be as in Remark \ref{newNi}. To show that $N_k$ is a divisor it is enough to show that $\hat{N}_k \cap \XX_{0}^{\text{trv}} \cap U_0$ is a divisor in $\XX_0^{\text{trv}} \cap U_0$.  Given a point $(\Lambda(\aaa), \Lambda(\bb)) \in \XX_0^{\text{trv}} \cap U_0$, the intersection of the two planes is $\Lambda(\aaa) \cap \Lambda(\bb) = V(\{\sum_{j=k}^{n}(a_{i,j}-b_{i,j})x_{j},y_i\}_{0\leq i \leq k-1})$. Thus the locus of points in $\XX_0^{\text{trv}} \cap U_0$ satisfying $(\Lambda(\aaa) \cap \Lambda(\bb)) \cap \Lambda_{2k-1} \ne \emptyset$ is cut out by the determinant of
$$
\begin{pmatrix} 
a_{0,n-k+1}-b_{0,n-k+1} &  \cdots & a_{k-1,n-k+1} - b_{k-1,n-k+1} \\
\vdots &  & \vdots \\
a_{0,n} -b_{0,n} &\cdots & a_{k-1,n} - b_{k-1,n}  \\
\end{pmatrix}
$$
Column reducing as in Proposition \ref{bigmatrix} (ii) and taking the strict transform gives item (iv).
\end{proof}


\section{Birational geometry of $\mathcal{H}_{n-k,n-k}^n$ for $n \geq 2k-1$} \label{Section Part I} \label{five}

This section is devoted to the proof of Proposition \ref{cones}. For the rest of the section we will assume  $n \geq 2k-1$. We begin by constructing  two families of curves and computing their intersection numbers with $D_i$ and $N_i$.

Roughly speaking, the first family of curves will fix a pair of planes and vary the embedded structures while the second family will vary the planes and fix the embedded structures.

\begin{definition} \label{family3}  For each $1 \leq j \leq k-1$, define the curve $C_j:\mathbf{P}^1 \to \mathcal{H}_{n-k,n-k}^n$ by
\begin{equation*} 
\begin{aligned}
C_j(s:t) =  I_{\Lambda}I_{\Lambda'} + (sx_{j-1}x_{n-k_j}-tx_jx_{n-k_{j-1}}) +  \sum_{p=0}^{j-2} x_p(x_{n-k_{p+1}},\dots,x_{n-k_j})
\end{aligned}
\end{equation*}
with $\Lambda = V(x_0,\dots,x_{k-1})$ and $\Lambda' = V(x_0,\dots,x_j,x_{j+1}+x_{n-k_{j+1}},\dots,x_{k-1}+x_{n})$.
\end{definition}

\begin{remark} \label{decompC}  Theorem \ref{primdecomp} shows that $C_j(s:t)$ is projectively equivalent to (\ref{projideal}) with 
$$
\mu_{1}=\cdots = \mu_{k-j-1}=1, 
	\quad \mu_{k-j} = 0,
	\quad \mu_{k-j+1} = \begin{cases} \frac{t}{s} \text{ if } s \ne 0 \\ 0 \text{ if } s=0 \end{cases}, 
	 \quad \mu_{k-j+2}=\cdots=\mu_k = 0.
$$
It also shows that for $j \leq k-2$, the general member of $C_j$ is a pair of $(n-k)$-planes meeting along a pencil of embedded $(n-2k+j+1)$-planes and containing fixed embedded $(n-2k+\ell)$-planes for all $1 \leq \ell \leq j-1$, while $C_{k-1}$ is a pencil of generically non-reduced $(n-k)$-planes. If $(s:t) = (1:0),(0:1)$, the corresponding subscheme has an embedded $(n-2k+j)$-plane.
\end{remark}

\begin{definition} \label{family2} Let $0 \leq j \leq k-1$. Let $\Lambda = V(x_0,\dots,x_{k-1})$ and consider the pencil of $(n-k)$-planes
$
\Lambda'(s:t) = V(x_0,\dots,x_{j-1},sx_{j}+tx_{n-k_j},x_{j+1}+x_{n-k_{j+1}},\dots,x_{k-1}+x_{n}). 
$
Define the curve $B_j: \mathbf{P}^1 \to \mathcal{H}_{n-k,n-k}^n$ by
\begin{equation*} 
\begin{aligned}
B_j(s:t)=I_{\Lambda}I_{\Lambda'(s:t)} + (x_px_{n-k_q}-x_qx_{n-k_p})_{0 \leq p < q \leq j-1} + (x_0,\dots,x_{j-1})x_{n-k_j}.
\end{aligned}
\end{equation*}
\end{definition}

\begin{remark} \label{decompB} Theorem \ref{primdecomp} shows that $B_j(s:t)$ is projectively equivalent to (\ref{projideal}) with 
$$
\mu_{1}=\cdots = \mu_{k-j-1}=1, 
 	\quad \mu_{k-j} = \begin{cases} \frac{t}{s} \text{ if } s \ne 0 \\ 1 \text{ if } s=0 \end{cases}, 
	\quad \mu_{k-j+1}  = 0,
	\quad \mu_{k-j+2} = \cdots = \mu_k=1.
$$

If $(s:t) \ne (1:0)$, then $B_0(s:t)$ is a pair of $(n-k)$-planes meeting transversely while $B_j(s:t)$ a pair of $(n-k)$-planes with a pure embedded $(n-2k+j)$-plane for $j> 0$. Moreover, the embedded $(n-2k+j)$-plane is fixed along the curve. 

If $(s:t) = (1:0)$, the corresponding subscheme has an embedded $(n-2k+j+1)$-plane. Note that $B_{k-1}(1:0)$ is, more precisely, a generically non-reduced $(n-k)$-plane. 
\end{remark}

Before we determine the intersection numbers we need to compute a few linear spans. We begin with notation that will be used a great deal in the following Lemmas.

\begin{notation} We use $C^{\dagger}_j(s:t)$ and $B^{\dagger}_j(s:t)$ to denote the subschemes of $\PP^n$ cut out by $C_j(s:t)$ and $B_j(s:t)$, respectively. Given an ideal $J \subseteq S$, let $\sat(J)$ denote its saturation with respect to $(x_0,\dots,x_n)$ and let $J(1)$ denote the ideal generated by the linear forms in $J$.
\end{notation}

\begin{lemma} \label{linearspan} Let $ 1 \leq i \leq j \leq k-1$ and let $\Lambda_{2k-i-1} = V(x_k,x_{k+1},\dots, x_{n-k_{i-2}},x_{n-k_i}-x_{n-k_{i-1}})$. For any $(s:t) \in \PP^1$, if $i \ne j$  the linear span of $C^{\dagger}_j(s:t) \cap \Lambda_{2k-i-1}$ is  
$$V(x_0,\dots,x_{i-1},x_k,\dots,x_{n-k_{i-2}},x_{n-k_i}-x_{n-k_{i-1}})
$$
and if $i = j$ the linear span of $C^{\dagger}_i(s:t)\cap \Lambda_{2k-i-1} $ is 
$$
V(x_0,\dots,x_{i-2},sx_{i-1}-tx_i,x_k,\dots,x_{n-k_{i-2}},x_{n-k_i}-x_{n-k_{i-1}}).
$$
\end{lemma}
\begin{proof} Let $\Lambda = \Lambda_{2k-i-1}$ and note that the linear span of $C^{\dagger}_j(s:t) \cap \Lambda$ is cut out by $\sat(C_j(s:t) + I_{\Lambda})(1)$. Assume $i < j$. It is straigthtforward to see that $x_{\ell}(x_0,\dots,x_n) \subseteq C_j(s:t) + I_{\Lambda}$ for every $0 \leq \ell \leq i-1$. Thus we have
\begin{align*}
\sat(C_j(s:t) + I_\Lambda) \supseteq & \quad I_{\Lambda}+ (x_0,\dots,x_{i-1})  +(x_i,\dots,x_{k-1})(x_i,\dots,x_j,x_{j+1}+x_{n-k_{j+1}},\dots,x_{k-1}+x_n) \\ 
			& \quad + (sx_{j-1}x_{n-k_j}-tx_jx_{n-k_{j-1}}) + \sum_{p=i}^{j-2}x_p(x_{n-k_{p+1}},\dots,x_{n-k_j}) \\
		 = & \quad \mathfrak{Q}.
\end{align*}
Moreover, it is clear that $\mathfrak{Q}(d) = (C_j(s:t) +I_\Lambda)(d)$ for all $d \geq 2$.  Thus if we show that $\mathfrak{Q}$ is saturated then $\mathfrak{Q} = \sat(C_j(s:t) +I_\Lambda)$, and this would give the desired linear span. If we write $ \mathfrak{Q} = I_\Lambda+ (x_0,\dots,x_{i-1}) + \mathfrak{Q}'$, it suffices to show that quadratic portion, $\mathfrak{Q}'$, is saturated. But notice that $ \mathfrak{Q}'$ is projectively equivalent to an ideal of the form  (\ref{projideal}) (for reasons similar to Remark \ref{decompC}). It follows from Lemma \ref{saturated} that $\mathfrak{Q}$ is saturated. The case of $i=j$ is analogous.
\end{proof}

\begin{remark} \label{spandim} Here are two simple facts about linear spans:
\begin{enumerate}
\item If $\Lambda_p$ and $\Lambda_q$ are disjoint linear spaces in $\PP^n$ then $\dim_{\kk} \spn (\Lambda_p \cup \Lambda_q) = p+q+1$.
\item $\spn(Y_1 \cup Y_2) = \spn(\spn Y_1 \cup \spn Y_2)$ for any subschemes $Y_1,Y_2 \subseteq \PP^n$. 
\end{enumerate}
The first fact is clear and the second follows from the following chain of equalities,
$$
I_{Y_1 \cup Y_2}(1) = (I_{Y_1} \cap I_{Y_2})(1) = (I_{Y_1}(1) \cap I_{Y_2}(1))(1).
$$
\end{remark}

\begin{lemma} \label{intersection1} Let $1 \leq i \leq k$ and $1 \leq j \leq k-1$. We have the following intersection numbers
\begin{enumerate}
\item $D_i \cdot C_j =0$ whenever $i \ne j$,
\item $D_i \cdot C_i = 1$ for all $i \leq k-1$.
\end{enumerate}
\end{lemma}
\begin{proof} Assume $i> j$. Since the dimension of any embedded subscheme of $C^{\dagger}_j(s:t)$ is at most $n-2k+j+1$, a generic $(2k-1-i)$-plane will not intersect any embedded subscheme of $C^{\dagger}_j(s:t)$. If $i < k$, the intersection of $C^{\dagger}_j(s:t)$ with a generic $\Lambda_{2k-1-i}$  is a pair of skew $(k-1-i)$-planes. Moreover, these skew planes are independent of $(s:t)$ and thus 
$$
\mathrm{span}\, (C^{\dagger}_j(s:t) \cap \Lambda_{2k-1-i}) \simeq \mathbf{P}^{2k-2i-1}
$$ is independent of $(s:t)$. As a consequence, we may choose an $(i-1)$-plane $\Lambda_{i-1} \subseteq \Lambda_{2k-1-i}$ that does not meet the $\mathbf{P}^{2k-2i-1}$. It follows from Remark \ref{spandim} that
$$
\dim_{\kk} \mathrm{span}\, (\Lambda_{i-1} \cup (C_j^{\dagger}(s:t) \cap \Lambda_{2k-1-i})) = 2k-1-i.
$$ 
If we use the flag $\{\Lambda_{i-1} \subseteq \Lambda_{2k-1-i}\}$ to define $D_i$ we see that $D_i \cdot C_j = 0$. Similarly, if $i = k$ and $\Lambda_{k-1}$ is generic we have that $C_j^{\dagger}(s:t) \cap \Lambda_{k-1} = \emptyset$. Thus $D_k \cdot C_j = 0$.

Assume $i <j$ and  let $\Lambda_{2k-i-1} = V(x_k,x_{k+1},\dots, x_{n-k_{i-2}},x_{n-k_i}-x_{n-k_{i-1}})$. By Lemma \ref{linearspan} we have that
$$
\mathrm{span}\, (C^{\dagger}_j(s:t) \cap \Lambda_{2k-1-i}) = V(x_0,\dots,x_{i-1},x_k,x_{k+1},\dots,x_{n-k_{i-2}},x_{n-k_i}-x_{n-k_{i-1}})\simeq \PP^{2k-2i-1}
$$ 
is fixed and independent of $(s:t)$. As done in the previous paragraph, if we choose a general $\Lambda_{i-1}$ inside $\Lambda_{2k-1-i}$ to define $D_i$, then $D_i \cdot C_j = 0$. This completes the proof of item (i).

Assume $i= j$ and let the flag $\{\Lambda_{i-1} \subseteq \Lambda_{2k-1-i}\}$ in (\ref{flagDi}) define $D_i$. By Lemma \ref{linearspan} we have that 
$$
\mathrm{span}\, (C_i^{\dagger}(s:t) \cap \Lambda_{2k-1-i}) = V(x_0,\dots,x_{i-2},sx_{i-1}-tx_i,x_k,x_{k+1},\dots,x_{n-k_{i-2}},x_{n-k_i}-x_{n-k_{i-1}})
$$
Thus, if $t \ne 0$, the linear span of $(C_i^{\dagger}(1:t) \cap \Lambda_{2k-i-1}) \cup \Lambda_{i-1}$ is all of $\Lambda_{2k-i-1}$. If $t=0$, the linear span of $(C_i^{\dagger}(1:0) \cap \Lambda_{2k-i-1}) \cup \Lambda_{i-1}$ is $\Lambda_{2k-i-1} \cap V(x_{i-1})$. Thus $D_i \cap C_i$ is supported on the point $Z_0= C_i(1:0)$.

Let $\tilde{C}_i$ denote the closure in $\XX_{k-1}$ of the curve, $\mathbf{A}^1 \hookrightarrow U_{k-1}$ obtained by setting $\lambda_{1},\dots,\lambda_{k-i-1} = 1$, $\lambda_{k-i+1} = t$ and all the other coordinates of Remark \ref{coordinates} to $0$. Since $\Xi(\tilde{C}_i)|_{U_{k-1}} = C_i(1:t)$ it follows that $\Xi(\tilde{C}_i) =C_i$. In particular $\tilde{C}_i \cap \Xi^{\star}(D_i)$ is supported at a unique point $\tilde{Z}_0 \in \Xi^{-1}(Z_0)$. Since $\Xi^{\star}(D_i)$ is linear in $\lambda_{k-i+1}$ (Lemma \ref{equationDi}), it follows that $\Xi^{\star}(D_i)$ and $\tilde{C}_i$ intersect transversely at $\tilde{Z}_0$. Using the push-pull formula we conclude that $C_i\cdot D_i = \Xi_{\star}\tilde{C}_i \cdot D_i = \Xi_{\star}(\tilde{C}_i \cdot \Xi^{\star}(D_i)) =1$.
\end{proof}

\begin{lemma} \label{intersection2} Let $1 \leq i \leq k$ and $0 \leq j \leq k-1$. We have the following intersection numbers
\begin{enumerate}
\item $D_i \cdot B_j =0$ for all $i \leq j$,
\item $D_i \cdot B_j = 1$ for all $i > j$.
\end{enumerate}
\end{lemma} 
\begin{proof} Assume $i \leq j$ and let $\Lambda_{2k-1-i} = V(x_k,\dots, x_{n-k_{i-2}},x_{n-k_i}-x_{n-k_{i-1}})$. Arguing as in Lemma \ref{linearspan} we see that 
$$
\mathrm{span}\, (\Lambda_{2k-1-i} \cap B^{\dagger}_j(s:t)) = V(x_0,\dots,x_{i-1},x_k,x_{k+1},\dots,x_{n-k_{i-1}}) \simeq \PP^{2k-2i-1}
$$ 
is independent of $(s:t)$. Arguing as in Lemma \ref{intersection1}  we deduce item (i). 

Assume that $j < i \leq k-1$ and let $\{\Lambda_{i-1} \subseteq \Lambda_{2k-1-i}\}$ be the flag (\ref{flagDi2}) defining $D_i$. Then 
$B_j^{\dagger}(s:t) \cap \Lambda_{2k-i-1}$ is a disjoint pair of $(k-i-1)$-planes defined by
\begin{align*}
&(x_0,\dots,x_{j-1},sx_j+tx_{n-k_j},x_{j+1},\dots,x_i,x_{i+1}+x_{n-k_{i+1}},\dots,x_{k-1}+x_n, \\
& \quad \quad \quad x_k,x_{k+1},\dots, x_{n-k_{j-2}},x_{n-k_{j}}-x_{n-k_{j-1}},x_{n-k_{j+1}},\dots,x_{n-k_i}) \cap \\
&\quad (x_0,\dots,x_{n-k_{j-2}},x_{n-k_{j}}-x_{n-k_{j-1}},x_{n-k_{j+1}},\dots,x_{n-k_i}).
\end{align*}
For  $t \ne 0$, the linear span of $(B_j^{\dagger}(s:t) \cap \Lambda_{2k-i-1}) \cup \Lambda_{i-1}$ is all of $\Lambda_{2k-i-1}$. On the other hand if $t= 0$, the linear span of $(B_j^{\dagger}(s:t) \cap \Lambda_{2k-i-1}) \cup \Lambda_{i-1}$ is  $\Lambda_{2k-1-i} \cap V(x_j)$. Thus  $D_i \cap B_j$ is supported at the point $Z_0 = B_j(1:0)$.

Let $\tilde{B}_j$ denote the closure in $\XX_{k-1}$ of the curve, $\mathbf{A}^1 \hookrightarrow U_{k-1}$ obtained by setting $\lambda_1=\dots =\lambda_{k-j-1} =1$, $\lambda_{k-j} = t$, $\lambda_{k-j+2} =\dots =\lambda_k =1$ and all the other coordinates of Remark \ref{coordinates} to $0$. Since $\Xi(\tilde{B}_j)|_{U_{k-1}} = B_j(1:t)$ we have $\Xi(\tilde{B}_j) = B_j$. Thus $\tilde{B}_j \cap \Xi^{\star}(D_i)$ is supported at a unique point $\tilde{Z}_0 \in \Xi^{-1}(Z_0)$. Since $\Xi^{\star}(D_i)$ is linear in $\lambda_{k-j}$ (Corollary \ref{equationDi2}), it follows that $\Xi^{\star}(D_i)$ and $\tilde{B}_j$ intersect transversely at $\tilde{Z}_0$. Using the push-pull formula we conclude that $B_j\cdot D_i = \Xi_{\star}\tilde{B}_j \cdot D_i = \Xi_{\star}(\tilde{B}_j \cdot \Xi^{\star}(D_i)) =1$. 

Now assume $j < i=k$ and let $\Lambda_{k-1} = V(x_{j}+x_{n-k_j},x_{k},\dots,x_{n-k_{j-1}},x_{n-k_{j+1}},\dots,x_n)$ be the plane defining $D_k$. It is evident that  $B_j\cap D_k$  is supported at the point $Z_{1,1} = B_j(1:1)$. Once again, $\tilde{B}_j$ (defined in the previous paragraph) and $\Xi^{\star}(D_k)$ will meet at a unique point $\tilde{Z}_{1,1} \in \Xi^{-1}(Z_{1,1})$. Since $\Xi^{\star}(D_k)$ is linear in $\lambda_{k-j}$ (Corollary \ref{equationDi2}) we see that $\tilde{B}_j$ meets $\Xi^{\star}(D_k)$ transversely at $\tilde{Z}_{1,1}$. Once again we conclude using the push-pull formula.
\end{proof}

\begin{lemma}  \label{intersection3} We have the following intersection numbers,
\begin{enumerate}
\item $N_i \cdot C_j = 0$ for each $1 \leq i \leq k-1$ and all $1 \leq j \leq k-i-1$,
\item $N_i \cdot B_j = 0$ for each $1 \leq i \leq k$ and all $j \ne k-i,k-i+1$,
\item $N_i \cdot C_{k-i+1} = 2$ for each $2 \leq i \leq k$,
\item $N_1 \cdot B_{k-1} = 2$ and $N_i \cdot B_{k-i} = 1$  for $2 \leq i \leq k$.
\end{enumerate}
\end{lemma} 
\begin{proof} Item (i) and item (ii), except for the case of $i=k$,  follow from the definition of the $N_i$ and the description of the embedded subschemes in Remark \ref{decompC} and Remark \ref{decompB}. We will deal with the case of $i=k$ in the last paragraph. For the rest of the proof let $Z_0 = C_{k-i+1}(1:0)$ and $Z_{\infty} = C_{k-i+1}(0:1)$. We will also use the curves $\tilde{C}_{k-i+1}$ and $\tilde{B}_j$ defined in Lemma \ref{intersection1}. In particular, let $\tilde{Z}_0, \tilde{Z}_{\infty} \in \tilde{C}_{k-i+1}$ be such that $\Xi(\tilde{Z}_0) = Z_0$ and $\Xi(\tilde{Z}_{\infty}) = Z_{\infty}$.

Assume $2 \leq i \leq k-1$. Since $N_{i}$ is the locus of subschemes containing an embedded $(n-k+1-i)$-plane, it meets the curve $C_{k-i+1}$ at $Z_0$ and $Z_{\infty}$. Thus $\tilde{C}_{k-i+1}$ meets $E_i$ at $\tilde{Z}_0$ and $\tilde{Z}_{\infty}$. Using Lemma \ref{equationNi} (ii), we obtain 
$$
N_i \cdot C_{k-i+1} =\Xi_{\star}(\tilde{C}_{k-i+1} \cdot \Xi^{\star}(N_i)) = \tilde{C}_{k-i+1} \cdot E_i 
	= (\tilde{C}_{k-i+1} \cdot E_i)|_{\tilde{Z}_0} + (\tilde{C}_{k-i+1} \cdot E_i)|_{\tilde{Z}_{\infty}}.
$$
Since $\tilde{Z}_0 \in U_{k-1}$ and $E_i$ is cut out by $\lambda_{i}$, $\tilde{C}_{k-i+1}$ meets $E_i$ transversely at $\tilde{Z}_0$. Symmetrically, $\tilde{C}_{k-i+1}$ will also meet $E_i$ transversally at $\tilde{Z}_{\infty}$. To see the latter statement, consider the projective transformation $g \in \GL(n+1)$ that interchanges $x_j$ with $x_{j-1}$, interchanges $x_{n-k_j}$ with $x_{n-k_{j-1}}$ and fixes the other coordinates. It follows from the definition that $g(C_{k-i+1}) = C_{k-i+1}$ and $g$ interchanges $Z_0$ with $Z_{\infty}$. Since intersection multiplicity is invariant under automorphisms of $\HH_{n-k,n-k}^n$ we obtain 
$$
(N_i \cdot C_{k-i+1})|_{Z_{\infty}} = \left(g(N_i) \cdot g(C_{k-i+1})\right)|_{g(Z_{\infty})} = N_i \cdot C_{k-i+1}|_{Z_0}  = (E_i \cdot \tilde{C}_{k-i+1})|_{\tilde{Z}_0} = 1.
$$
This proves item (iii) for $i \ne k$.

Since $N_1$ is the locus of generically non-reduced subschemes, it meets the curve $B_{k-1}$ at $B_{k-1}(1:0)$.  Using Lemma \ref{equationNi} (i) we obtain $N_1 \cdot B_{k-1} =\Xi_{\star}(\tilde{B}_{k-1} \cdot \Xi^{\star}(N_1)) =  2\tilde{B}_{k-1}\cdot E_1= 2$.  Similarly,  using Lemma \ref{equationNi} we obtain $N_i \cdot B_{k-i} = 1$  for all $2 \leq i \leq k-1$. This finishes item (iv) for $i \ne k$

Finally, assume $i=k$ and let $\Lambda_{2k-1} = V(x_k,\dots,x_{n-k})$ be the plane defining $N_k$ (if $n > 2k-1$). By Lemma \ref{equationNi} (iii), (iv) we see that $\Xi^{\star}(N_k)$ meets $\tilde{C}_1$ at $Z_0$ and possibly also at $Z_{\infty}$ (since the latter does not lie in $U_{k-1}$). Moreover, $\Xi^{\star}(N_k)$ meets $\tilde{C}_1$ transversely at  $\tilde{Z}_0$. We may argue as in the previous paragraph to show that $\Xi^{\star}(N_k)$ also meets $\tilde{C}_1$ transversely at $\tilde{Z}_{\infty}$. Indeed, the projective transformation $g$ fixes $N_k$. This is clear if $n = 2k-1$ and the case of $n>2k-1$ follows from the fact that $g$ fixes $\Lambda_{2k-1}$. Thus $N_k \cdot C_{1} = (N_k \cdot C_{1})|_{Z_{0}} + (N_k \cdot C_{1})|_{Z_{\infty}} = 2(N_k \cdot C_{1})|_{Z_0} = 2$, completing the proof of item (iii).  For items (ii) and (iv) we argue similarly using the following projective transformation: $g' \in \GL(n+1)$ that maps $x_{n-k_j} \mapsto x_{n-k_j}+ x_j$ and fixes the other coordinates. It is straightforward to verify that $g'(B_j) = B_j$, $g'(B_j(0:1)) = B_j(1:1)$ and $g'$ fixes $N_k$ (since $g'$ fixes $\Lambda_{2k-1}$). This implies
$$
(N_k\cdot B_j)|_{B_j(0:1)} = \left(g'(N_k) \cdot g'(B_j)\right)|_{g'(B_j(0:1))} = (N_k \cdot B_j)|_{B_j(1:1)} 
= 0
$$
for $j \ne 1$. Thus, we may compute $\Xi^{\star}(N_k) \cdot \tilde{B}_j$ along $U_{k-1}$ to obtain the desired results.
\end{proof}

\begin{prop} \label{lincomb} Let $1 \leq i \leq k$. Then we have
\begin{itemize}[itemsep=0.5ex]
\item $N_1 = 2D_k-2D_{k-1}$, 
\item $N_i = 2D_{k-i+1} - D_{k-i} - D_{k-i+2}$ \, for all $2 \leq i \leq k-1$,
\item $N_k = 2D_1 - D_2$.
\end{itemize}
\end{prop} 
\begin{proof}
By Lemma \ref{picardgroup1}, Lemma \ref{intersection1} and Lemma \ref{intersection2} we see that $N^1(\mathcal{H}_{n-k,n-k}^n)$ is generated by $\{D_1,\dots,D_k\}$. This allows us to write $N_i = \sum_{\ell=1}^k\epsilon_{i,\ell}D_\ell$ for some $\epsilon_{i,\ell} \in \ZZ$. Using Lemmas \ref{intersection1} - \ref{intersection3} we obtain
\begin{itemize}[itemsep=0.5ex]
\item $N_1 \cdot C_\ell = \epsilon_{1,\ell} = 0$ for  $\ell \leq k-2$, 
\item $N_1\cdot B_{k-1} = \epsilon_{1,k} =2$ and $N_1\cdot B_{k-2} = \epsilon_{1,k-1} + \epsilon_{1,k} = 0$.
\end{itemize}
This immediately implies $N_1 = 2D_k-2D_{k-1}$. For each $2 \leq i \leq k$ we obtain
\begin{itemize}[itemsep=0.5ex]
\item $N_i \cdot B_j = \sum_{\ell = j+1}^k\epsilon_{i,\ell} = 0$ for  $j \ne k-i, k-i+1$
\item $N_i \cdot B_{k-i} = \sum_{\ell = k-i+1}^k\epsilon_{i,\ell} =1$ and $N_i \cdot C_{k-i+1} = \epsilon_{i,k-i+1} = 2$.
\end{itemize}
If $i \ne k $, we obtain $\epsilon_{i,k-i} = -1, \epsilon_{i,k-i+1} = 2, \epsilon_{i,k-i+2} = -1$, and $\epsilon_{i,\ell}=0$ for other $\ell$. If $i=k$ we obtain $\epsilon_{k,1} = 2$, $\epsilon_{k,2}=-1$ and $\epsilon_{i,\ell} =0$ for other $\ell$. This completes the proof.
\end{proof}



\begin{prop} \label{cones} Let $k\geq 2$ and $n \geq 2k-1$. Then we have
$$
\emph{Eff}(\HH_{n-k,n-k}^n) = \langle N_1,\dots,N_k \rangle \quad \text{and} \quad 
\emph{Nef}(\HH_{n-k,n-k}^n) = \langle D_1,\dots,D_k \rangle.
$$
Moreover, $\HH_{n-k,n-k}^n$ is Fano if and only if either $k=3$ and $n =5$, or $k \ne 3$ and $n \in \{2k-1,2k\}$. 
\end{prop}
\begin{proof} 
It is clear that the divisors $N_1,\dots,N_k$ are effective and generate $N^1(\HH_{n-k,n-k}^n)$. To conclude that the effective cone is generated by $N_1,\dots,N_k$, it is enough to show that any $\mathbf{R}$-divisor $N = \sum_{i=1}^{k} \epsilon_iN_i$, with some $\epsilon_j <0$, is not effective. Let $A_j: \PP^1 \hookrightarrow \HH_{n-k,n-k}^n$ denote any curve such that for $(s:t) \ne (1:0)$, $A_j(s:t)$ is a pair of $(n-k)$-planes meeting transversely while $A_j(1:0)$ it is a pair of $(n-k)$-planes with a pure embedded $(n-k+1-j)$-plane if $j>1$ and generically non-reduced if $j=1$. Clearly, $A_j \cdot N_i =0$ for $i \ne j$ and $A_j \cdot N_j >0$. Since $N \cdot A_j = \epsilon_j < 0$ and $A_j$ is not contained in the support of $N$, we see that $N$ cannot be an effective divisor.

By varying the flags it is easy to see that each of the $D_i$ is base point free; thus it is also nef. Similar to the previous paragraph, to show that the nef cone gone is generated by $D_1,\dots,D_k$, it is enough to show that any $\mathbf{R}$-divisor $D = \sum_{i=1}^{k} \epsilon_iD_i$, with some $\epsilon_j <0$, is not nef. If $j \ne k$, we have $D \cdot C_j = \epsilon_j < 0$ and if $j=k$ we have $D \cdot B_{k-1} = \epsilon_k < 0$. Thus $D$ is not nef. 

We will now compute the canonical divisor of $\HH_{n-k,n-k}^n$ using the branched cover $\Xi: \XX_{k-1} \to \HH_{n-k,n-k}^n$. By  \cite[Exercise 8.5b]{hartshorne} and \cite[Exercise 10.10]{eisenbud} we may write
\begin{equation*}
K_{\mathcal{X}_{k-1}} = \sum_{j=1}^{k-1}((k-j+1)(n-k-j+2)-1)E_j - (n+1)\hat{D}_k
\end{equation*}
where $\hat{D}_k$ is the strict transform of $\mathcal{O}_{\XX_0}(1,1)$ (Remark \ref{newDi}). Note that the canonical divisor of $\XX_0$ is $\mathcal{O}_{\XX_0}(-n-1,-n-1)$. Let $K_{\mathcal{H}_{n-k,n-k}^n} = \epsilon_1N_1 + \cdots +\epsilon_{k-1}N_{k-1} + \epsilon_kD_k$ for some $\epsilon_i \in \QQ$. Hurwitz's theorem implies that $K_{\mathcal{X}_{k-1}} = \Xi^{\star}(K_{\mathcal{H}_{n-k,n-k}^n})+E_1$. Using this and Lemma \ref{equationNi} we obtain
\begin{align*}
2\epsilon_1E_1 + \sum_{j=2}^{k-1}\epsilon_j E_j + \epsilon_k\hat{D}_k = \Xi^{\star}(K_{\mathcal{H}_{n-k,n-k}^n}) & =
			(k(n-k+1)-2)E_1+  \\
			& \quad \quad \sum_{j=2}^{k-1}((k-j+1)(n-k-j+2)-1)E_j - (n+1)\hat{D}_k.
\end{align*}
Let $\tilde{\epsilon}_j = (k-j+1)(n-k-j+2)-1$ and using Proposition \ref{lincomb} we obtain
$$
K_{\HH_{n-k,n-k}^n} =\frac{1}{2}(\tilde{\epsilon}_1 -1)(2D_k-2D_{k-1})
			+ \sum_{j=2}^{k-1}\tilde{\epsilon}_j(2D_{k-j+1}-D_{k-j}-D_{k-j+2}) -(n+1)D_k. 
$$
For $k =2,3$ the above expression simplifies to
$$
K_{\HH_{n-2,n-2}^n} = (4-2n)D_1+ (n-5)D_2, \quad
K_{\HH_{n-3,n-3}^n}  = (7-2n)D_1+(n-6)D_2-2D_3.
$$
If $k\geq 4$ we can rewrite the expression as follows:
\begin{align*}
K_{\HH_{n-k,n-k}^n} &= (\tilde{\epsilon}_1-1)(D_k-D_{k-1}) - (n+1)D_k +
							 \sum_{j=2}^{k-3}(2\tilde{\epsilon}_{j+1}-\tilde{\epsilon}_j-\tilde{\epsilon}_{j+2})D_{k-j} \\
	&   \quad\quad\quad	-\tilde{\epsilon}_2D_k + (2\tilde{\epsilon}_2-\tilde{\epsilon}_3)D_{k-1} 
							+(2\tilde{\epsilon}_{k-1}-\tilde{\epsilon}_{k-2})D_2 - \tilde{\epsilon}_{k-1}D_1\\
	&=  (\tilde{\epsilon}_1-\tilde{\epsilon}_2-n-2)D_k + (2\tilde{\epsilon}_2-\tilde{\epsilon}_3-\tilde{\epsilon}_1+1)D_{k-1} 
					+ \sum_{j=2}^{k-3}(2\tilde{\epsilon}_{j+1}-\tilde{\epsilon}_j-\tilde{\epsilon}_{j+2})D_{k-j}  \\
	& \quad\quad\quad +(2\tilde{\epsilon}_{k-1}-\tilde{\epsilon}_{k-2})D_2 - \tilde{\epsilon}_{k-1}D_1.
\end{align*}
Since $2\tilde{\epsilon}_{j+1}-\tilde{\epsilon}_j-\tilde{\epsilon}_{j+2} = -2$ for all $j$ we obtain
\begin{equation*}
K_{\HH_{n-k,n-k}^n} = (4k-5-2n)D_1 + (n-2k-1)D_2 - 2D_3 - 2D_4 -\cdots -2D_{k-2} - D_{k-1} - 2D_k.
\end{equation*}

Since the ample cone is the interior of the nef cone, we see that $-K_{\HH_{n-2,n-2}^n}$ is ample if and only if $n=3,4$ and that $-K_{\HH_{n-3,n-3}^n}$ is ample precisely when $n = 5$. If $k \geq 4$, $-K_{\HH_{n-k,n-k}^n}$ is ample if and only if $n = 2k-1,2k$.
\end{proof}

\section{Birational geometry of $\mathcal{H}_{k-1,k-1}^n$ for $n > 2k-1$} \label{Section Part II} \label{six}

This section is devoted to the proof of Theorem \ref{MORI}. We will show that  $\HH_{k-1,k-1}^{n}$ is Fano, and thus a Mori dream space. By constructing a contraction from $\HH_{k-1,k-1}^{n}$ to $\HH_{n-k,n-k}^n$ (Proposition \ref{hilbertchow}) we will also deduce that  $\HH_{n-k,n-k}^n$ is a Mori dream space.
\smallskip

\begin{notation} \label{c&d} In this section we will primarily be interested in the case when the pair of planes do not span all of $\PP^n$. By swapping the roles of codimension and dimension, the components we are interested in are of the form $\HH_{k-1,k-1}^n$ with $n > 2k-1$. 
\end{notation}

Corollary \ref{mainthree} states that for $n > 2k-1$, the morphism $\rho:\mathcal{H}_{k-1,k-1}^n \longrightarrow \mathbf{Gr}(2k-1,n)$ that sends a scheme to its linear span is smooth; the fiber over a point $\Lambda$ is $\mathcal{H}_{k-1,k-1}(\Lambda)$.

\begin{remark} \label{upthere}  Let $W = \Spec \kk[f_{2k,j},\dots,f_{n,j}]_{0 \leq j \leq 2k-1}$ be a neighbourhood of $\Lambda = V(x_{2k},\dots,x_n) \in \mathbf{Gr}(2k-1,n)$ such that its $\kk$-points are identified with 
$$
V(x_{2k} + \sum_{j=0}^{2k-1}f_{2k,j}x_j,\dots,x_{n} + \sum_{j=0}^{2k-1}f_{n,j}x_j). 
$$
Then the open subset $\rho^{-1}(W)$ is naturally isomorphic to $W \times \HH_{k-1,k-1}(\Lambda)$.
\end{remark}

\begin{lemma} \label{picardgroup2} Let $n > 2k-1$. Then $N^1(\mathcal{H}_{k-1,k-1}^{n}) = \ZZ^{k+1}$.
\end{lemma}
\begin{proof} As explained in Lemma \ref{picardgroup1}, since $\HH_{k-1,k-1}^n$ is rational and smooth, it suffices to compute $N^1(\HH_{k-1,k-1}^n) \otimes \QQ$ which equals $\Pic(\mathcal{H}_{k-1,k-1}^{n})\otimes \QQ = H^{2}(\mathcal{H}_{k-1,k-1}^n,\mathbf{Q})$. By Corollary \ref{mainthree} we have a smooth morphism $\mathcal{H}_{k-1,k-1}^n \longrightarrow \mathbf{Gr}(2k-1,n)$ with fibers isomorphic to $\mathcal{H}_{k-1,k-1}^{2k-1}$. Since the base of this morphism is simply connected, we may apply the Leray-Hirsch theorem \cite[Theorem 7.33]{voisin} and Lemma \ref{picardgroup1} to deduce that $H^{2}(\mathcal{H}_{k-1,k-1}^n,\mathbf{Q}) \simeq  \mathbf{Q}^{k+1}$.
\end{proof}

Using the fibration $\rho$ and Remark \ref{upthere} one can easily verify that the loci $D_i', N_i', F$ are divisorial. We now define the curves inside $\HH_{k-1,k-1}^n$; all but two of them come from curves lying inside $\HH_{k-1,k-1}^{2k-1}$.

\begin{definition} Let  $\Lambda = V(x_{2k},\dots,x_n)$. For each relevant $j$, let $A_j', B_j', C_j'$ be the images of $A_j,B_j, C_j$ (Definition \ref{family3}, \ref{family2}, Proposition \ref{cones}) under the inclusion $\rho^{-1}(\Lambda) = \mathcal{H}_{k-1,k-1}(\Lambda) \hookrightarrow \HH_{k-1,k-1}^n $, respectively.
\end{definition}

\begin{definition} \label{Ys}  Let $\Lambda' = V(x_k,\dots,x_{n})$ and let
$
\Lambda(s:t) = V(x_0,\dots,x_{k-1},sx_{2k}+tx_k,x_{2k+1},\dots,x_n)
$
be a pencil of $(k-1)$-planes disjoint from $\Lambda'$. 
Define the curve $Y_1:\mathbf{P}^1 \to \HH_{k-1,k-1}^n$ by $(s:t) \mapsto \Lambda(s:t) \cup \Lambda'$. Explicitly
$$
Y_1(s:t) = (sx_{2k}+tx_{k},x_{2k+1},\dots,x_n) + (x_0,\dots,x_{k-1})(x_{k},\dots,x_{2k-1}).
$$
Define the curve $Y_2:\mathbf{P}^1 \to \HH_{k-1,k-1}^n$ by
\begin{align*}
Y_2(s:t)  & = (sx_{2k}+tx_{0},x_{2k+1},\dots,x_n) + (x_1,\dots,x_{k-1})(x_{k+1},\dots,x_{2k-1}) \\
	& \quad \quad + (x_0,x_{2k})^2 + (x_0,x_{2k})(x_1,\dots,x_{k-1},x_{k+1},\dots,x_{2k-1}).
\end{align*}
\end{definition}

\begin{remark}  Let $\Lambda = V(x_0,\dots,x_{k-1},x_{2k},\dots,x_n)$ and $\Lambda' = V(x_{0},x_{k+1},\dots,x_n)$ be a pair of $(k-1)$-planes meeting along a point. Then we have
$$
Y_2(s:t) = I_{\Lambda} \cap I_{\Lambda'}\cap ((x_0,x_{2k})^2,sx_{2k}+tx_0,x_1,\dots,x_{k-1},x_{k+1},\dots,x_{2k-1},x_{2k+1},\dots, x_{n}).
$$
In particular, $Y_2$ is a  pair of fixed $(k-1)$-planes with a pencil of embedded points.
\end{remark}

\begin{lemma} \label{movingY2} $Y_2$ is a moving curve in $N_k'$ i.e. its deformations span $N_k'$.
\end{lemma}
\begin{proof} The general subscheme parameterized by $N_k'$ is a pair of $(k-1)$-planes meeting along an embedded point. By Corollary \ref{mainthree} and Theorem \ref{primdecomp}, up to projectively equivalence, such a subscheme is cut out by 
$$
(x_0,\dots,x_{k-1},x_{2k},\dots,x_n) \cap (x_0,x_{k+1},\dots,x_n) \cap (x_0^2,x_1,\dots,x_{k-1},x_{k+1},\dots,x_n) = Y_2(1:0)
$$
In particular, the $\GL(n+1)$ orbit of $Y_2$ covers a dense subset of $N_k'$. 
\end{proof}

\begin{lemma} \label{interstellar} For all pairs of relevant indices $i,j$ (the ones appearing in Lemma \ref{intersection1}, \ref{intersection2}, \ref{intersection3}), the intersection numbers of $D_i',N_i'$ with $B_j',C_j'$ are the same as the intersection numbers of $D_i,N_i$ with $B_j,C_j$, respectively.
\end{lemma}
\begin{proof} We will only verify $D_i' \cdot C_j' = D_i \cdot C_j$ for $ 1 \leq i,j  \leq k-1 $; the other cases are analogous. Let $\Lambda = V(x_{2k},\dots,x_n)$ be a fixed $(2k-1)$-plane. Let $D_i'$ be defined by a flag $\mathcal{F}_i' = \{\Lambda_{n-2k+i} \subseteq \Lambda_{n-i}\}$, where the flag is chosen to satisfy the following two properties:
\begin{itemize}
\item $\Lambda$ is transverse to each element of the flag $\mathcal{F}_i'$,
\item Let $D_i \subseteq \HH_{k-1,k-1}(\Lambda)$ be defined by the flag $\mathcal{F}_i = \{\Lambda_{n-2k+i} \cap \Lambda \subseteq \Lambda_{n-i} \cap \Lambda\}$. Then either $D_i \cap C_j = \emptyset$ if $i \ne j$ or $D_i$ is transverse to $C_j$ if $i=j$.
\end{itemize}

Let $W$ be the open neighbourhood of $\Lambda$ from Remark \ref{upthere}. The first bullet point implies that every element of $W$ is transverse to the flag $\mathcal{F}_i'$.  It follows that $D_i'|_{\rho^{-1}(W)} = W \times D_i$ and $C_j' = \{\Lambda \} \times C_j$. Thus we have  $D_i' \cdot C_j' = D_i'|_{\rho^{-1}(W)} \cdot C_j' = D_i \cdot C_j$.
\end{proof}

\begin{lemma} \label{intersection5} We have the following intersection numbers
\begin{enumerate}
\item  $D_i' \cdot Y_2 = N_i' \cdot Y_1 =0$ for all $1 \leq i \leq k$,
\item $N_i' \cdot Y_2 = 0$ for all $1 \leq i \leq k-1$,
\item $D_i'\cdot Y_1 =1$ for all $1 \leq i \leq k$,
\item $F \cdot Y_1 = F \cdot Y_2 = 1$.
\end{enumerate}
\end{lemma}
\begin{proof} Items (i) and (ii) are clear from the definition of the divisors. 

Let $1 \leq i \leq k$, $\Lambda = V(x_{2k},\dots,x_n)$ and $W$ be as in Remark \ref{upthere}. We may choose a flag $\mathcal{F}_i'$ to define $D_i'$ so that the following properties are satisfied:
\begin{itemize}
\item $\Lambda$ is transverse to each element of the flag $\mathcal{F}_i'$,
\item $D_i' \cap Y_1$ is supported at $Z_0 = Y_1(1:0)$.
\end{itemize}

Let $W' =  \Spec \kk[\epsilon_1,\dots,\epsilon_{k^2}] \subseteq \HH_{k-1,k-1}(\Lambda)$ be any affine open containing the image of $Z_0$ in $\HH_{k-1,k-1}(\Lambda)$. Then $W \times W'$ is identified with an open neighbourhood of $Z_0 \in \HH_{k-1,k-1}^n$. Along this open set, $Y_1$ is the curve obtained by setting $f_{2k,k} = t$, $f_{i,j} = 0 $ for other $i,j$, and $\epsilon_i = \delta_i$ for some constants $\delta_i \in \kk$. On the other hand, $D_i' = W \times (D_i \cap W')$ where $D_i$ is the divisor defined by the flag $\mathcal{F}_i' \cap \Lambda$. It immediately follows that $D_i'$ meets $Y_1$ transversely at $Z_0$  inside $W \times W'$; this proves item (iii).

For item (iv), we will only verify $F\cdot Y_1 = 1$ as the other case is similar. Let $F$ be defined by the $(n-2k)$-plane, $V(x_0,\dots,x_{k-1},x_{k+1},\dots,x_{2k})$. It follows that $F \cap Y_1$ is also supported at $Z_0$. Moreover, along $W \times W'$, $F$ is cut out by the function $f_{2k,k}$. Combining this with the equation of $Y_1$ along $W \times W'$ we see that $F$ meets $Y_1$ transversely at $Z_0$.
\end{proof}

\begin{prop} \label{lincomb2} Let $k \geq 2$ and $n > 2k-1$. Then we have,
$$
\emph{Eff}(\HH_{k-1,k-1}^n) = \langle N_1',\dots,N_k',F \rangle \quad \text{and} \quad 
\emph{Nef}(\HH_{k-1,k-1}^n)  = \langle D_1',\dots,D_k',F \rangle.
$$
Moreover we have, 
\begin{itemize}[itemsep=0.5ex]
\item $N_1' = 2D_k'-2D_{k-1}'$,
\item $N_i' = 2D_{k-i+1}' - D_{k-i}' - D_{k-i+2}'$ for all  $2 \leq i \leq k-1$,
\item $N_k' = 2D_1' - D_2' - F$.
\end{itemize}
\end{prop}
\begin{proof} Using the intersection numbers with the curves $\{C_1',\dots,C_k',Y_2\}$ and arguing as in Proposition \ref{lincomb}, \ref{cones} we see that $N^1(\HH_{k-1,k-1}^n) $ and $\Nef(\HH_{k-1,k-1}^n)$ are both generated by $D_1',\dots,D_k',F$. Using the curves $\{A_1',\dots,A_k',Y_1\}$ and arguing as in Proposition \ref{cones}, we see that $N_1',\dots,N_k',F$ generate the effective cone.

By Proposition \ref{lincomb} and Remark \ref{interstellar} there exists $\epsilon_i \in \QQ$ such that
\begin{itemize}[itemsep=0.5ex]
\item $N_1' =2D_k' -2D_{k-1}' + \epsilon_1F$,
\item $N_i' = 2D_{k-i+1}' - D_{k-i}' - D_{k-i+2}' +\epsilon_iF$ \, for all \, $2 \leq i \leq k-1$, 
\item $N_k' = 2D_1'-D_2'+\epsilon_kF$.
\end{itemize}
Intersecting these divisors with $Y_1,Y_2$ and using Lemma \ref{intersection5} we obtain $\epsilon_1,\dots,\epsilon_{k-1}=0$ and $\epsilon_k=-1$. \end{proof}

We are now ready to relate $\HH_{k-1,k-1}^n$ with $\HH_{n-k,n-k}^n$.

\begin{prop} \label{hilbertchow} There is a morphism  $\Psi: \HH_{k-1,k-1}^n \longrightarrow \HH_{n-k,n-k}^n$ with exceptional locus $N_k'$. Moreover, $N_k'$ is a $\PP^{n-2k+1}$-fibration over $\Psi(N_k')$. Geometrically, $\Psi$ "forgets" the embedded points.
\end{prop}
\begin{proof} Given an $(n+1)$-dimensional vector space $V$, let 
$$
\Gamma_i(\PP V) = \{(\Lambda,\Lambda'): \dim (\Lambda \cap \Lambda') \geq k-i \} \subseteq \Gr(k-1,\PP V)^2.
$$
The Hilbert-Chow morphism induces a birational morphism, $\HH_{k-1,k-1}(\PP V) \longrightarrow \Sym^2 \Gr(k-1,\PP V)$ \cite[Theorem 6.3]{kollar}. Let $\widebar{\Gamma}_i(\PP V)$ denote the image of $\Gamma_i(\PP V)$ in $\Sym^2 \Gr(k-1,\PP V)$.  Since the pullback of each $\widebar{\Gamma}_i(\PP V)$ is $N_i'$, we obtain a morphism 
$$
\Psi_1:\HH_{k-1,k-1}^n \longrightarrow \Bl_{\widebar{\Gamma}_{k-1}(\PP V)}\cdots\Bl_{\widebar{\Gamma}_{1}(\PP V)}\Sym^2 \Gr(k-1,\PP V).
$$ 

There is an isomorphism $\mathbf{Gr}(k-1,\PP V)^2 \simeq \Gr(n-k, (\PP V)^{\star})^2$ induced by map $\Lambda \mapsto \Lambda^{\star}$ that sends a linear space to its dual variety. This isomorphism maps $\Gamma_i(\PP V)$ to $\Gamma_i$ (Definition \ref{gamma}) and thus maps $\widebar{\Gamma}_i(\PP V)$ to $\widebar{\Gamma}_i$ after quotienting by $\sss_2$. Therefore we obtain an isomorphism
\begin{align*}
\Psi_2:\Bl_{\widebar{\Gamma}_{k-1}(\PP V)}\cdots\Bl_{\widebar{\Gamma}_{1}(\PP V)}\Sym^2 \Gr(k-1,\PP V) 
		& \xrightarrow{\simeq} \Bl_{\widebar{\Gamma}_{k-1}} \cdots \Bl_{\widebar{\Gamma}_1} \Sym^2 \mathbf{Gr}(n-k,n) \\
		& = \HH_{n-k,n-k}((\PP V)^{\star}).
\end{align*}
Let $\Psi = \Psi_2 \circ \Psi_1$. One can directly check that $\Psi^{\star}(D_i) = D_i'$ for all $i$ and $\Psi^{\star}(N_i) = N_i'$ for $1 \leq i \leq k-1$. 

To show that $\Psi$ contracts $N_k'$, it is enough to show that $\Psi$ contracts $Y_2$ (Lemma \ref{movingY2}). Using Lemma \ref{intersection5} we obtain $\Psi_{\star}Y_2 \cdot D_i =\Psi_{\star}(Y_2 \cdot \Psi^{\star}(D_i)) = \Psi_{\star}(Y_2 \cdot D_i')= 0$ for all $i$. Since $D_1,\dots,D_k$ generates the nef-cone of $\HH_{n-k,n-k}^n$ we must have $\Psi_{\star}Y_2 = 0$, i.e. $\Psi$ contracts $Y_2$. 

Conversely, let $C$ be any curve contracted by $\Psi$. If $C \cdot D_i' \ne 0$ for some $i$, we would have $\Psi_{\star}C \cdot D_i = \Psi_{\star}(C \cdot D_i')\ne 0$, proving that $\Psi$ does not contract $C$. Thus we may assume $C \cdot D_i' = 0$ for all $i$. Since $\{D_i'\}_{i}\cup F$ generates the nef-cone of $\HH_{k-1,k-1}^n$ we must have $F \cdot C > 0$. Using Proposition \ref{lincomb2} we obtain $N_k' \cdot C = - F \cdot C < 0$, i.e. $C$ lies inside $N_k'$.

Lastly, we need to verify that $N_k'$ is a $\PP^{n-2k+1}$-fibration over $\Psi(N_k')$. Up to projective equivalence, it is enough to verify that the fiber of $\Psi_1$ over 
$
Z =V(x_0,\dots,x_{k-1},x_{2k},\dots,x_n) \cup V(x_0,x_{k+1},\dots,x_n)
$
is isomorphic to $\PP^{n-2k+1}$, c.f. Example \ref{exampledual}. Let $H = \text{span}_{\kk}\{x_0,x_{2k},\dots,x_n\}$. Similar to the proof of Lemma \ref{movingY2}, any subscheme parameterized by $\HH_{k-1,k-1}^n$ and supported on $Z$ is cut out by
\begin{equation} \label{Hideal}
\resizebox{.93\hsize}{!}{$
(x_0,\dots,x_{k-1},x_{2k},\dots,x_n) \cap (x_0,x_{k+1},\dots,x_n) \cap \left[(x_1,\dots,x_{k-1},x_{k+1},\dots,x_{2k-1})+(H') + (H'')^2 \right]$}
\end{equation}
where $H' \in \text{Gr}(n-2k+1,H)$ and $H'' \subseteq H$ is chosen so that $H' \oplus H'' = H$. Notice that for a fixed $H'$, all choices of $H''$ give the same ideal as (\ref{Hideal}). It follows that the $\Psi_1^{-1}(Z)$ is paramaterized by $\Gr(n-2k,\PP H) \simeq \PP^{n-2k+1}$.
\end{proof}

\begin{example} \label{exampledual} Consider $X \subseteq \PP^{4}$ cut out by $(x_0,x_1,x_4) \cap (x_0,x_3,x_4) \cap (x_0^2,x_1,x_3,x_4)$. This is a pair of lines meeting along an embedded point. Let $x_0^{\star},\dots,x_4^{\star}$ be the dual coordinates on $(\PP^4)^{\star}$. We can trace the image of $X$ under the map $\Psi:\HH_{1,1}(\PP^4) \to \HH_{2,2}((\PP^4)^{\star})$ as follows:
\begin{align*}
(x_0,x_1,x_4) \cap (x_0,x_3,x_4) \cap (x_0^2,x_1,x_3,x_4) 
						& \overset{\Psi_1}{\mapsto} \, (x_0,x_1,x_4) \cap (x_0,x_3,x_4) \\
	& \overset{\Psi_2}{\mapsto} \, \text{point in } \HH_{2,2}^4 \text{ corresponding to } (x^{\star}_2,x^{\star}_3) \cap (x^{\star}_1,x^{\star}_2)\\
						 & = \, (x^{\star}_2,x^{\star}_3) \cdot(x^{\star}_1,x^{\star}_2)  \\
						 & = \,  (x^{\star}_2,x^{\star}_3) \cap (x^{\star}_1,x^{\star}_2) \cap ((x^{\star}_2)^2,x^{\star}_1,x^{\star}_3).
\end{align*}
\end{example}

\begin{prop}  \label{fanothree} Let $k \geq 2$ and $n > 2k-1$. The component $\HH_{k-1,k-1}^n$ is Fano.
\end{prop}
\begin{proof} Using Proposition \ref{hilbertchow} and the canonical divisor in Proposition $\ref{cones}$ we deduce that
\begin{align*}
K_{\HH_{k-1,k-1}^n} & =  \Psi^{\star}K_{\HH_{n-k,n-k}^n}+(n-2k+1)N_k' \\
			& = \Psi^{\star}K_{\HH_{n-k,n-k}^n}+ (n-2k+1)(2D_1'-D_2'-F) \\
			& = \begin{cases}
				 -3D_1' -2D_2' - 2D_3' -\cdots -2D_{k-2}' - D_{k-1}' - 2D_k' -(n-2k+1)F & \text{if } \, k \geq 4, \\
		 		-3D_1' -D_2' - 2D_3' - (n-5)F & \text{if } \, k = 3, \\ 
				-2D_1' -2D_2' -(n-3)F & \text{if } \, k = 2.
	 	\end{cases}
\end{align*}
The first equality is a modification of \cite[Exercise 8.5]{hartshorne} combined with the fact that the codimension of $\Psi(N_k')$ in $\HH_{n-k,n-k}^n$ is $n-2k+2$. It follows from Proposition \ref{lincomb2} that $-K_{\HH_{k-1,k-1}^n}$ is ample in all cases;  thus $\HH_{k-1,k-1}^n$ is always Fano. 
\end{proof}
Here is the the main theorem of the paper:

\begin{thm} \label{MORI} The components $\HH_{k-1,k-1}^n$ and $\HH_{n-k,n-k}^n$ are Mori dream spaces. 
\end{thm}
\begin{proof} This follows immediately from Proposition \ref{cones}, \ref{hilbertchow} , \ref{fanothree} and the subsequent two facts:
\begin{enumerate}
\item A smooth Fano variety is a Mori dream space \cite[Corollary 4.9]{mckernan}, 
\item Let $f:X \to Y$ be a surjective morphism of smooth, projective varieties. If $X$ is a Mori dream space, then so is $Y$ \cite[Theorem 1.1]{okawa}. \qedhere
\end{enumerate} 
\end{proof}


\section{Birational geometry of $\HH_{n-c,n-d}^n$ and $\HH_{c-1,d-1}^n$} \label{last bro} \label{seven}
In this section we explain how the proofs of Section \ref{Section Divisors}, \ref{Section Part I} and \ref{Section Part II} carry over, almost identically, to the case when the pair of planes are of different dimension. In particular, the definition of the divisors and curves, and computations of their  intersection numbers, including transversality, are very similar. Thus we will omit most of the proofs and indicate all the required modifications.



We begin by defining divisors analogous to the ones in Definition \ref{Di} and \ref{Ni} when the pair of planes span $\PP^n$.

\begin{definition}  Let $n \geq c+d-1$. For each $1 \leq i \leq c-1$ and a choice of a flag of linear spaces $\{\Lambda_{i-1} \subseteq \Lambda_{c+d-1-i}\}$, let $D_i$ denote the divisor class of the locus of subschemes $Z \in \HH_{n-c,n-d}^n$, for which the linear span of $\Lambda_{i-1} \cup (Z\cap \Lambda_{c+d-1-i})$ has dimension less than $c+d-1-i$.
\end{definition}

\begin{definition} Let $n \geq c+d-1$.  Let $D_c^{(1)}$ denote the closure of the locus of subschemes supported on two distinct planes for which the $(n-d)$-plane meets a fixed $\Lambda_{d-1}$. Let $D_c^{(2)}$ denote the closure of the locus of subschemes supported on two distinct planes for which the $(n-c)$-plane meets a fixed $\Lambda_{c-1}$.
\end{definition}

\begin{remark} The divisors $D_c^{(1)}$ and $D_c^{(2)}$ are the Weil divisors associated to the strict transforms, under $\Xi$, of $\mathcal{O}_{\XX_0}(0,1)$ and $\mathcal{O}_{\XX_0}(1,0)$, respectively. Here $\XX_0 = \Gr(n-c,n) \times \Gr(n-d,n)$.

\end{remark}

\begin{definition} Let $n \geq c+d-1$.  For each $1 \leq i \leq c-1$, let $N_i$ denote the divisor class of the locus of subschemes in $\HH_{n-c,n-d}^n$ with an embedded $(n-d+1-i)$-plane. If $n =c+d-1$ let $N_c$ denote the divisor class of the locus of subschemes with an embedded point. If $n > c+d-1$ let $N_c$ denote the class of the closure of the locus of pairs of planes meeting transversely, where the intersection of the two planes meets a fixed $\Lambda_{c+d-1}$.
\end{definition}

We can easily modify the curves in Definition \ref{family3}, \ref{family2} to obtain curves in $\HH_{n-c,n-d}^n$. However, this time we can have two variations, depending on whether the $(n-c)$-plane or $(n-d)$-plane is fixed along the curve. 

\begin{definition} For each $1 \leq j \leq c-1$, define the curve $C_j:\mathbf{P}^1 \to \mathcal{H}_{n-c,n-d}^n$ by
\begin{equation*} 
\begin{aligned}
C_j(s:t) =  I_{\Lambda}I_{\Lambda'} + (sx_{j-1}x_{n-c_j}-tx_jx_{n-c_{j-1}}) +  \sum_{p=0}^{j-2} x_p(x_{n-c_{p+1}},\dots,x_{n-c_j})
\end{aligned}
\end{equation*}
with $\Lambda = V(x_0,\dots,x_{d-1})$ and $\Lambda' = V(x_0,\dots,x_j,x_{j+1}+x_{n-c_{j+1}},\dots,x_{c-1}+x_{n})$.\footnote{Analogous to the notation $k_j$, we define $c_j = c-1-j$ and $d_j = d-1-j$.}
\end{definition}

\begin{definition} For each $0 \leq j \leq c-1$ consider the pencils
$$
\Lambda(s:t) = V(x_0,\dots,x_{j-1},sx_{j}+tx_{n-c_j},x_{j+1}+x_{n-c_{j+1}},\dots,x_{c-1}+x_n)
$$ 
and 
$$
\Lambda'(s:t) = V(x_0,\dots,x_{j-1},sx_{j}+tx_{n-d_j},x_{j+1}+x_{n-d_{j+1}},\dots,x_{d-1}+x_n).
$$

\smallskip
Define the curves $B_j^{(1)}: \mathbf{P}^1 \to \mathcal{H}_{n-c,n-d}^n$ and $B_{j}^{(2)}:\mathbf{P}^1 \to \mathcal{H}_{n-c,n-d}^n$ by
\begin{equation*} 
\begin{aligned}
B_{j}^{(1)}(s:t)=(x_0,\dots,x_{c-1})I_{\Lambda'(s:t)} + (x_px_{n-d_q}-x_qx_{n-d_p})_{0 \leq p < q \leq j-1} + (x_0,\dots,x_{j-1})x_{n-d_j}.
\end{aligned}
\end{equation*}
and
\begin{equation*} 
\begin{aligned}
B_{j}^{(2)}(s:t)=I_{\Lambda(s:t)}(x_0,\dots,x_{d-1}) + (x_px_{n-c_q}-x_qx_{n-c_p})_{0 \leq p < q \leq j-1} + (x_0,\dots,x_{j-1})x_{n-c_j}.
\end{aligned}
\end{equation*}
\end{definition}

Here are the analogues of Lemmas \ref{intersection1} - \ref{intersection3}. To prove these, one first constructs an open set on $\HH_{n-c,n-d}^n$ analogous to $U_{k-1}$ as described in Section \ref{three}.
Then we proceed as in Section \ref{Section Divisors} and describe equations for $D_i$ and $N_i$ along this open set.

\begin{lemma} Let $1 \leq i \leq c-1$ and $0 \leq j \leq c-1$. We have the following intersection numbers,
\begin{enumerate}
\item $D_i \cdot C_i = 1$ and $D_i \cdot C_j =0$ for all $i \ne j$,
\item $D_c^{(1)} \cdot C_i = D_c^{(2)} \cdot C_i = 0$ for all $i$,
\item $D_i \cdot B_j^{(1)} = D_i \cdot B_j^{(2)} =0$ for all $i \leq j$,
\item $D_i \cdot B_j^{(1)} = D_i \cdot B_j^{(2)} = 1$ for all $i > j$,
\item $D_c^{(1)} \cdot B_j^{(1)} = D_c^{(2)}\cdot B_j^{(2)} = 1$ and $D_c^{(1)} \cdot B_j^{(2)} = D_c^{(2)} \cdot B_j^{(1)} = 0$ for all $j$.
\end{enumerate}
\end{lemma}

\begin{lemma} We have the following intersection numbers,
\begin{enumerate}
\item $N_i \cdot C_j = 0$ for each $1 \leq i \leq c-1$ and all $1 \leq j \leq c-1-i$,
\item $N_i \cdot B_j^{(1)} = N_i \cdot B_j^{(2)} = 0$ for each $1 \leq i \leq c$ and all $j \ne c-i,c-i+1$,
\item $N_i \cdot C_{c-i+1} = 2$ for each $2 \leq i \leq c$,
\item $N_i \cdot B_{c-i}^{(1)} = N_i \cdot B_{c-i}^{(2)} = 1$  for each $1 \leq i \leq c$.
\end{enumerate}
\end{lemma}

Just as in Proposition \ref{lincomb} we have,

\begin{prop} Let $1 \leq i \leq c$. If $c=2$ we have,
$$
N_1 = D_2^{(1)} + D_2^{(2)} -D_{1} \, \, \text{and} \,  \, N_2 = 2D_1 - D_2^{(1)} - D_2^{(2)}.
$$
If $c \geq 3$ we have,
\begin{enumerate}[itemsep=0.5ex]
\item $N_1 = D_c^{(1)} + D_c^{(2)} -D_{c-1}$, 
\item $N_2 = 2D_{c-1}-D_{c-2} - D_c^{(1)}-D_c^{(2)}$,
\item $N_i = 2D_{c-i+1} - D_{c-i} - D_{c-i+2}$ \, for all $3 \leq i \leq c-1$,
\item $N_c = 2D_1 - D_2$.
\end{enumerate}
\end{prop}

\begin{prop} \label{cones2} Let $d >c \geq 2$ and $n \geq c+d-1$. We have,
$$
\emph{Eff}(\HH_{n-c,n-d}^n) = \langle N_1,\dots,N_c,D_{c}^{(1)},D_{c}^{(2)} \rangle \quad \text{and} \quad \emph{Nef}(\HH_{n-c,n-d}^n) = \langle D_1,\dots,D_{c-1},D_{c}^{(1)},D_{c}^{(2)} \rangle.
$$
Moreover, 
\begin{enumerate}
\item If $c=2$, then only $\HH_{d-1,1}^{d+1}, \HH_{d,2}^{d+2},\dots,\HH_{2d-3,d-1}^{2d-1}$ are Fano,
\item If $c \geq 3$, then only $\mathcal{H}_{d-1,c-1}^{c+d-1}$ and $\mathcal{H}_{d,c}^{c+d}$ are Fano.
\end{enumerate}
\end{prop}
\begin{proof} 
The verification of the effective and nef cone is similar to Proposition \ref{cones}. Using the formula of the canonical divisor of a blowup and arguing as in  Proposition \ref{lincomb}  we obtain, $K_{\HH_{n-2,n-d}^n} =(2d-2n-1)D_1+ (n-2d)(D_2^{(1)}+D_2^{(2)})$ and 
$$
K_{\HH_{n-c,n-d}^n} = (2c+2d-2n-5)D_1 + (n-d-c-1)D_2- \sum_{j=3}^{c-1}2D_{c-j} + (c-d-1)(D_c^{(1)}+D_c^{(2)})
$$
for $c \geq 3$.  Therefore $\HH_{n-c,n-d}^n$ is Fano for $n \in \{c+d-1,c+d\}$ if $c \geq3$, and $n \in \{d+1,\dots,2d-1\}$ if $c=2$.
\end{proof}

We move on to the case when the pair of linear spaces do not span $\PP^n$.  

\begin{definition}Let $n > c+d-1$. For each $1 \leq i \leq c-1$ and a choice of flag $\{\Lambda_{n-c-d+i} \subseteq \Lambda_{n-i}\}$, let $D_{i}'$ denote the divisor class of the locus of subschemes $Z \in \HH_{c-1,d-1}^n$, for which the linear span of $\Lambda_{n-c-d+i} \cup (\Lambda_{n-i} \cap Z)$ has dimension less than $n-i$. Let $D_c'^{(1)}$ denote the class of the closure of the locus of subschemes supported on two distinct planes for which the $(d-1)$-plane meets a fixed $\Lambda_{n-d}$. Let $D_c'^{(2)}$ denote class of the closure of the locus of subschemes supported on two distinct planes for which the $(c-1)$-plane meets a fixed $\Lambda_{n-c}$. 

Let $F$ denote the class of the locus of subschemes $Z$ such that its linear span meets a fixed $\Lambda_{n-c-d}$.
\end{definition}

\begin{definition} Let $n > c+d-1$. For each $1 \leq i \leq c$, let $N_i'$ denote the divisor class of the locus of subschemes with an embedded $(c-i)$-plane.
\end{definition}

By lifting the curves $C_i,B_i^{(1)},B_i^{(2)}$ to $\HH_{c-1,d-1}^n$ (c.f. Lemma \ref{interstellar}) and defining curves $Y_1^{(1)}, Y_1^{(2)}, Y_2$ analogous to Definition \ref{Ys}, we obtain the following proposition. Since computations of the intersection numbers are exactly the same as Lemma \ref{intersection5} and Proposition \ref{lincomb2}, we omit the proof.

\begin{prop} Let $c \geq 2$ and $n > c+d-1$. Then we have,
$$
\emph{Eff}(\HH_{c-1,d-1}^n) = \langle N_1',\dots,N_c',D_{c}'^{(1)},D_{c}'^{(2)},F \rangle \quad \text{and} \quad 
\emph{Nef}(\HH_{c-1,d-1}^n)  = \langle D_1',\dots,D_{c-1}',D_{c}'^{(1)},D_{c}'^{(2)},F \rangle.
$$
Moreover, if $c=2$ we have $N_1 = D_2'^{(1)} + D_2'^{(2)} -D'_{1} \, \, \text{and} \,  \, N_2 = 2D_1' - D_2'^{(1)} - D_2'^{(2)} - F$.
If $c \geq 3$ we have
\begin{enumerate}[itemsep=0.5ex]
\item $N_1' = D_c'^{(1)} + D_c'^{(2)} -D_{c-1}'$, 
\item $N_2' = 2D_{c-1}'-D_{c-2}' - D_c'^{(1)}-D_c'^{(2)}$,
\item $N_i' = 2D_{c-i+1}' - D_{c-i}' - D_{c-i+2}'$ \, for all $3 \leq i \leq c-1$,
\item $N_c' = 2D_1' - D_2'-F$.
\end{enumerate}
\end{prop}

Here is the analogue of Proposition \ref{fanothree}.

\begin{prop} Let $c \geq 2$ and $n > c+d-1$. The component $\HH_{c-1,d-1}^n$ is Fano.
\end{prop}
\begin{proof} Similar to Proposition \ref{hilbertchow}, there is a morphism, $\Psi:\HH_{c-1,d-1}^n \to \HH_{n-c,n-d}^n$ with exceptional locus $N_c'$. As explained in the Proposition  \ref{fanothree} we deduce 
\begin{align*}
K_{\HH_{c-1,d-1}^n} &= \Psi^{\star}K_{\HH_{n-c,n-d}^n}+(n-c-d+1)N_c'  \\
				&=  \begin{cases} \Psi^{\star}K_{\HH_{n-2,n-d}^n} +(n-d-1)(2D'_1 - D_2'^{(1)} - D_2'^{(2)}-F) & \text{ if } c = 2 \\ 
							\Psi^{\star}K_{\HH_{n-c,n-d}^n} + (n-c-d+1)(2D_1'-D_2'-F) &\text{ if } c \geq 3. \end{cases}
\end{align*}
Using the expression for $K_{\HH_{n-c,n-d}^n}$ in Proposition \ref{cones2}, it follows that $-K_{\HH_{c-1,d-1}^n}$ is ample.
\end{proof}

Thus we deduce the main theorem of this section,

\begin{thm} \label{fanofour} The components $\HH_{c-1,d-1}^n$ and $\HH_{n-c,n-d}^n$ are Mori dream spaces. 
\end{thm}

\section*{Acknowledgement} We would like to thank Dawei Chen, Michael Christianson, David Eisenbud and Frank-Olaf Schreyer for helpful discussions. We are especially thankful to Michael Christianson for reading a draft and suggesting improvements. The author is partially supported by an NSERC PGSD scholarship.


\end{document}